\documentclass[11pt]{article}
\usepackage[utf8]{inputenc}
\usepackage{geometry}
\usepackage{etoolbox}
\usepackage{lmodern}
\usepackage[affil-it]{authblk} 
\usepackage{wrapfig}
\usepackage{graphicx}
\usepackage{listings}
\usepackage[dvipsnames]{xcolor}
\usepackage{tcolorbox}
\usepackage{subfiles}
\usepackage{esvect}
\usepackage{amsrefs}
\usepackage{url}

\usepackage{systeme}
\usepackage{enumerate}
\usepackage{amssymb, amsmath, amsthm, bm}
\usepackage{mathtools}
\usepackage{commath}
\usepackage[export]{adjustbox}
\usepackage{caption}
\usepackage{subcaption}
\usepackage[english]{babel}
\usepackage{dsfont}
\usepackage{hyperref}
\usepackage{cleveref}
\usepackage{thmtools}
\usepackage{thm-restate}

\newcommand{\comm}[1]{}  

\theoremstyle{plain}
\newtheorem{thm}{Theorem}[section]
\newtheorem{lemma}[thm]{Lemma}
\newtheorem{cor}[thm]{Corollary}
\newtheorem{prop}[thm]{Proposition}
\newtheorem{rmk}[thm]{Remark}
\newtheorem{eg}[thm]{Example}
\newtheorem{claim}[thm]{Claim}
\theoremstyle{definition}
\newtheorem{defn}[thm]{Definition}

\DeclareMathOperator{\trace}{tr}

\DeclareMathOperator{\polylog}{polylog}

\DeclareMathOperator{\dx}{dx}

\DeclareMathOperator{\dz}{dz}
\DeclareMathOperator{\dw}{dw}
\DeclareMathOperator{\val}{val}
\DeclareMathOperator{\var}{Var}

\def \a{\alpha}

\def \C{\mathcal{C}}

\def \Nb{\mathbb{N}}
\def \Eb{\mathbb{E}}
\def \Pbb{\mathbb{P}}

\def \({\left(}
\def \){\right)}
\def \[{\left[}
\def \]{\right]}
\def \<{\langle}
\def \>{\rangle}
\def \={\longleftrightarrow}

\def \mzshape{\alpha_{Z(m)}}

\numberwithin{equation}{section}
\numberwithin{figure}{section}

\title{The Spectrum of the Singular Values of Z-Shaped Graph Matrices}
\author{Wenjun Cai and Aaron Potechin}
\date{\today}

\begin{document}
\graphicspath{{images/}}

\maketitle

\begin{abstract}
Graph matrices are a type of matrix which has played a crucial role in analyzing the sum of squares hierarchy on average case problems. However, except for rough norm bounds, little is known about graph matrices. In this paper, we take a step towards better understanding graph matrices by determining the limiting distribution of the spectrum of the singular values of Z-shaped graph matrices. We then give a partial generalization of our results for $m$-layer Z-shaped graph matrices.
\end{abstract}
\newpage
\tableofcontents

\newpage
\section{Introduction}\label{section:intro}

\setlength{\parskip}{1.5mm}
\setlength{\baselineskip}{1.5em}
Graph matrices (see Definition \ref{defn:graph matrix}) are a type of matrix which is a powerful tool for analyzing problems on random inputs. Graph matrices were first developed in order to prove sum of squares lower bounds for planted clique \cites{meka2015sum, pmlr-v40-Deshpande15, 10.1145/3178538, BHKKMP16}. Since then, graph matrices have played a crucial role in further sum of squares lower bounds for average case problems \cites{hop17, MRX20, GJJPR20, Kunisky20, 9719766, potechin2022subexponential, 10.1145/3564246.3585221}. Recently, graph matrices have been used to analyze power-sum decompositions of polynomials \cite{9996747}, to analyze the ellipsoid fitting conjecture \cites{PTVW23, hsieh_et_al:LIPIcs.ICALP.2023.78}, and to analyze a class of first-order iterative algorithms including belief propagation and approximate message passing \cite{jones2024iterative}.

Currently, we only have a partial understanding of graph matrices. Ahn, Medarametla and Potechin \cites{MP16,AMP20} proved rough norm bounds on all graph matrices where the random input is $G(n,\frac{1}{2})$, generalized these norm bounds to apply to a much wider range of dense random inputs, and showed that several technical lemmas in the literature on analyzing random subspaces and analyzing tensor decompositions can be proved directly using graph matrices. Later on, rough norm bounds were shown for graph matrices where the random input is a sparse random graph \cites{9719766, rajendran2022concentration}. However, except for these rough norm bounds, little is known about graph matrices.

A natural question is whether we can analyze graph matrices more precisely. For example, can we make the norm bounds on graph matrices tight up to a factor of $(1 + o(1))$ rather than $\polylog(n)$? More ambitiously, can we determine the spectrum of the singular values and/or eigenvalues of graph matrices? In other words, can we find analogues of Wigner's semicircle law \cites{wigner1958distribution, wigner1993characteristic} and/or Girko's Circular Law \cite{girko} for graph matrices?

In this paper, we take a step towards this goal by analyzing the spectrum of the singular values of a family of graph matrices which we call multi-Z-shaped graph matrices. In particular, we determine the limiting distribution of the spectrum of singular values of the Z-shaped graph matrix (see \Cref{section:zshape-spectrum}) when the random input is $G(n,\frac{1}{2})$. We then describe how our techniques can be generalized to determine the limiting distribution of the spectrum of the singular values of multi-layer Z-shaped graph matrices.

In a follow-up paper \cite{CP22}, we show how to find the spectrum of the singular values of multi-layer Z-shaped graph matrices in the more general setting where the entries of the random input have an arbitrary probability distribution $\Omega$ with bounded moments where the odd moments are all $0$.\\

\subsection{Definitions}
In order to state our results, we need a few definitions.

\subsubsection{Graph matrices}
First, we need to define graph matrices.
\begin{defn}[Fourier characters]
Given a graph $G$ and a set $\displaystyle{E \subseteq { \binom{V(G)}{2}}}$ of possible edges, define the Fourier character $\chi_{E}(G)$ to be $\chi_{E}(G) = (-1)^{|E \setminus E(G)|} = \prod_{e \in E} e(G)$ where the \textit{edge variable} $e(G)$ is $e(G)=1$ if $e \in E(G)$ and $-1$ otherwise.

If $E$ is a multi-set then we define $\chi_{E}(G) = \prod_{e \in E}{e(G)}$. Note that $\chi_{E}(G) = \chi_{E_{red}}(G)$ where $\displaystyle E_{red} \subseteq \binom{V(G)}{2}$ is the set of edges which appear in $E$ with odd multiplicity.
\end{defn}

\begin{defn}[Shapes]
We define a \emph{shape} $\alpha$ to be a graph with vertices $V(\alpha)$, edges $E(\alpha)$, and distinguished tuples of vertices $U_{\alpha} = \(u_1,\ldots,u_{|U_{\alpha}|}\)$ and $V_{\alpha} = \(v_1,\ldots,v_{|V_{\alpha}|}\)$.
\end{defn}

\begin{defn}[Bipartite Shapes]
We say that a shape $\alpha$ is \emph{bipartite} if $U_{\alpha} \cap V_{\alpha} = \emptyset$, $V(\alpha) = U_{\alpha} \cup V_{\alpha}$, and all edges in $E(\alpha)$ are between $U_{\alpha}$ and $V_{\alpha}$.
\end{defn}

\begin{defn}[Graph Matrices]\label{defn:graph matrix}
Given a shape $\alpha$, we define the \emph{graph matrix} $M_{\alpha}$ (which depends on the input graph $G$) to be the $\frac{n!}{\(n-|U_{\alpha}|\)!}\times\frac{n!}{\(n-|V_{\alpha}|\)!}$ matrix with rows indexed by tuples of $|U_{\a}|$ distinct vertices, columns indexed by tuples of $|V_{\a}|$ distinct vertices, and entries
\begin{equation*}
    M_{\a}(A,B) = \sum_{\substack{\sigma: V(\alpha) \to V(G): \  \sigma \text{ is injective}, \\ \sigma(U_{\a}) = A, \sigma(V_{\a}) = B}}{\chi_{\sigma(E(\a))}(G)}.
\end{equation*}
\end{defn}
\begin{rmk}
In this paper, we only consider the special case where $\a$ is bipartite. When $\a$ is bipartite, if $A \cap B = \emptyset$ then there is unique injective map $\sigma: V(\a) \to V(G)$ such that $\sigma(U_{\alpha}) = A$ and  $\sigma(V_{\alpha}) = B$ so $M_{\alpha}(A,B) = \chi_{\sigma(E(\a))}(G)$. If $A \cap B \neq \emptyset$ then there are no injective maps $\sigma: V(\a) \to V(G)$ such that $\sigma(U_{\alpha}) = A$ and  $\sigma(V_{\alpha}) = B$ so $M_{\alpha}(A,B) = 0$.
\end{rmk}
\begin{rmk}
While we only consider the setting where the random input is $G(n,\frac{1}{2})$ in this paper, graph matrices can also be defined for other random inputs. For details, see \cite{AMP20}.
\end{rmk}

\subsubsection{Spectrum of the singular values of a matrix}
We now recall the singular value decomposition of a matrix and define the spectrum of the singular values of a matrix.
\begin{defn}[Singular value decomposition] Given a rank $r$ matrix $M$ with real entries, the singular value decomposition (SVD) of $M$ is $\displaystyle M = \sum_{i=1}^{r}{{\sigma_i}{u_i}{v_i}^T}$ where the vectors $u_1,\ldots,u_r$ are orthonormal, the vectors $v_1,\ldots,v_r$ are orthonormal, and $\sigma_1 \geq \ldots \geq \sigma_r > 0$ are the nonzero singular values of $M$.
\end{defn}
\begin{rmk}
While the singular value decomposition of $M$ may not be unique, the nonzero singular values $\sigma_1,\ldots,\sigma_r$ are unique.
\end{rmk}
In this paper, we analyze the spectrum of the singular values of $N \times N$ matrices (where $N$ is a function of $n$).
\begin{defn}
Given an $N \times N$ matrix $M$, we define the spectrum $D_M$ of the singular values of $M$ to be the uniform distribution on the singular values of $M$ (including singular values which are $0$). In other words, if the nonzero singular values of $M$ are $\sigma_1 \geq \sigma_2 \cdots \geq \sigma_r > 0$ then $D_M$ is the probability distribution where each positive $\sigma \in \mathbb{R}$ occurs with probability $\frac{1}{N}\abs{\{i \in [r]: \sigma_i = \sigma\}}$ and $0$ occurs with probability $\frac{N-r}{N}$.
\end{defn}
The key fact which we need about the singular values of $M$ is that we can compute the $2k$th moment of $D_M$ by computing $\trace\(\(M{M^T}\)^k\)$.
\begin{prop}
For all $k \in \mathbb{N}$, 
$\displaystyle \trace\(\(M{M^T}\)^k\) = \sum_{i=1}^{r}{\sigma_i^{2k}}$
\end{prop}

\subsubsection{Generalized Catalan numbers}
Finally, we need the generalized Catalan numbers (\cite{FussCatalan}, \cite{fuss1791solutio}).
\begin{defn}\label{defn:generalized-calatan-number}
The generalized Catalan number $C_{m,r}(k)$ is $\displaystyle C_{m,r}(k)=\frac{r}{mk+r}\binom{mk+r}{k}$.
\end{defn}
\begin{rmk}
Note that $\displaystyle C_{2,1}(k) = \frac{1}{2k+1}\binom{2k+1}{n} = \frac{1}{k+1}\binom{2k}{k}$ gives the usual Catalan numbers. 
\end{rmk}
The generalized Catalan number $\displaystyle C_{3,1}(k) = \frac{1}{3k+1}\binom{3k+1}{k} = \frac{1}{2k+1}\binom{3k}{k}$ plays a key role in our analysis, so we use the shorthand $C'_k$ for it.
\begin{defn}
We define $\displaystyle C'_k = C_{3,1}(k) = \frac{1}{2k+1}\binom{3k}{k}$.
\end{defn}

\subsection{Our results}
\begin{figure}[hbt!]
    \centering
    \includegraphics[scale=0.35]{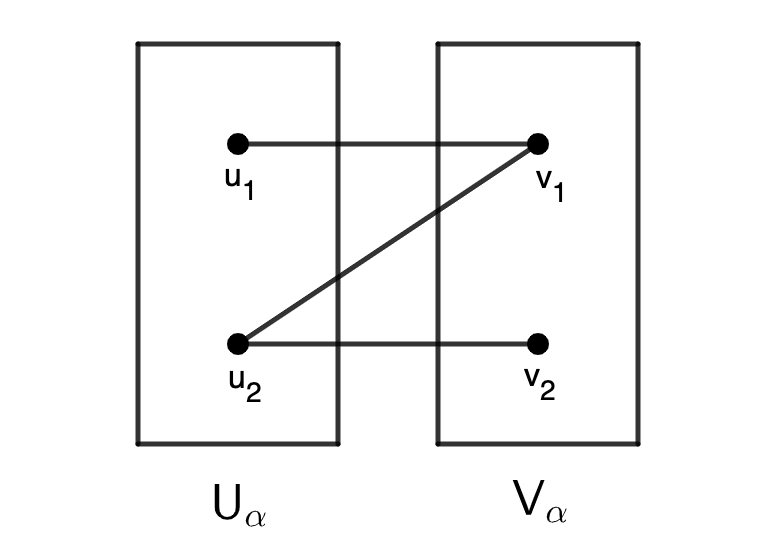}
    \caption{Z-shape $\a_Z$}
    \label{fig:graph H}
\end{figure}

\begin{figure}[hbt!]
      \centering
      \includegraphics[scale = 0.54]{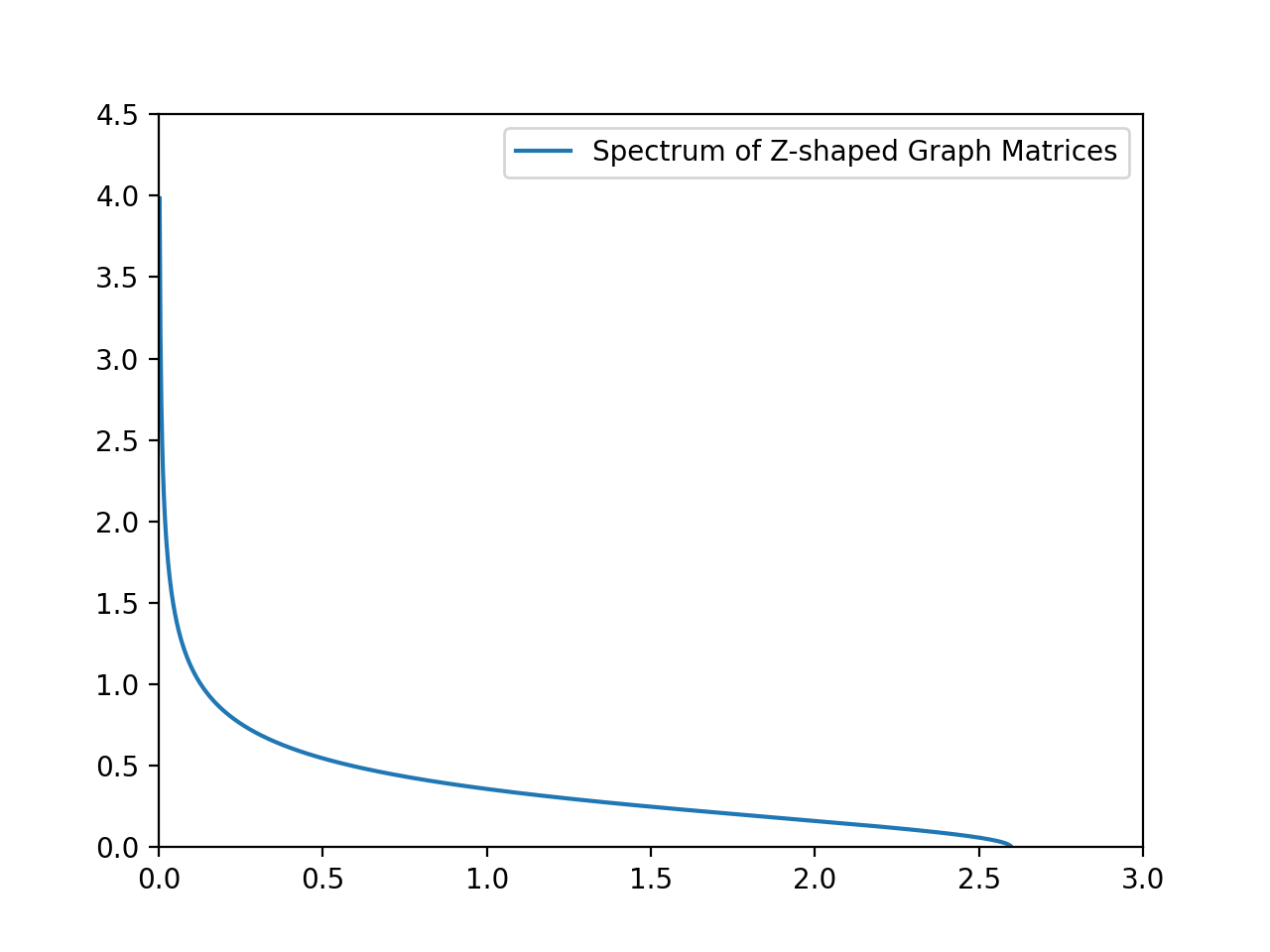}
      \caption{The limiting distribution of the singular values of $\frac{1}{n}M_{\alpha_Z}$ as $n \to \infty$}
      \label{introspectrumfigure}
    \end{figure}

Our main result is determining the limiting distribution of the spectrum of the singular values of Z-shaped graph matrices as $n \to \infty$.
\begin{defn}\label{defn:z-shape}
    Let $\a_{Z}$ be the bipartite shape with vertices $V(\a_{Z})=\{u_1,u_2,v_1,v_2\}$, edges $E(\a_{Z})=\left\{\{u_1,v_1\},\{u_2,v_1\},\{u_2,v_2\}\right\}$, and distinguished tuples of vertices $U_{\a_{Z}}=(u_1,u_2)$ and $V_{\a_{Z}}=(v_1,v_2)$. See Figure \ref{fig:graph H} for an illustration. We call $\a_{Z}$ the \emph{Z-shape}.
\end{defn} 
\begin{rmk}
The direct definition of the Z-shaped graph matrix $M_{\a_Z}$ is as follows. When $\{a_1,a_2\} \cap \{b_1,b_2\} = \emptyset$, $M_{\a_Z}\left((a_1,a_2),(b_1,b_2)\right) = e_1(G)e_2(G)e_3(G)$ where $e_1 = \{a_1,b_1\}$, $e_2 = \{a_2,b_1\}$, and $e_3 = \{a_2,b_2\}$. When $\{a_1,a_2\} \cap \{b_1,b_2\} \neq \emptyset$, $M_{\a_Z}\left((a_1,a_2),(b_1,b_2)\right) = 0$.
\end{rmk}
\begin{defn}
Let $a = \dfrac{3\sqrt{3}}{2}$ and define $g_{\alpha_Z}: \mathbb{R} \setminus \{0\} \to \mathbb{R}$ to be the function such that 
\begin{equation*} 
g_{\alpha_Z}(x) = \dfrac{i}{\pi}\cdot \(\sqrt{3}\cdot \sin\(\dfrac{1}{3}\cdot\arctan\(\dfrac{3}{\sqrt{4x^2/3-9}}\)\)+ \cos\(\dfrac{1}{3}\cdot\arctan\(\dfrac{3}{\sqrt{4x^2/3-9}}\)\)\)
\end{equation*}
if $x \in (0,a]$ and $g_{\alpha_Z}(x) = 0$ if $x < 0$ or $x > a$. See Figure \ref{introspectrumfigure} for an illustration.
\end{defn}
\begin{thm}
As $n \to \infty$, the spectrum of the singular values of $\frac{1}{n}M_{\alpha_Z}$ weakly converges to $g_{\alpha_Z}$.
\end{thm}
After proving this result, we apply our techniques to give a partial analysis of the spectrum of the singular values of multi-layer Z-shaped graph matrices.

\subsection{Proof overview and paper outline}
In order to obtain the limiting distribution of the singular values of the normalized Z-shape graph matrices $\dfrac{1}{n}M_{\a_Z}$ as $n\to\infty$, we use the trace power method. For this method, we observe (see Section \ref{section:prelim} for the precise details) that if $\{M_N\}$ is a sequence of random $N \times N$ matrices and the spectrum of the singular values of $M_N$ weakly converges to the distribution with probability density function $g(x)$ then $\displaystyle \int_{-\infty}^{\infty}{x^{2k}g(x) \,dx} = \lim_{N \to \infty}{\frac{1}{N}\Eb\Big[\trace\left(\left({M_N}{M_N}^T\right)^k\right)\Big]}$.

In particular, taking $M_N = \dfrac{M_{\a_Z}}{n}$, if the spectrum of the singular values of $M_N = \dfrac{M_{\a_Z}}{n}$ weakly converges to the distribution with probability density function $g_{\a_Z}(x)$ then we have that for all $k \in \mathbb{N}$,
\begin{equation}
\int_{-\infty}^{\infty}{x^{2k}g_{\a_Z}(x)dx} = \lim_{n \to \infty}{\frac{1}{n(n-1)}\Eb\Bigg[\trace\(\(\frac{M_{\a_Z}M_{\a_Z}^T}{n^2}\)^k\)\Bigg]}
\end{equation}
Thus, by computing the leading order term of $\Eb\Big[\trace\left(\left(M_{\a_Z}M_{\a_Z}^T\right)^k\right)\Big]$ for each $k \in \mathbb{N}$, we can compute the even moments of $g_{\a_Z}(x)$. Once we have the even moments of $g_{\a_Z}(x)$, we can use these moments to recover $g_{\a_Z}(x)$. In the next few sections, we implement this approach as follows.

In \Cref{section:prelim}, we reduce the problem of computing the leading order term of $\Eb\Big[\trace\(\(M_{\a_Z}M_{\a_Z}^T\)^k\)\Big]$ to a combinatorics problem, namely counting the number of dominant constraint graphs on $H(\a_Z,2k)$ (see \Cref{section:zshape}). 
    
In \Cref{section:zshape}, we show that the number of dominant constraint graphs on $H(\a_Z,2k)$ is $\displaystyle C'_k = \dfrac{1}{2k+1}\binom{3k}{k}$. This implies that $\displaystyle\int_{-\infty}^{\infty}{x^{2k}g_{\a_Z}(x)\, dx} = C'_k = \dfrac{1}{2k+1}\binom{3k}{k}$.
    
In \Cref{section:zshape-spectrum}, we use the recurrence relation $C'_k = \dfrac{3(3k-1)(3k-2)}{2k(2k+1)}C'_{k-1}$ for $C'_k$ to derive the differential equation $\displaystyle (4x^4-27x^2)g''_{\a_Z}(x)+(8x^3-27x)g'_{\a_Z}(x)+3g_{a_Z}(x)=0$ for $g_{\a_Z}$. Solving this differential equation gives $g_{\a_Z}$.

In \Cref{section:mzshape}, we generalize the results in \Cref{section:zshape} to calculate the number of dominant constraint graphs for multi-layer Z-shaped graph matrices. In \Cref{section:mzshape-spectrum}, we use this result to determine a differential equation for the limiting distribution of the singular values of one multi-layer Z-shaped graph matrix.

\section{Techniques}\label{section:prelim}

\setlength{\parskip}{1.5mm}
\setlength{\baselineskip}{1.5em}
In this section, we describe the techniques we will use for our analysis. We first describe more precisely what it means for the spectrum of the singular values of a random matrix to weakly converge to a distribution. We then describe how this convergence can be proved using the trace power method. Finally, we describe how the trace power method can be implemented using constraint graphs.

\subsection{Convergence of matrix spectra}
There are many equivalent definitions of weak convergence. Since we are considering probability distributions over $\mathbb{R}$, we can use the following definition.
\begin{defn}[Weak Convergence]
Let $\{D_n: n \in \mathbb{N}\}$ be a sequence of probability distributions on $\mathbb{R}$ and let $F_n(x) = \Pbb_{X \sim D_n}(X \leq x)$ be the cumulative distribution function of $D_n$. Let $D$ be another probability distribution on $\mathbb{R}$ and let $F(x) = \Pbb_{X \sim D}(X \leq x)$ be the cumulative distribution function of $D$. We say that the probability distributions $\{D_n: n \in \mathbb{N}\}$ weakly converge to the probability distribution $D$ if $\lim_{n \to \infty}{F_n(x)} = F(x)$ for all $x$ such that $F(x)$ is continuous.
\end{defn}

For our setting, the situation is more subtle because our probability distributions (the uniform distribution over the singular values of $\frac{M_{\a_Z}}{n}$) are themselves chosen from a probability distribution as $M_{\a_Z}$ is a function of a random graph $G \sim G(n,\frac{1}{2})$. We handle this by showing that if we take a sequence of matrices of the form $\frac{M_{\a_Z}}{n}$ then with probability $1$, the spectra of the singular values of these matrices weakly converge to the distribution whose probability density function is $g_{\alpha_Z}(x)$. Our precise theorem statement is as follows.
\begin{thm}\label{thm:precisemainresult}
If $\{G_n:n \in \mathbb{N}\}$ is a sequence of $G\(n,\frac{1}{2}\)$ graphs and $M_n = \frac{1}{n}{M_{\alpha_Z}(G_n)}$, then with probability $1$, the probability distributions $\{D_{M_n}: n \in \mathbb{N}\}$ weakly converge to $D_{\a_Z}$ where $D_{\a_Z}$ is the probability distribution whose probability density function is 
\begin{equation*} 
g_{\alpha_Z}(x) = \dfrac{i}{\pi}\cdot \(\sqrt{3}\cdot \sin\(\dfrac{1}{3}\cdot\arctan\(\dfrac{3}{\sqrt{4x^2/3-9}}\)\)+ \cos\(\dfrac{1}{3}\cdot\arctan\(\dfrac{3}{\sqrt{4x^2/3-9}}\)\)\)
\end{equation*}
if $x \in (0,\frac{3\sqrt{3}}{2}\,]$ and $g_{\alpha_Z}(x) = 0$ if $x < 0$ or $x > \frac{3\sqrt{3}}{2}$ where we take $\arctan(ix) = \dfrac{i}{2}\ln\(\dfrac{1+x}{1-x}\)$ for all real $x$, we take $ln(-x) = {\pi}i + ln(x)$ for all $x > 0$, and we take $sin(x) = \frac{e^{ix} - e^{-ix}}{2i}$ and $cos(x) = \frac{e^{ix} + e^{-ix}}{2}$ for all complex $x$.
\end{thm}

\subsection{Proving Theorem \ref{thm:precisemainresult} via the trace power method} \label{subsection:prelim-trace-power}
We now show how to prove Theorem \ref{thm:precisemainresult} using the trace power method. The key result which we need to prove Theorem \ref{thm:precisemainresult} is the following theorem for the trace power of $M_{\a_Z}$ which we prove in \Cref{section:zshape}.

\begin{thm}\label{thm:zshape-trace-convergence}
For all $k \in \mathbb{N} \cup \{0\}$,
\begin{enumerate}
    \item $\Eb\Big[\trace\(\(M_{\a_Z}M_{\a_Z}^T\)^k\)\Big] = C'_k\cdot n^{2k + 2} \pm O\(n^{2k+1}\)$.
    \item $\var \Big(\trace\(\(M_{\a_Z}M_{\a_Z}^T\)^k\)\Big)$ is $O\(n^{4k+2}\)$.
\end{enumerate}
where $\displaystyle C'_k = \dfrac{1}{2k+1}\binom{3k}{k}$ and the constants inside the big $O$ notation depend on $k$ but not on $n$.
\end{thm}

Using \Cref{thm:zshape-trace-convergence}, we can show that as $n \to \infty$, $D_{M_n}$ and $D_{\a_Z}$ have the same even moments. More precisely, we have the following lemma.
\begin{lemma}\label{lem:evenmomentsmatch}
With probability $1$, for all $k \in \mathbb{N} \cup \{0\}$, $\displaystyle \lim_{n\to\infty}{\Eb_{X \sim D_{M_n}} [X^{2k}]} = \Eb_{X \sim D_{\a_Z}} [X^{2k}] = C'_k$.
\end{lemma}
\begin{proof}
In \Cref{section:zshape-spectrum}, we verify directly that $\Eb_{x \sim D_{\a_Z}} [x^{2k}] = C'_k$ (see Theorem \ref{thm:solution to ODE for Z-shape}). To show that $\displaystyle \lim_{n\to\infty}{\Eb_{x \sim D_{M_n}} [x^{2k}]} = C'_k$, we use the following proposition.
\begin{prop}\label{prop:rv-as-convergence}
    Let $x_n$ be a sequence of random variables and let $c$ be a constant. If $\displaystyle \lim_{n\to \infty} \Eb[x_n] = c$ and $\displaystyle \sum_{n}\, \var(x_n) <\infty$, then with probability $1$, $\displaystyle \lim_{n\to\infty} x_n = c$.
\end{prop}
Given $k \in \mathbb{N}$, let $x_n = \Eb_{x \sim D_{M_n}} [x^{2k}] = \dfrac{1}{n(n-1)}\trace\(\(\frac{M_{\a_Z}M_{\a_Z}^T}{n^2}\)^k\)$ and let $c = C'_k$. By Theorem \ref{thm:zshape-trace-convergence}, $\displaystyle \lim_{n \to \infty}{x_n} = C'_k = c$ and $\displaystyle \sum_{n=1}^{\infty}{\var(x_n)}$ is $\displaystyle O\big(\sum_{n=1}^{\infty}{\frac{1}{n^2}}\big) < \infty$. By Proposition \ref{prop:rv-as-convergence}, with probability $1$, $\displaystyle \lim_{n \to \infty}{\Eb_{x \sim D_{M_n}}[x^{2k}]} = \lim_{n \to \infty}{x_n} = 
c = C'_k$.
\end{proof}

We now use \Cref{lem:evenmomentsmatch} to prove \Cref{thm:precisemainresult}.
\begin{proof}[Proof of Theorem \ref{thm:precisemainresult}]
Let $D'_n$ be the probability distribution obtained by taking a sample $x$ from $D_{M_n}$ and then taking $x$ with probability $\frac{1}{2}$ and $-x$ with probability $\frac{1}{2}$. In other words, $D'_n$ is the probability distribution whose even moments are the same as $D_{M_n}$ and whose odd moments are $0$. Similarly, let $D'$ be the probability distribution obtained by taking a sample $x$ from $D_{\a_Z}$ and then taking $x$ with probability $\frac{1}{2}$ and $-x$ with probability $\frac{1}{2}$. 

To show that $\{D_{M_n}: n \in \mathbb{N}\}$ weakly converges to $D_{\a_Z}$, it is sufficient to show that $\{D'_n: n \in \mathbb{N}\}$ weakly converges to $D'$. To show this, we use the following lemma which can be shown using Lemmas B.1 and B.3 in \cite{bai2010spectral}.
\begin{lemma}\label{lem:moment-convergence}
    Let $\{D_n: n \in \mathbb{N}\}$ be a sequence of probability distributions over $\mathbb{R}$ and let $D$ be another probability distribution over $\mathbb{R}$. $\{D_n: n \in \mathbb{N}\}$ weakly converges to $D$ if the following statements are true:
    \begin{enumerate}
        \item $D_n$ has finite moments of all orders for all $n \in \mathbb{N}$. In other words, for all $n \in \mathbb{N}$ and all $k \in \mathbb{N} \cup \{0\}$, $\Eb_{x \sim D_n}[x^k]$ is finite.
        \item For each $k \in \mathbb{N} \cup \{0\}$, $\displaystyle \beta_k := \lim_{n\to\infty}{\Eb_{x \sim D_n} [x^k]}$ exists and is finite. Moreover, for each $k \in \mathbb{N} \cup \{0\}$, $\Eb_{x \sim D} [x^k] = \beta_k$.
        \item (Carleman Condition) $\displaystyle \sum_{k=1}^{\infty}{ {\beta_{2k}}^{-1/2k}} = \infty$.
    \end{enumerate}
\end{lemma}
We now make the following observations
\begin{enumerate}
\item The condition that for all $n \in \mathbb{N}$ and all $k \in \mathbb{N} \cup \{0\}$, $\Eb_{x \sim D'_n}[x^k]$ is finite is trivial as all finite matrices have finite singular values. 
\item If $k \in \mathbb{N}$ is odd then for all $n \in \mathbb{N}$, $\Eb_{x \sim D_n} [x^k] = \Eb_{X \sim D} [x^k] = 0$.
\item By Lemma \ref{lem:evenmomentsmatch}, if $k \in \mathbb{N} \cup \{0\}$ is even then with probability $1$, $\displaystyle \lim_{n\to\infty}{\Eb_{x \sim D_n} [x^k]} = \Eb_{x \sim D} [x^k] = C'_{k/2}$.
\item Since $\displaystyle C'_k = \dfrac{1}{2k+1}\binom{3k}{k} \leq 2^{3k}$, $\displaystyle \sum_{k=0}^{\infty}{\(C'_k\)^{-\frac{1}{2k}}} \geq \sum_{k=0}^{\infty}{\frac{1}{2^{3/2}}} = \infty$.
\end{enumerate}
Thus, with probability $1$, $\{D'_n: n \in \mathbb{N}\}$ weakly converges to $D'$ which implies that with probability $1$, $\{D_{M_n}: n \in \mathbb{N}\}$ weakly converges to $D_{\a_Z}$, as needed.
\end{proof}


\subsection{Analyzing \texorpdfstring{$\trace\(\(M_{\a_Z}M_{\a_Z}^T\)^k\)$}{tr((Ma MaT)k)} using constraint graphs} \label{subsection:prelim-constraint-graph}
To analyze $\trace\(\(M_{\a_Z}M_{\a_Z}^T\)^k\)$, we expand out $\trace\(\(M_{\a_Z}M_{\a_Z}^T\)^k\)$ as a sum of many different terms and use constraint graphs to visualize these terms. Many of the following definitions and results are taken from Section 3 of \cite{AMP20}.

\subsubsection{Expressing \texorpdfstring{$\trace\(\(M_{\a_Z}M_{\a_Z}^T\)^k\)$}{tr((Ma MaT)k)} using the graph \texorpdfstring{$H(\a_Z,2k)}{H(aZ,2k)}$}

We first show how to express $\trace\(\(M_{\a_Z}M_{\a_Z}^T\)^k\)$ in terms of maps $\phi$ from the vertices of a graph $H(\a_Z,2k)$ to $[n]$. We will then use constraint graphs to visualize these maps.

\begin{figure}[hbt!]
	    \centering
	    \includegraphics[scale=0.32]{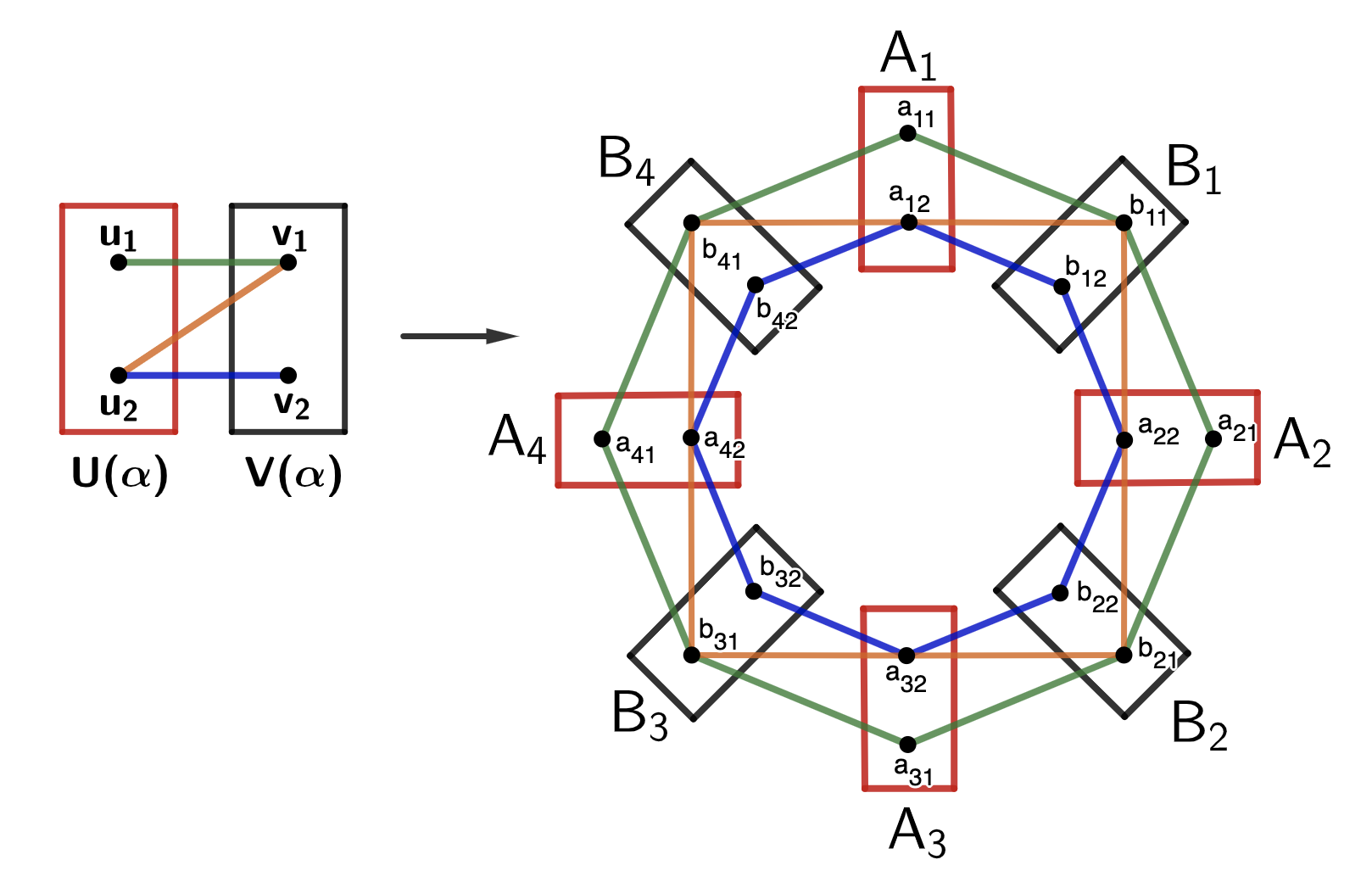}
	    \caption{Illustration of \Cref{def:copies}: $\a$ and $H(\a,2k)$ where $k=4$.}
	    \label{fig:H-a-2k}
	\end{figure}

\begin{defn}[Definition 3.2 of \cite{AMP20}]\label{def:copies}
    Given a shape $\a$ and a $k \in \mathbb{N}$, we define $H(\a,2k)$ to be the multi-graph which is formed as follows:
    \begin{enumerate}
		\item We take $k$ copies $\a_1,\ldots,\a_{k}$ of $\a$ and $k$ copies $\a^{T}_1,\ldots,\a^{T}_k$ of $\a^{T}$, where $\a^{T}$ is the shape obtained from $\a$ by switching the role of $U_{\a}$ and $V_{\a}$.
		\item We compose the shapes $\a_1,\a^{T}_1,\ldots,\a_k,\a^{T}_k$ by setting $V_{\a_{i}} = U_{\a^{T}_i}$ for all $i \in [k]$ and setting $V_{\a^{T}_i} = U_{\a_{i+1}}$ for all $i \in [k-1]$. We then set $V_{\a^{T}_k} = U_{\a_{1}}$.
    \end{enumerate}
    We define $V(\alpha,2k) = V(H(\alpha,2k))$ and we define $E(\alpha,2k) = E(H(\alpha,2k))$. See \Cref{fig:H-a-2k} for an illustration.
\end{defn}
\begin{rmk}
    $H(\a,2k)$ is defined as a multi-graph because edges will be duplicated if $U_{\a}$ or $V_{\a}$ contains one or more edges or $k = 1$ and $\abs{E\(H(\a,2k)\)} > 0$. If $U_{\a}$ and $V_{\a}$ do not contain any edges and $k > 1$ then $H(\a,2k)$ will be a graph.
\end{rmk}
	
	\begin{defn}[Definition 3.4 of \cite{AMP20}: Piecewise injectivity]
		We say that a map $\phi:V\(\alpha,2k\) \to [n]$ is \emph{piecewise injective} if $\phi$ is injective on each piece $V(\alpha_i)$ and each piece $V\(\alpha^{T}_i\)$ for all $i \in [k]$. In other words, $\phi(u) \neq \phi(v)$ whenever $u,v \in V\(\alpha_i\)$ for some $i \in [k]$ or $u,v \in V\(\alpha^T_i\)$ for some $i \in [k]$.
	\end{defn}
	
As observed in \cite{AMP20}, with these definitions $\Eb\left[\trace\left((M_{\alpha}M_{\alpha}^T)^k\right)\right]$ can be reexpressed as follows
\begin{prop}[Proposition 3.5 of \cite{AMP20}]
		For all shapes $\alpha$ and all $k \in \mathbb{N}$,
		\begin{equation*}
		\Eb\left[\trace\left((M_{\alpha}M_{\alpha}^T)^k\right)\right] = \sum_{\substack{\phi: V(\alpha,2k) \to [n]:
		\\ \phi \text{ is piecewise injective}}}{\Eb\[\chi_{\phi\(E(\alpha,2k)\)}(G)\]}\,.
		\end{equation*}
\end{prop}

\subsubsection{Constraint graphs}
To analyze this expression, we use constraint graphs $C$ where an edge $\{u,v\} \in E(C)$ represents the fact that $\phi(u) = \phi(v)$.

\begin{defn}
We define a constraint graph $C$ on a set of vertices $V$ to be an acyclic graph with vertices $V(C) = V$. 
\end{defn}

\begin{defn}
Given a set of vertices $V$ and a constraint graph $C$ on $V$, we say that two vertices $u,v \in V$ are constrained together by $C$ if $u = v$ or there is a path from $u$ to $v$ in $C$. We denote this by $u \=_{C} v$. When $C$ is clear from context, we just write $u \= v$.

Similarly, we say that two edges $e = \{u,v\}$ and $e' = \{u',v'\}$ are constrained together by $C$ if either $u \=_{C} u'$ and $v \=_{C} v'$ or $u \=_{C} v'$ and $v \=_{C} u'$. We denote this by $e \=_{C} e'$. When $C$ is clear from context, we just write $e \= e'$.
\end{defn}

\begin{defn}
We say that two constraint graphs $C,C'$ are \emph{equivalent} if for any pair of vertices $u,v \in V$, $u \=_C v$ if and only if $u \=_{C'} v$.
We denote this by $C \equiv C'$.
\end{defn}

\begin{defn}[Induced constraint graphs] \label{defn:induced-constraint-graphs}
Given a multi-graph $H$, a constraint graph $C$ on $V(H)$, and a set of vertices $V \subseteq V(H)$, we define the \emph{induced constraint graph $C'$ on $V$} to be the constraint graph such that $V(C') = V$ and for all $u,v \in V$, $u \= v$ in $C'$ if and only if $u \= v$ in $C$.
\end{defn}

Each map $\phi: V(\alpha,2k) \to [n]$ has an associated constraint graph $C_{\phi}$.
\begin{defn}
Given a map $\phi: V(\alpha,2k) \to [n]$, we define $C_{\phi}$ to be the constraint graph such that $u \=_{C_{\phi}} v$ if and only if $\phi(u) = \phi(v)$. Equivalently, $C_{\phi}$ consists of a disjoint union of trees where $u$ and $v$ are in the same tree if and only if $\phi(u) = \phi(v)$ (note that all such constraint graphs are equivalent).
\end{defn}

Note that if $\phi,\phi'$ are two maps with the same constraint graph then $\phi$ and $\phi'$ are the same up to permuting the indices in $[n]$.

\begin{defn}[Definition 3.8 of \cite{AMP20}: Constraint graphs on $H(\alpha,2k)$]\label{defn:constraint-graph-on-H}
    We define $\mathcal{C}_{(\alpha,2k)} = \left\{C_{\phi}: \text{The map } \phi: V(\alpha,2k) \to [n] \text{ is piecewise injective}\right\}$ to be the set of all possible constraint graphs on $V(\alpha,2k)$ which come from a piecewise injective map $\phi: V(\alpha,2k) \to [n]$.
		
    Given a constraint graph $C \in \mathcal{C}_{(\alpha,2k)}$, we make the following definitions:
    \begin{enumerate}
	\item We define $N(C) = \abs{\left\{\phi: V(\alpha,2k) \to [n]: \phi \text{ is piecewise injective}, C_{\phi} \equiv C\right\}}$.
	\item We define $\val(C) = \Eb\big[\chi_{\phi\(E(\alpha,2k)\)}(G)\big]$ where $\phi: V(\alpha,2k) \to [n]$ is any piecewise injective map such that $C_{\phi} \equiv C$.
    \end{enumerate}
    We say that a constraint graph $C \in \mathcal{C}_{(\alpha,2k)}$ is \textit{nonzero-valued} if $\val(C)\neq 0$.
\end{defn}

We now define a simple shape which we call the \emph{line shape} and use it to give examples for \Cref{defn:constraint-graph-on-H}.
\begin{defn}\label{defn:line-shape}
    Let $\a_{0}$ be the bipartite shape with vertices $V(\a_{0})=\{u,v\}$ and a single edge $\{u,v\}$ with distinguished tuples of vertices $U_{\a_{0}}=(u)$ and $V_{\a_{0}}=(v)$. We call $\a_0$ the \textit{line shape}.
    
    \begin{figure}[hbt!]
        \centering
        \includegraphics[scale=0.32]{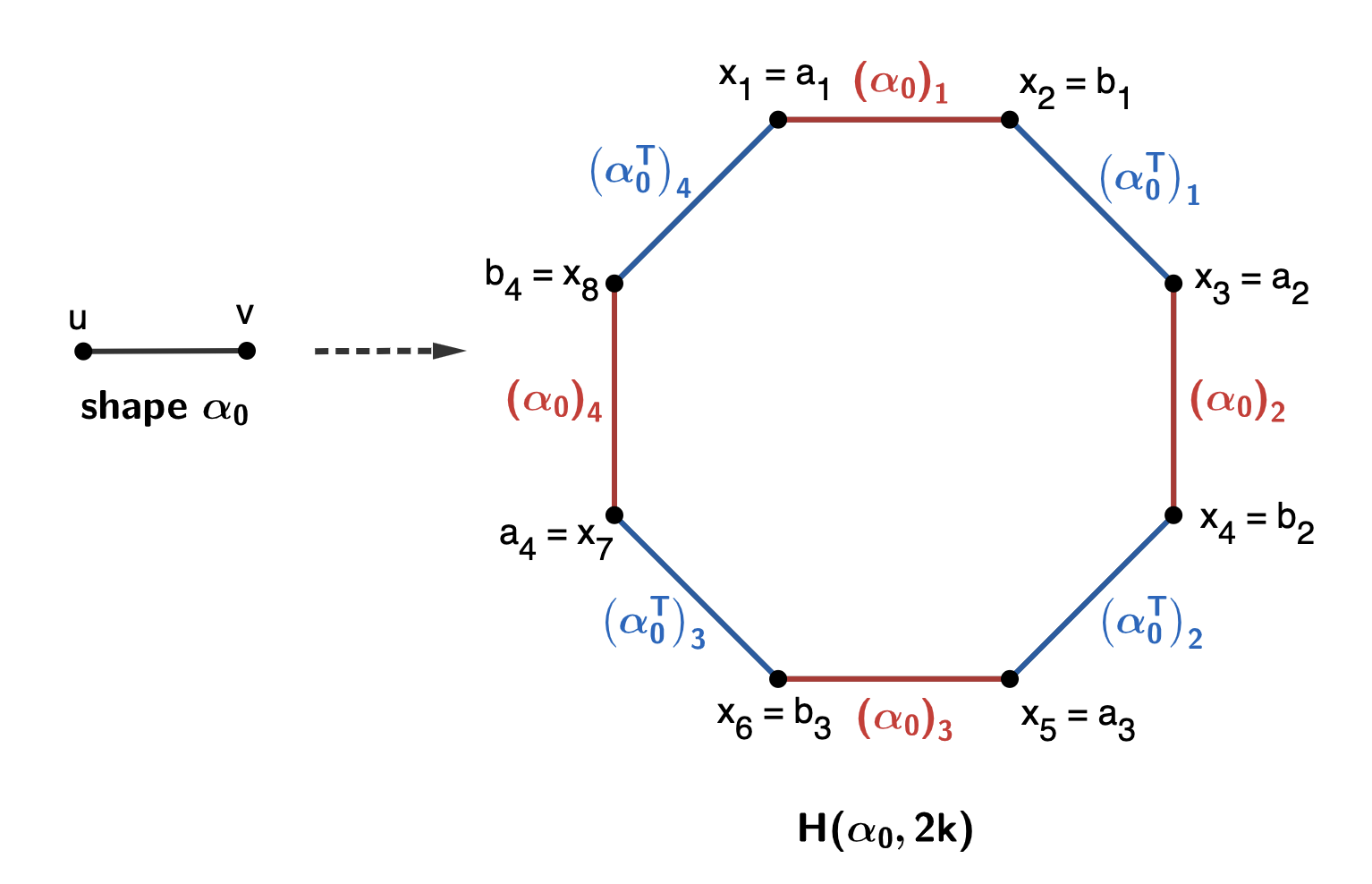}
        \caption{Illustration of \Cref{defn:line-shape} and \Cref{defn:line-shape-H-label}: $\a_0$ is the line shape. $H(\a_0,2k)$ is a cycle of length $2k$. Here $k=4$.}
        \label{fig:line-shape}
    \end{figure}
\end{defn}

\begin{defn}[Labelings of $H(\a_0,2k)$]\label{defn:line-shape-H-label}
Let $H(\a_0,2k)$ be the multi-graph given by \Cref{def:copies}. We label the vertices of $H(\a_0,2k)$ in two different but closely related ways.
\begin{enumerate}
\item Letting $(\a_0)_1,\ldots,(\a_0)_k$ be the copies of $\a_0$ and letting $\(\a_0^T\)_1,\ldots,\(\a_0^T\)_k$ be the copies of $\a_0^T$, we label the vertices of $V\((\a_0)_i\)$ as $\{u_i,v_i\}$ and the vertices of $V\big(\(\a_0^T\)_i\big)$ as $\{a_i,b_{i+1}\}$ where we take $a_{k+1} = a_1$. This gives the labeling $a_1,b_1,a_2,b_2,\ldots,a_k,b_k$ for the vertices of $H(\a_0,2k)$.
\item For convenience, we also label the vertices of $H(\a_0,2k)$ as $x_1,\ldots,x_{2k}$. In other words, for all $j \in [k]$, we take $x_{2j-1} = a_j$ and $x_{2j} = b_j$.
\end{enumerate}
\end{defn}

\begin{eg}\label{eg:val-C}
\begin{figure}[hbt!]
    \centering
    \includegraphics[scale=0.32]{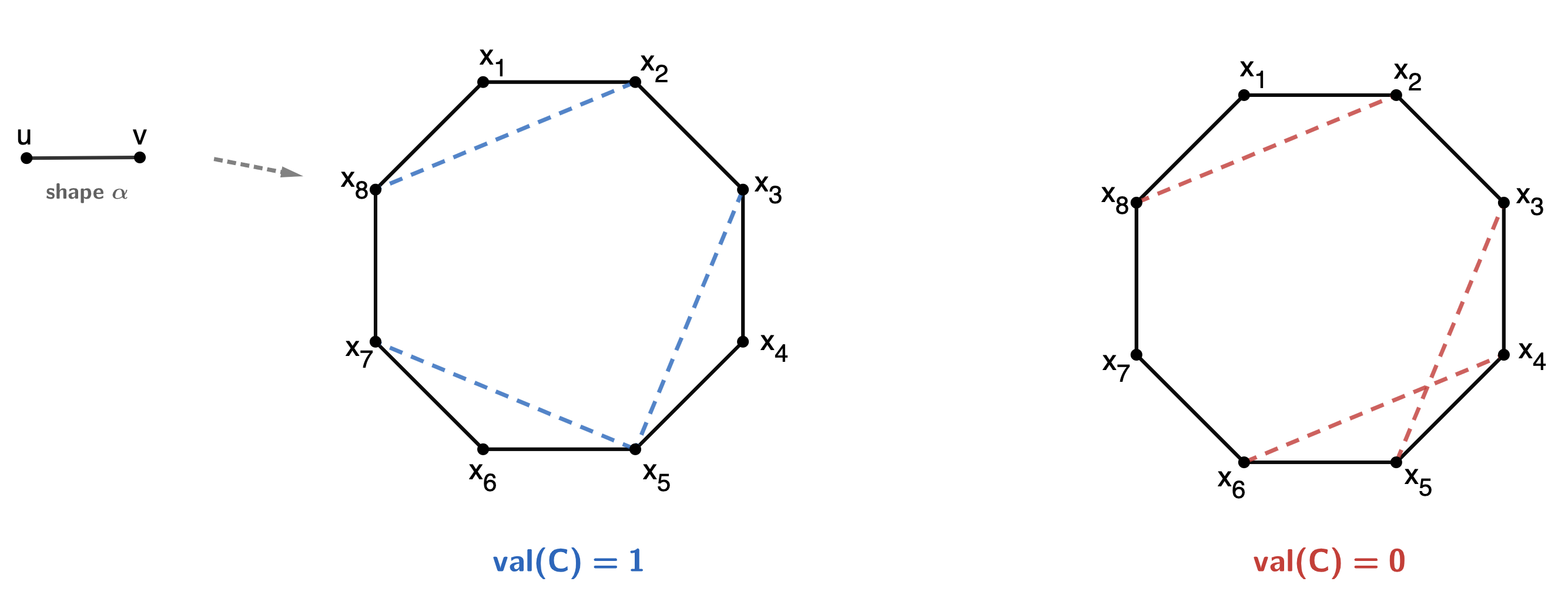}
    \caption{Illustration of \Cref{eg:val-C}.}
    \label{fig:example-valC}
\end{figure}
    Let $\a$ be the \emph{line shape} with $U(\a) = \{u\}$, $V(\a) = \{v\}$ and $E(\a) = \{\{u,v\}\}$ and consider its corresponding graph $H(\a,2k)$. Let $k=4$. 
    \begin{enumerate}
        \item If $E(C) = \left\{\{x_2,x_8\}, \{x_3,x_5\}, \{x_5, x_7\}\right\}$ then $\val(C) = 1$ and $N(C) = \dfrac{n!}{(n-5)!}$ since there are $8-\abs{E(C)} = 8-3 = 5$ distinct vertices whose values $\phi(v)$ can be chosen from $[n]$.
        \item If $E(C) = \left\{\{x_2,x_8\}, \{x_3,x_5\}, \{x_4,x_6\}\right\}$ then $\val(C) = 0$ and $N(C) = \dfrac{n!}{(n-5)!}$.
    \end{enumerate}
\end{eg}

As observed in \cite{AMP20}, with these definitions $\Eb\Big[\trace\left(\(M_{\alpha}M_{\alpha}^T\)^k\right)\Big]$can be re-expressed as follows.
\begin{prop}[Proposition 3.9 of \cite{AMP20}]\label{prop:exp value diffrep}
    For all shapes $\alpha$ and all $k \in \mathbb{N}$,
    \begin{equation*}
        \Eb\left[\trace\left((M_{\alpha}M_{\alpha}^T)^k\right)\right] = \sum_{C \in \mathcal{C}_{(\alpha,2k)}}{N(C)\val(C)}.
    \end{equation*}
\end{prop}
Given a constraint graph $C \in \mathcal{C}_{(\alpha,2k)}$, it can be useful to consider the graph obtained by starting with $H(\a,2k)$ and contracting the edges of $C$.

\begin{defn}\label{defn:H/C}
    \begin{figure}[hbt!]
        \centering
        \includegraphics[scale=0.34]{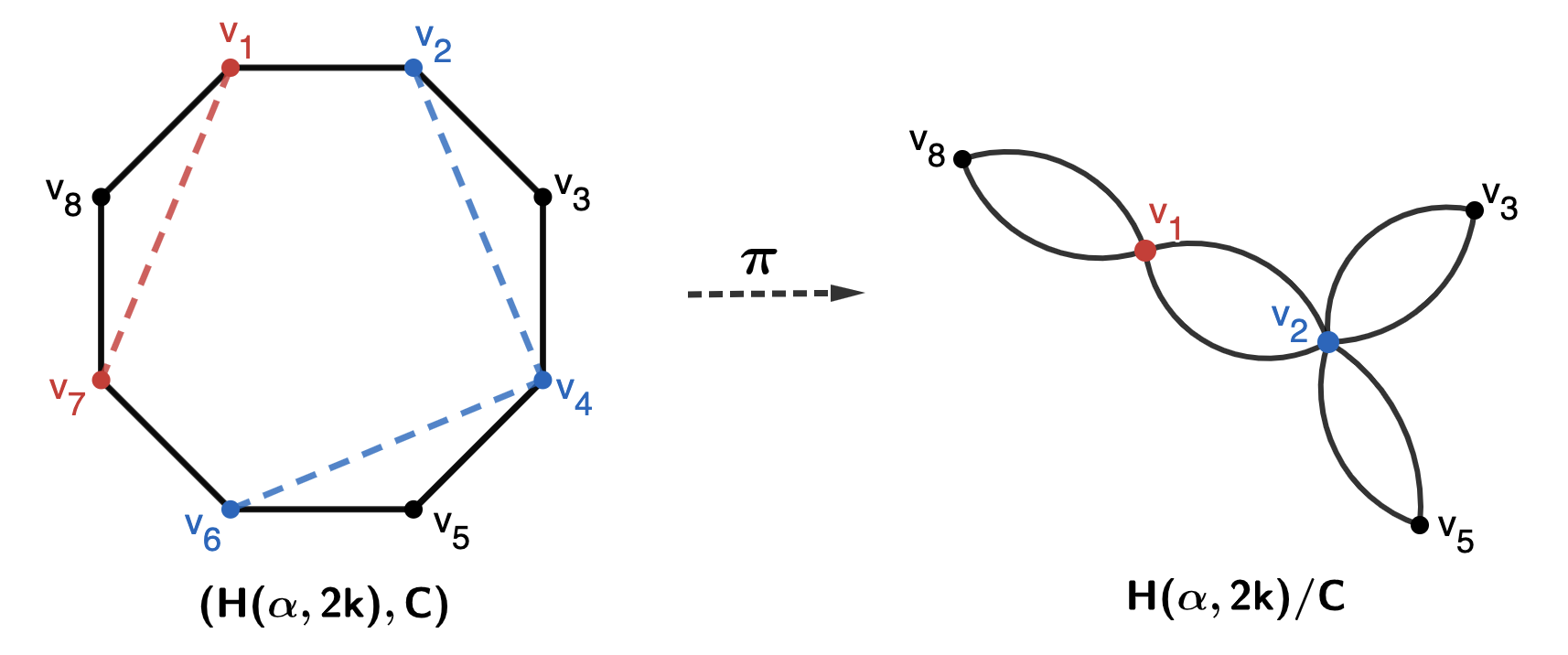}
        \caption{Illustration of \Cref{defn:H/C}: $H(\a,2k)/C$}
        \label{fig:H-constraint-graph-apply}
    \end{figure}
Given a multi-graph $H$ and a constraint graph $C$ on $V(H)$, we define $H/C$ to be the graph obtained from $H$ by contracting the edges in $C$. More precisely, we can define $H/C$ using the following map. Letting $V(H) = v_1,\ldots,v_m$ be the vertices of $H$, define 
$\pi: V(H) \to V(H)$ to be the map such that $\pi(v_i) = v_{\min\{j \in [m]:v_i \= v_j \text{ or } i = j\}}$. We define $H/C$ to be the multi-graph with vertices $V(H/C) = \{v_i: \pi(v_i) = v_i\}$ and edges $E(H/C) = \left\{\{\pi(v_i),\pi(v_j)\}: \{v_i,v_j\} \in E(G)\right\}$. See \Cref{fig:H-constraint-graph-apply} for an illustration.
\end{defn}

For technical reasons, the following definition is convenient.
\begin{defn}\label{def:appearingmultiplicity}
    Let $H$ be a multi-graph and let $C$ be a constraint graph on $V(H)$. Given two vertices $u,v \in V(H)$, we say that the edge $e = \{u,v\}$ appears with odd multiplicity in $H/C$ if $\{\pi(u),\pi(v)\}$ appears with odd multiplicity in $H/C$. Otherwise, we say that $e = \{u,v\}$ appears with even multiplicity in $H/C$.
\end{defn}
Note that \Cref{def:appearingmultiplicity} applies even if $\{u,v\} \notin E(H)$ or $u,v \notin V(H/C)$.

\subsubsection{Dominant constraint graphs}
We now analyze $\val(C)$ and $N(C)$ for constraint graphs $C \in \mathcal{C}_{(\a,2k)}$. We observe that the sum $\displaystyle \sum_{C \in \mathcal{C}_{(\alpha,2k)}}{N(C)\val(C)}$ is dominated by nonzero-valued constraint graphs with the minimum possible number of edges.
\begin{prop}[Proposition 3.10 of \cite{AMP20}]\label{prop:nonzero iff even multi-edges}
	For every constraint graph $C \in \mathcal{C}_{(\a,2k)}$, $\val(C) = 1$ if every edge in $\phi\(E(\a,2k)\)$ appears an even number of times and $\val(C) = 0$ otherwise (where $\phi: V(\a,2k) \to [n]$ is any piecewise injective map such that $C_{\phi} \equiv C$). Equivalently, $\val(C)=1$ if every edge in $H(\a,2k)/C$ appears an even number of times and $\val(C) = 0$ otherwise.
\end{prop}

\begin{prop}\label{prop: num of phi maps under a constr graph}
    For every constraint graph $C \in \mathcal{C}_{(\alpha,2k)}$, $\displaystyle N(C) = \frac{n!}{\(n - \abs{V(\alpha,2k)} + \abs{E(C)}\)!}$. 
\end{prop}
\begin{proof}
    Observe that choosing a piecewise injective map $\phi$ such that $C(\phi) \equiv C$ is equivalent to choosing a distinct element of $[n]$ for each of the $\abs{V(\alpha,2k)} - \abs{E(C)}$ connected components of $C$.
\end{proof}

Since the number of constraint graphs in $\mathcal{C}_{(\alpha,2k)}$ depends on $k$ but not on $n$, as $n \to \infty$ we only care about the nonzero-valued constraint graphs in $\mathcal{C}_{(\alpha,2k)}$ which have the minimum possible number of edges. We call such constraint graphs dominant constraint graphs.

\begin{defn}[Dominant Constraint Graphs]\label{defn:dominant constraint graph}
we say a constraint graph $C \in \mathcal{C}_{(\a,2k)}$ is a \textit{dominant constraint graph} if $\val(C) \neq 0$ and $
\abs{E(C)}=\min\left\{\abs{E(C')}: C' \in \mathcal{C}_{(\a,2k)}, \val(C') \neq 0\right\}$. More generally, given a multi-graph $H$ and a constraint graph $C$ on $V(H)$, we say that $C$ is a dominant constraint graph on $H$ if every edge of $H/C$ appears with even multiplicity and $C$ has the minimum possible number of edges (i.e., for all constraint graphs $C'$ on $V(H)$ such that $\abs{E(C')} < \abs{E(C)}$, 
\end{defn}

In order to determine the order of magnitude of $\displaystyle \Eb\Big[\trace\left(\(M_{\alpha}M_{\alpha}^T\)^k\right)\Big] = \sum_{C \in \mathcal{C}_{(\alpha,2k)}}{N(C)\val(C)}$, we need to know how many edges dominant constraint graphs in $\mathcal{C}_{(\a,2k)}$ have. It turns out that the number of edges in dominant constraint graphs in $\mathcal{C}_{(\a,2k)}$ only depends on $V(\a)$, $|U_{\a}|$, $|V_{\a}|$, $|U_{\a} \cap V_{\a}|$, and the minimum size of a vertex separator of $\a$.

\begin{defn}[Vertex Separators]
    We say that $S \subseteq V(\a)$ is a \emph{vertex separator} of $\a$ if every path from a vertex $u \in U_{\a}$ to a vertex $v \in V_{\a}$ contains at least one vertex in $S$. Note that this includes paths of length $0$ so we automatically have that $U_{\a} \cap V_{\a} \subseteq S$.
\end{defn}

\begin{defn}
	Given a shape $\alpha$, define $s_{\alpha}$ to be the minimum size of a vertex separator of $\alpha$.
\end{defn}

\begin{lemma}[Follows from Lemma 6.4 of \cite{AMP20}]\label{lem:min constr edges for bipartite shape}
    For any shape $\alpha$, for any nonzero-valued $C \in \mathcal{C}_{(\alpha,2k)}$, $\abs{E(C)} \geq k\abs{V\(\a\) \setminus \(U_{\a} \cup V_{\a}\)} + (k-1)(s_{\a} - \abs{U_{\a} \cap V_{\a}})$. Moreover, there exists a nonzero-valued $C \in \mathcal{C}_{(\alpha,2k)}$ such that $\abs{E(C)} = k\abs{V\(\a\) \setminus \(U_{\a} \cup V_{\a}\)} + (k-1)(s_{\a} - \abs{U_{\a} \cap V_{\a}})$ so this bound is tight.
\end{lemma}
\begin{rmk}
    In \cite{AMP20}, this result was only proved for well-behaved constraint graphs (see Definition \ref{defn:wellbehaved}). That said, using the ideas in Appendix B of \cite{AMP20}, it can be shown for all constraint graphs $C \in \mathcal{C}_{(\a,2k)}$.
\end{rmk}
\begin{cor}\label{cor:dominantnumberofedges}
For all bipartite shapes $\alpha$, for all dominant constraint graphs $C \in \mathcal{C}_{(\alpha,2k)}$, $\abs{E(C)} = (k-1)s_{\alpha}$.
\end{cor}
We can now reduce the problem of estimating $\Eb\Big[\trace\(\(M_{\a}M_{\a}^T\)^k\)\Big]$ to the problem of counting the number of dominant constraint graphs in $\mathcal{C}_{(\a,2k)}$.

\begin{cor}\label{cor:expected trace for bipartite shapes}
For all bipartite shapes $\alpha$ and all $k \in \mathbb{N}$,
    \begin{equation*}
        \Eb\left[\trace\(\(M_{\a}M_{\a}^T\)^k\)\right] = \abs{\left\{C \in \mathcal{C}_{(\a,2k)}: C \text{ is dominant}\right\}}n^{k\(|V(\a)| - s_{\a}\) + s_{\a}} \pm O\(n^{k\(|V(\a)| - s_{\a}\) + s_{\a} - 1}\).
    \end{equation*}
\end{cor}
\begin{proof}

By \Cref{prop:exp value diffrep}, $\displaystyle \Eb\Big[\trace\left(\(M_{\alpha}M_{\alpha}^T\)^k\right)\Big] = \sum_{C \in \mathcal{C}_{(\alpha,2k)}}{N(C)\val(C)}$. By Propositions \ref{prop:nonzero iff even multi-edges} and \ref{prop: num of phi maps under a constr graph}, for each dominant $C \in \mathcal{C}_{(\alpha,2k)}$, $\val(C) = 1$ and 
\begin{equation*}
N(C) = \frac{n!}{\(n - k\abs{V(\a)} + (k-1)s_{\alpha}\)!} = n^{k\(|V(\a)| - s_{\a}\) + s_{\a}} - O\(n^{k\(|V(\a)| - s_{\a}\) + s_{\a} - 1}\)
\end{equation*}
as $|V(\a,2k)| = k|V(\a)|$ and $\abs{E(C)} = (k-1)s_{\alpha}$. Thus, the contribution to $\Eb\Big[\trace\(\(M_{\a}M_{\a}^T\)^k\)\Big]$ from $C$ is $n^{k\(|V(\a)| - s_{\a}\) + s_{\a}} - O\(n^{k\(|V(\a)| - s_{\a}\) + s_{\a} - 1}\)$. For each non-dominant $C \in \mathcal{C}_{(\alpha,2k)}$, either $\val(C) = 0$ or $N(C)$ is $O\(n^{k\(|V(\a)| - s_{\a}\) + s_{\a} - 1}\)$ so the contribution to $\Eb\Big[\trace\(\(M_{\a}M_{\a}^T\)^k\)\Big]$ from $C$ is $O\(n^{k\(|V(\a)| - s_{\a}\) + s_{\a} - 1}\)$.
\end{proof}

Applying \Cref{cor:expected trace for bipartite shapes} to $\a_Z$ gives the following estimate for $\Eb\Big[\trace\(\(M_{\a_Z}M_{\a_Z}^T\)^k\)\Big]$.
\begin{cor}\label{cor:Zshapeexpectedtrace}
    $\Eb\Big[\trace\(\(M_{\a_Z}M_{\a_Z}^T\)^k\)\Big] = \abs{\left\{C \in \mathcal{C}_{(\a_Z,2k)}: C \text{ is dominant}\right\}}n^{2k + 2} \pm O\(n^{2k+1}\)$.
\end{cor}

\subsection{Partitioning \texorpdfstring{$H$}{H} via constraint edges}
In order to count the number of dominant constraint graphs for $H(\a,2k)$, it is useful to find sets of constraint edges which partition $H(\a,2k)$ into smaller graphs which can be analyzed recursively. In this subsection, we develop tools for doing this.

\begin{defn}
Given a multi-graph $H$, let $E_{odd}(H)$ be the set (not the multi-set) of edges of $H$ which have odd multiplicity and let $E_{even}(H)$ be the set of edges of $H$ which have even multiplicity.
\end{defn}
\begin{defn}\label{defn:Hpartitions}
Let $H$ be a multi-graph, let $E$ be a set of constraint edges which does not have any cycles, and let $C_E$ be the constraint graph on $V(H)$ with edges $E(C_E) = E$. We say that the set of constraint edges $E$ partitions $H$ into subgraphs $H_1,\ldots,H_j$ (which may actually be multi-graphs) if the following conditions hold:
\begin{enumerate}
\item $H_1,\ldots,H_j$ are subgraphs of $H/{C_E}$ and $E_{odd}(H_1),\ldots,E_{odd}(H_j)$ give a partition of $E_{odd}(H/{C_E})$. 

More precisely, letting $\pi: V(H) \to V(H)$ be the map describing $H/{C_E}$ (see \Cref{defn:H/C}) and letting $V(H/{C_E}) = \{v \in V(H): \pi(v) = v\}$ be the vertices of $H/{C_E}$, each subgraph $H_i$ is specified by a map $\pi_i: V(H) \to V(H/{C_E}) \cup \emptyset$ and these maps $\{\pi_i: i \in [j]\}$ satisfy the following properties:
\begin{enumerate}
\item[1.] For each vertex $u \in V(H/{C_E})$ and each $i \in [j]$, we have that either $\pi_i(u') = u$ for all vertices $u' \in V(H)$ such that $u' \=_{C_E} u$  or $\pi_i(u') = \emptyset$ for all vertices $u' \in V(H)$ such that $u' \=_{C_E} u$.
\item[2.] For each edge $e = \{u,v\} \in E(H)$ such that $e$ appears with odd multiplicity in $H/C_E$, there is exactly one $i \in [j]$ such that $\pi_i(u) \neq \emptyset$ and $\pi_i(v) \neq \emptyset$ (i.e., $\pi_i(e)$ exists). When this property is satisfied, we say that $e$ appears in $E(H_i)$.
\end{enumerate}
We then take $V(H_i) = \{v \in V(H/{C_E}): \pi_i(v) \neq \emptyset\}$ and $E(H_i) = \{\{\pi_i(u),\pi_i(v)\}: e = \{u,v\} \in E(H), \pi_i(u) \neq \emptyset, \pi_i(v) \neq \emptyset\}$ (note that $E(H_i)$ is a multi-set).
\item There is a bijection between dominant constraint graphs $C$ on $H$ such that $E \subseteq E(C)$ and dominant constraint graphs $C_1,\ldots,C_j$ on $H_1,\ldots,H_j$. 

More precisely, for all dominant constraint graphs $C_1,\ldots,C_j$ on $H_1,\ldots,H_j$, if we take $ E(C) = \(\bigcup_{i=1}^{j}{E(C_i)}\) \cup E$ then $E(C)$ does not contain any cycles, $\abs{E(C)} = |E| + \sum_{i=1}^{j}{\abs{E(C_i)}}$, and if we take $C$ to be the constraint graph on $H$ with constraint edges $E(C)$ then $C$ is a dominant constraint graph on $H$. Conversely, if $C$ is a dominant constraint graph on $H$ and $E \subseteq E(C)$ then if we take $C_i = \pi_i(C)$ to be the constraint graph on $H_i$ such that for all pairs of vertices $u,v \in V(H_i)$, $u \=_{C_i} v$ if and only if $u \=_{C} v$, we have that for all $i \in [j]$, $C_i$ is a dominant constraint graph on $H_i$.
\end{enumerate}
\end{defn}

\begin{eg}\label{eg:partition}
As observed in \Cref{defn:line-shape-H-label}, $H\(\a_0,2k\)$ is the graph with vertices $x_1,\ldots,x_{2k}$ and edges $\{\{x_i,x_{i+1}\}: i \in [2k-1]\} \cup \{x_{2k},x_1\}$.
\begin{enumerate}
\begin{figure}[hbt!]
    \centering
    \begin{subfigure}{.48\textwidth}
      \includegraphics[width=1\linewidth]{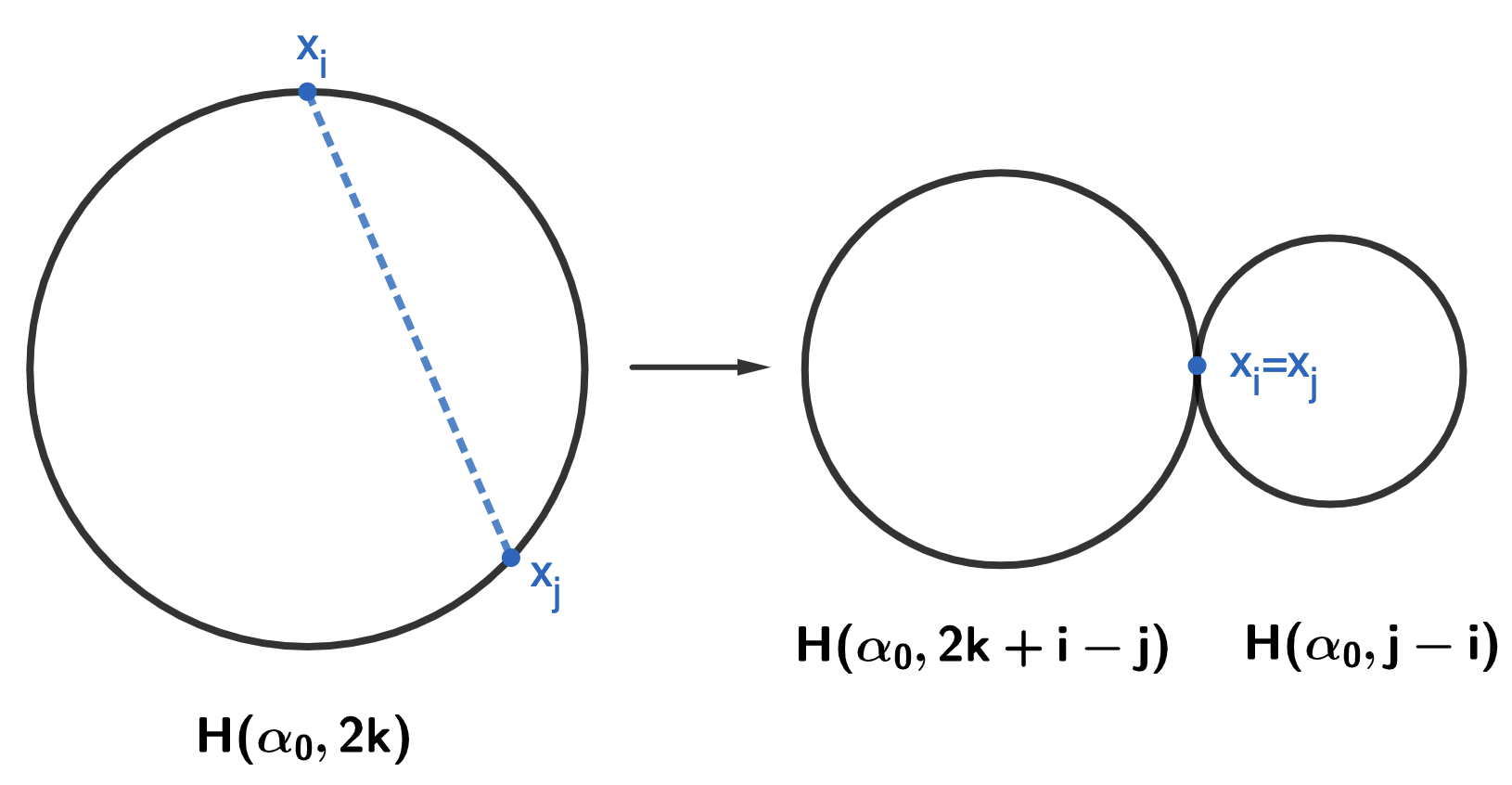}
      \caption{Example 1}
      \label{fig:Example1}
    \end{subfigure}\hspace{0.3cm}
    \begin{subfigure}{.49\textwidth}
      \includegraphics[width=1\linewidth]{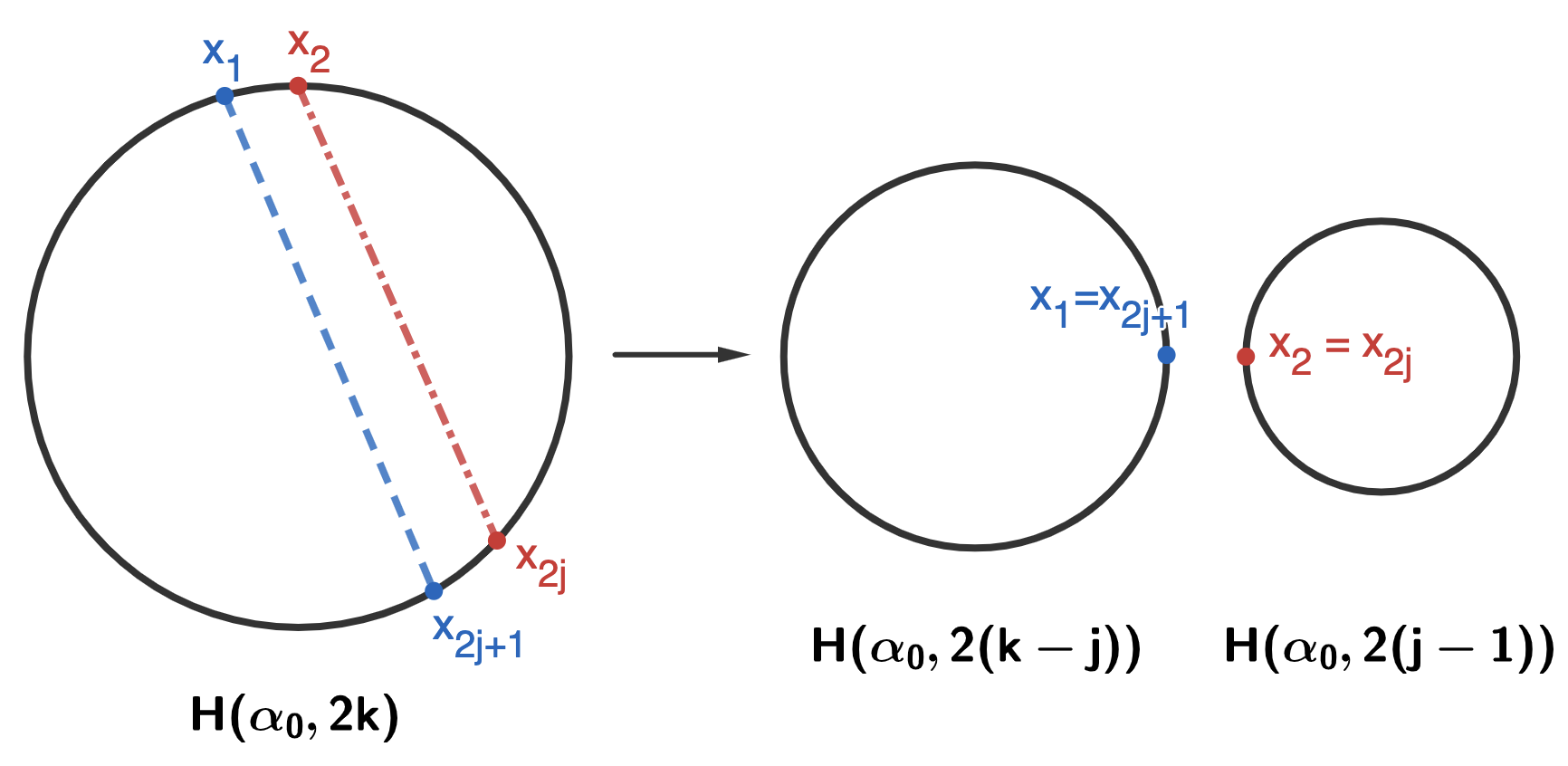}
      \caption{Example 2}
      \label{fig:Example2}
    \end{subfigure}
    \caption{Illustration of \Cref{eg:partition}}
    \label{fig:Example}
\end{figure}

    \item For all $i,j \in [2k]$ such that $i < j$ and $j-i$ is even, the constraint edge $\{x_i,x_j\}$ partitions $H\(\a_0,2k\)$ into $H_1 \simeq H\(\a_0,2k+i-j\)$ and $H_2 \simeq H\(\a_0,j-i\)$ where $V(H_1) = \{x_1,\ldots,x_{i}\} \cup \{x_{j+1},\ldots,x_{2k}\}$ and $V(H_2) = \{x_{i},x_{i+1},\ldots,x_{j-1}\}$.
    
    \item For all $j \in [2,k-1]$, the set of constraint edges $E = \{\{x_1,x_{2j+1}\}, \{x_2,x_{2j}\}\}$ partitions $H\(\a_0,2k\)$ into $H_1 \simeq H\(\a_0,2(j-1)\)$ and $H_2 \simeq H_{\a_0,2(k-j)}$ where $V(H_1) = \{x_2,\ldots,x_{2j-1}\}$ and $V(H_2) = \{x_{2j+2},x_{2j+3}\ldots,x_{2k}\} \cup \{x_1\}$.
\end{enumerate}
\end{eg}
In order to use Definition \ref{defn:Hpartitions}, it is helpful to have sufficient conditions for ensuring that the second condition of Definition \ref{defn:Hpartitions} holds, i.e., there is a bijection between dominant constraint graphs $C$ on $H$ such that $E \subseteq E(C)$ and dominant constraint graphs $C_1,\ldots,C_j$ on $H_1,\ldots,H_j$.
\begin{lemma}\label{lem:Hpartitionconditions}
Let $H$ be a multi-graph, let $E$ be a set of constraint edges, let $C_E$ be the constraint graph on $H$ with edges $E(C_E) = E$, and let $H_1,\ldots,H_j$ be subgraphs of $H/{C_E}$ satisfying the first condition of Definition \ref{defn:Hpartitions}.

If the following four conditions are satisfied then the second condition of Definition \ref{defn:Hpartitions} is satisfied as well.
\begin{enumerate}
\item Letting $H'_{final}$ be the graph with vertices $V(H/C_E)$ and edges $E_{odd}(H/C_E)$, $H'_{final}$ can be constructed by starting with $H' = H_i$ for some $i \in [j]$ and then iteratively applying the following operation. Take the current multi-graph $H'$ and choose $i,i'$ in $[j]$ such that 
\begin{enumerate}
\item[1.] $H_{i}$ has been added to $H'$ but $H_{i'}$ has not yet been added to $H'$.
\item[2.] $V(H') \cap V(H_{i'}) = V(H_i) \cap V(H_{i'})$.
\end{enumerate}
We then add the vertices and edges of $H_{i'}$ to $H'$.
\item For all $i,i' \in [j]$ and all dominant constraint graphs $C_i$ and $C_{i'}$ on $H_i$ and $H_{i'}$, for any pair of distinct vertices $u,v \in V(H_i) \cap V(H_{i'})$, there is no path from $u$ to $v$ in either $C_i$ or $C_{i'}$.
\item For all dominant constraint graphs $C$ on $H$ such that $E \subseteq E(C)$, for all $i,i' \in [j]$ and all $u,v \in V(H_i) \cap V(H_{i'})$, there is no path from $u$ to $v$ in $C$. 
\item For all dominant constraint graphs $C$ on $H$ such that $E \subseteq E(C)$, if $e_1,e_2 \in E(H)$ appear with odd multiplicity in $H/C_E$ and $e_1 \=_{C} e_2$ then $e_1$ and $e_2$ appear in the same subgraph $H_i$ (i.e, if we let $i,i' \in [j]$ be the unique indices such that $\pi_i(e_1)$ and $\pi_{i'}(e_2)$ exist then $i' = i$).
\end{enumerate}
\end{lemma}
\begin{rmk}
When we use \Cref{lem:Hpartitionconditions}, we will generally have that the first condition of \Cref{lem:Hpartitionconditions} is easy to check directly and the second and third conditions of \Cref{lem:Hpartitionconditions} follow from the fact that all dominant constraint graphs on $H/C_E$, $H_i$, and $H_{i'}$ are well-behaved (see \Cref{defn:wellbehaved}). Thus, we will just need to verify the fourth condition of \Cref{lem:Hpartitionconditions}.
\end{rmk}
\begin{proof}
To prove this lemma, we first show the following statement.
\begin{lemma}\label{lem:Hpartitionconditionssublemma}
Let $C_1,\ldots,C_j$ be constraint graphs on $H_1,\ldots,H_j$ such that for all $i,i' \in [j]$ and all pairs of distinct vertices $u,v \in V(H_i) \cap V(H_{i'})$, there is no path from $u$ to $v$ in either $C_i$ or $C_{i'}$. Letting $E(C') = \bigcup_{i \in [j]}{E(C_i)}$, $E(C')$ does not contain any cycles and $\abs{E(C')} = \sum_{i \in [j]}{\abs{E(C_i)}}$. Moreover, if we take $C'$ to be the constraint graph on $H/C_E$ with constraint edges $E(C')$, for all $i \in [j]$ and distinct vertices $u,v \in V(H_i)$, $C'$ contains a path between $u$ and $v$ if and only if $C_i$ contains a path between $u$ and $v$.
\end{lemma}
\begin{proof}
To prove this lemma, we use induction to show that at each step in the construction of $H/C_E$, letting $E(C') = \bigcup_{i: H_i \text{ has been added to } H'}{E(C_i)}$,
\begin{enumerate}
\item $\abs{E(C')} = \sum_{i: H_i \text{ has been added to } H'}{\abs{E(C_i)}}$ and $E(C')$ does not contain any cycles.
\item Taking $C'$ to be the constraint graph on $H'$ with constraint edges $E(C')$, for all $i$ such that $H_i$ has been added to $H'$, for all vertices $u,v \in V(H_i)$, $C'$ contains a path between $u$ and $v$ if and only if $C_i$ contains a path between $u$ and $v$.
\end{enumerate}

The base case $H' = H_i$ and $C' = C_i$ is trivial. For the inductive step, assume that the statement holds for the current $H'$. Let $i,i' \in [j]$ be the indices such that the next step of the construction attaches $H_{i'}$ to $H'$ at the vertices $V(H') \cap V(H_{i'}) = V(H_i) \cap V(H_{i'})$. Let $H''$ be the resulting graph after adding $H_{i'}$ to $H'$ and let $E(C'') = E(C') \cup E(H_{i'})$.

For all $u,v \in V(H') \cap V(H_{i'}) = V(H_i) \cap V(H_{i'})$, there is no path from $u$ to $v$ in $E(C_{i'})$ so we must have that $E(C') \cap E(C_{i'}) = \emptyset$ and thus $\abs{E(C'')} = \abs{E(C_{i'})} + \abs{E(C')}$. By the inductive hypothesis, $\abs{E(C')} = \sum_{i: H_i \text{has been added to } H'}{\abs{E(C_i)}}$ so
\begin{equation*}
\abs{E(C'')} = \abs{E(C_{i'})} + \sum_{i: H_i \text{ has been added to } H'}{\abs{E(C_i)}} = \sum_{i: H_i \text{ has been added to } H''}{\abs{E(C_i)}}.
\end{equation*}
To show that $E(C'')$ does not contain a cycle, observe that since $C'$ and $C_{i'}$ do not contain any cycles, any cycle in $E(C'')$ would have to contain a path in $C_{i'}$ between two distinct vertices $u,v \in V(H') \cap V(H_{i'}) = V(H_i) \cap V(H_{i'})$. However, $C_{i'}$ contains no such paths.

We now take $C''$ to be the constraint graph on $V(H'')$ with constraint edges $E(C'')$.

For all vertices $u,v \in V(H')$, $C''$ contains a path between $u$ and $v$ if and only if $C'$ contains a path from $u$ to $v$. To see this, assume that there are vertices $u,v \in V(H')$ such that $C''$ contains a path $P$ from $u$ to $v$ but $C'$ does not contain a path from $u$ to $v$. If so, then $P$ must contain a path in $C_{i'}$ between two distinct vertices in $V(H') \cap V(H_{i'}) = V(H_{i}) \cap V(H_{i'})$. This gives a contradiction as $C_{i'}$ contains no such paths. 

Similarly, for all vertices $u,v \in V(H_{i'})$, $C''$ contains a path between $u$ and $v$ if and only if $C_{i'}$ contains a path from $u$ to $v$. To see this, assume that there are vertices $u,v \in V(H_{i'})$ such that $C''$ contains a path $P$ from $u$ to $v$ but $C_{i'}$ does not contain a path from $u$ to $v$. If so, then $P$ must contain a path in $C'$ between two distinct vertices $u'$ and $v'$ in $V(H') \cap V(H_{i'}) = V(H_{i}) \cap V(H_{i'})$. This gives a contradiction as there is no path from $u'$ to $v'$ in $C_i$ and thus by the inductive hypothesis, there is no path from $u'$ to $v'$ in $C'$.
\end{proof}
\begin{cor}\label{cor:Hpartitionconditionssublemmacorollary}
    Let $C_1,\ldots,C_j$ be constraint graphs on $H_1,\ldots,H_j$ such that for all $i,i' \in [j]$ and all pairs of distinct vertices $u,v \in V(H_i) \cap V(H_{i'})$, there is no path from $u$ to $v$ in either $C_i$ or $C_{i'}$. Letting $E(C) = \left(\bigcup_{i \in [j]}{E(C_i)}\right) \cup E$, $E(C)$ does not contain any cycles and $\abs{E(C)} = |E| + \sum_{i \in [j]}{\abs{E(C_i)}}$. Moreover, if we take $C$ to be the constraint graph on $H$ with constraint edges $E(C)$, for all $i \in [j]$ and distinct vertices $u,v \in V(H_i)$, $C$ contains a path between $u$ and $v$ if and only if $C_i$ contains a path between $u$ and $v$.
\end{cor}
\begin{proof}
Let $C'$ be the constraint graph given by Lemma \ref{lem:Hpartitionconditionssublemma}. Observe that we can obtain $C$ from $C'$ by simply adding the vertices $\{v \in V(H): \pi(v) \neq v\}$ (i.e., the vertices of $H$ which are not in $H/C_E$) and the edges $E$. This does not add any cycles and does not affect whether there is a path between any two vertices $u,v \in V(H/C_E)$.
\end{proof}
We now use Lemma \ref{cor:Hpartitionconditionssublemmacorollary} to prove Lemma \ref{lem:Hpartitionconditions}.
\begin{lemma}\label{lem:Hpartitionconditionsdirectionone}
Under the conditions of Lemma \ref{lem:Hpartitionconditions}, if $C$ is a dominant constraint graph on $H$ such that $E \subseteq E(C)$ then if we take $C_i$ to be the constraint graph on $H_i$ such that for all vertices $u,v \in V(H_i)$, $u \=_{C_i} v$ if and only if $u \=_C v$, 
\begin{enumerate}
\item Taking $E(C') = \(\bigcup_{i=1}^{j}{E(C_i)}\) \cup E$, $E(C')$ has no cycles, $\abs{E(C')} = |E| + \sum_{i=1}^{j}{\abs{E(C_i)}}$, and if we take $C'$ to be the constraint graph on $H$ with constraint edges $E(C')$ then $C \equiv C'$.
\item For all $i \in [j]$, $C_i$ is a dominant constraint graph on $H_i$.
\end{enumerate}
\end{lemma}
\begin{proof}
We first observe that condition 4 of Lemma \ref{lem:Hpartitionconditions} implies that for all $i \in [j]$, $C_i$ is a nonzero-valued constraint graph on $H_i$. 

To see this, assume that there is an $i \in [j]$ and an edge $e = \{u,v\} \in E(H_i)$ which appears with odd multiplicity in $H_i/C_i$. Since $C$ is a nonzero-valued constraint graph on $H$, $e$ must appear with even multiplicity in $H/C$ which implies that there must be an edge $e' = \{u',v'\} \in E(H)$ such that $e' \=_{C} e$ and $e'$ appears with odd multiplicity in $H/C_E$ but $\pi(u') \notin V(H_i)$ or $\pi(v') \notin V(H_i)$. However, this contradicts condition 4 of \Cref{lem:Hpartitionconditions}.

By condition 3 of \Cref{lem:Hpartitionconditions}, for all $i,i' \in [j]$ and all $u,v \in V(H_i) \cap V(H_{i'})$, there is no path from $u$ to $v$ in $C$ so there is no path from $u$ to $v$ in either $C_i$ or $C_{i'}$. By Corollary \ref{cor:Hpartitionconditionssublemmacorollary}, $E(C') = \left(\bigcup_{i \in [j]}{E(C_i)}\right) \cup E$ has no cycles and $\abs{E(C')} = |E| + \sum_{i=1}^{j}{\abs{E(C_i)}}$. Since $C_i$ is a nonzero-valued constraint graph on $H_i$ for all $i \in [j]$, $C'$ is a nonzero-valued constraint graph on $H$.

Now observe that for any $e = \{u,v\} \in E(C')$, $u \=_{C} v$. This implies that for all $u,v \in V(H)$, if $u \=_{C'} v$ then $u \=_{C} v$. If there were $u,v \in V(H)$ such that $u \=_{C} v$ but $C'$ does not contain a path from $u$ to $v$ then we would have that $\abs{E(C')} < \abs{E(C)}$ which gives a contradiction as $C$ is a dominant constraint graph on $H$. Thus, we must have that $u \=_{C'} v$ if and only if $u \=_{C} v$ and thus $C' \equiv C$.

To show that $C_i$ is a dominant constraint graph on $H_i$ for all $i \in [j]$, assume there is an $i' \in [j]$ and a dominant constraint graph $C'_{i'}$ on $H_{i'}$ such that $\abs{E(C'_{i'})} < \abs{E(C_{i'})}$. Letting $ E(C'') = \left(\bigcup_{i \in [j] \setminus \{i'\}}{E(C_i)}\right) \bigcup E\(C'_{i'}\) \bigcup E$, by \Cref{cor:Hpartitionconditionssublemmacorollary}, $E(C'')$ does not contain any cycles and $ \abs{E(C'')} = \abs{E} + \abs{E(C'_{i'})} + \sum_{i \in [j] \setminus \{i\}}{\abs{E(C_i)}} < \abs{E(C)}$. Since $C_i$ is a nonzero-valued constraint graph on $H_i$ for all $i \in [j] \setminus \{i\}$ and $C'_{i'}$ is a nonzero-valued constraint graph on $H_{i'}$, $C''$ is a nonzero-valued constraint graph on $H$. However, this contradicts the assumption that $C \equiv C'$ is a dominant constraint graph on $H$.
\end{proof}
\begin{lemma}
Under the conditions of Lemma \ref{lem:Hpartitionconditions}, if $C_1,\ldots,C_j$ are dominant constraint graphs on $H_1,\ldots,H_j$ and we take $E(C) = \(\bigcup_{i=1}^{j}{E(C_i)}\) \cup E$ then 
\begin{enumerate}
\item $E(C)$ does not contain any cycles and $\displaystyle \abs{E(C)} = |E| + \sum_{i=1}^{j}{\abs{E(C_i)}}$.
\item The constraint graph $C$ with constraint edges $E(C)$ is a dominant constraint graph on $H$.
\end{enumerate}
\end{lemma}
\begin{proof}
The first statement follows immediately from Corollary \ref{cor:Hpartitionconditionssublemmacorollary}. To show that $C$ is a nonzero-valued constraint graph on $H$, observe that since $C_i$ is a nonzero-valued constraint graph on $H_i$ for all $i \in [j]$, all edges in $H/C = (H/C_E)/C'$ appear with even multiplicity which implies that $C$ is a nonzero-valued constraint graph on $H$.

To show that $C$ is a dominant constraint graph, assume that there is a dominant constraint graph $C'$ on $H$ such that $\abs{E(C')} < \abs{E(C)}$. By Lemma 
\ref{lem:Hpartitionconditionsdirectionone}, if we take $C'_i$ to be the constraint graph on $H_i$ such that for all vertices $u,v \in V(H_i)$, $u \=_{C'_i} v$ if and only if $u \=_C' v$ then $ \abs{E(C')} = \abs{E} + \sum_{i=1}^{j}\, {\abs{E(C'_i)}}$ and $C'_i$ is a dominant constraint graph on $H_i$ for all $i \in [j]$. Now observe that 
\begin{equation*}
|E| + \sum_{i=1}^{j}{\abs{E(C_i)}} = \abs{E(C)} < \abs{E(C')} = |E| + \sum_{i=1}^{j}{\abs{E(C'_i)}}
\end{equation*}
so there must be some $i \in [j]$ such that $\abs{E(C'_i)} < \abs{E(C_i)}$ which contradicts the assumption that $C_i$ is a dominant constraint graph on $H_i$.
\end{proof}
\end{proof}

\section{Counting dominant constraint graphs on \texorpdfstring{$H(\a_Z,2k)$}{H(aZ, 2k)}} \label{section:zshape}

\setlength{\parskip}{1.5mm}
\setlength{\baselineskip}{1.3em}
In this section, we show that the number of dominant constraint graphs on $H(\a,2k)$ is $C'_k$.

\begin{thm}\label{thm:num-zshape-constraint-graph}
For all $k \in \mathbb{N}$, $\displaystyle \abs{\left\{C \in \mathcal{C}_{(\a_Z,2k)}: C \text{ is dominant}\right\}} = C'_k = \frac{1}{2k+1}\binom{3k}{k}$. 
\end{thm}
Combined with \Cref{cor:Zshapeexpectedtrace}, this shows that $\Eb\Big[\trace\(\(M_{\a_Z}M_{\a_Z}^T\)^k\)\Big] = C'_k\cdot n^{2k + 2} \pm O\(n^{2k+1}\)$, which is the first part of \Cref{thm:zshape-trace-convergence}.

After proving \Cref{thm:num-zshape-constraint-graph}, we show that $\var \Big(\trace\(\(M_{\a_Z}M_{\a_Z}^T\)^k\)\Big)$ is $O\(n^{4k+2}\)$, which completes the proof of \Cref{thm:zshape-trace-convergence}.

\subsection{Recurrence relation for \texorpdfstring{$C'_k$}{Cn'}}

One of the key ingredients for proving \Cref{thm:num-zshape-constraint-graph} is the following recurrence relation on $C'_k$.
\begin{thm}\label{thm:catalan2-recurrence-relation}
For all $k \in \mathbb{N} \cup \{0\}$,
    \begin{equation}
        C'_{k+1} = \sum_{a,b,c\geq 0: a+b+c=k} {C'_a}{C'_b}{C'_c}
               =\sum_{a=0}^{k}{ {C'_a} \(\sum_{b=0}^{k-a} {C'_b}{C'_{k-a-b}}\)}.
    \end{equation}
\end{thm}

To prove this recurrence relation, we consider walks on grids. This proof is a generalization of the third proof in the Wikipedia article on Catalan numbers (\cite{rukavicka2011generalized}). 
\begin{defn}[Grid Walk]\label{defn grid walk}
    Given $x,x',y,y' \in \mathbb{Z}$, a \textit{grid walk from $(x,y)$ to $(x',y')$} is a sequence of $(x'-x) + (y'-y) + 1$ coordinates $\(z_0,z_1,\ldots,z_{(x'-x) + (y'-y)}\)$ where
    \begin{enumerate}
        \item $z_i=(x_i,y_i)$ where $x_i\in [x,x']$ and $y_i\in [y,y']$ for each $i \in [x+y]$.
        \item $z_0=(x,y)$ and $z_{(x'-x) + (y'-y)}=(x',y')$. 
        \item $z_{i}-z_{i-1}=(1,0)$ or $(0,1)$ for all $i\in[(x'-x) + (y'-y)]$.
    \end{enumerate} 
    Pictorially, a grid is a walk from $(x,y)$ to $(x',y')$ that steps on integer coordinates and only moves straight up or straight right.
\end{defn}
\begin{defn}
    Given a line $l$ of the form $y = kx + b$ where $k > 0$, we say that a grid walk $w = (z_0,z_1,\ldots,z_{j})$ is weakly below $l$ if $y_i \leq kx_i + b$ for all $i \in [j] \cup \{0\}$.

    An important special case is that a grid walk $(z_0,z_1,\dots, z_{x+y})$ from $(0,0)$ to $(x,y)$ is weakly below the diagonal if ${y_i/x_i\leq y/x}$ for all $i \in [x+y]$.
    
\end{defn}
It is often useful to apply translations to grid walks.
\begin{defn}[Grid walk translations]
Given a grid walk $w = (z_0,z_1,\ldots,z_{j})$ and $x,y \in \mathbb{Z}$, we define $w + (x,y)$ to be the grid walk $w + (x,y) = (z_0 + (x,y),z_1 + (x,y),\ldots,z_{j} + (x,y))$
\end{defn}
With these definitions, we can now prove \Cref{thm:catalan2-recurrence-relation}.
\begin{proof}[Proof of \Cref{thm:catalan2-recurrence-relation}]

\setlength{\parindent}{1.5em}
\setlength{\parskip}{1.2mm}
\renewcommand{\baselinestretch}{1.2}

Let $W_k$ be the set of all grid walks from $(0,0)$ to $(k,2k)$ weakly below the diagonal and let $d_k=|W_k|$. We will prove that $\displaystyle d_k = C'_k = \frac{1}{2k+1}\binom{3k}{k}$ and $d_k$ satisfies the recurrence relation $\displaystyle d_{k+1}=\sum_{a,b,c\geq 0: a+b+c=k} {d_a}{d_b}{d_c}$.

\begin{lemma}\label{lem:gridwalkrecurrencerelation}
For all $k \in \mathbb{N}$, 
    $\displaystyle d_{k}= \sum_{a=1}^{k}{\sum_{b=1}^{a}{d_{b-1}d_{a-b}d_{k-a}}} = \sum_{\substack{a,b,c\geq 0:\\ a+b+c=k-1}} {d_a}{d_b}{d_c} $
\end{lemma}
\begin{proof}
    
    We will establish a bijection between $W_{k}$ and $W'_{k}:=\bigcup_{a,b: 1 \leq b \leq a \leq k}{W_{b-1} \times W_{a-b} \times W_{k-a}}$.

    Let $w=(z_0,z_1,\dots,z_{3k})$ be a grid walk from $(0,0)$ to $(k,2k)$ weakly below the diagonal. We construct the walks $w'_1,w'_2,w'_3$ in reverse order.
    
    Consider the first point $(a,2a)$ after the starting point $(0,0)$ where $w$ touches the diagonal. In other words, let $a$ be the smallest element of $[k]$ such that $z_{3a} = (a,2a)$. Then $w_3=(z_{3a},z_{3a+1},\dots,z_{3k})$ is a grid walk from $(a,2a)$ to $(k,2k)$ weakly below the diagonal. Taking $w'_3 = w_3 - (a,2a)$, $w'_3\in W_{k-a}$. 
    
    Let $l$ be the line parallel to the diagonal which passes through $(a,2a-1)$. In other words, $l$ is the line $y = 2x - 1$. Since $z_{3a}$ is the first point touching the diagonal, $(z_1,\dots, z_{i-1})$ is weakly below $l$. Let $z_{3b-1}=(b,2b-1)$ be the first point touching $l$. Then $w_2=(z_{3b-1},\dots,z_{3a-1})$ is a grid walk from $(b,2b-1)$ to $(a,2a-1)$ weakly below $l$. Taking $w'_2 = w_2 - (b,2b-1)$, $w'_2\in W_{a-b}$.
    
    Let $l'$ be the line parallel to the diagonal which passes through $(b,2b-2)$. In other words, $l$ is the line $y = 2x - 2$. Since $z_{3b-1}$ is the first point touching $l$ and $z_1 = (1,0)$, $(z_1,\dots, z_{3b-2})$ is weakly below $l'$. Observe that $w_3=(z_1,\dots, z_{3b-2})$ is a grid walk from $(1,0)$ to $(b,2b-2)$ weakly below $l'$. Taking $w'_1 = w_1 - (1,0)$, $w'_1 \in W_{b-1}$.
    
    Thus from $w\in W_k$ we get a tuple $(w'_1,w'_2,w'_3)\in W_{b-1} \times W_{a-b} \times W_{k-a}$ where $a,b$ are uniquely determined by $w$ and $1 \leq b \leq a \leq k$.
    Conversely, given $a,b$ such that $1 \leq b \leq a \leq k$ and $(w'_1,w'_2,w'_3)\in W_{b-1} \times W_{a-b} \times W_{k-a}$, let
    \begin{equation*}
        w=\((0,0), w_1+(1,0), w_2+(b,2b-1), w_3+(a,2a)\)
    \end{equation*}
    It is not hard to check that $w \in W_k$ and this is a bijection. 
     \begin{figure}[hbt!]
        \centering
        \begin{subfigure}{.3\textwidth}
            \includegraphics[width=1\linewidth]{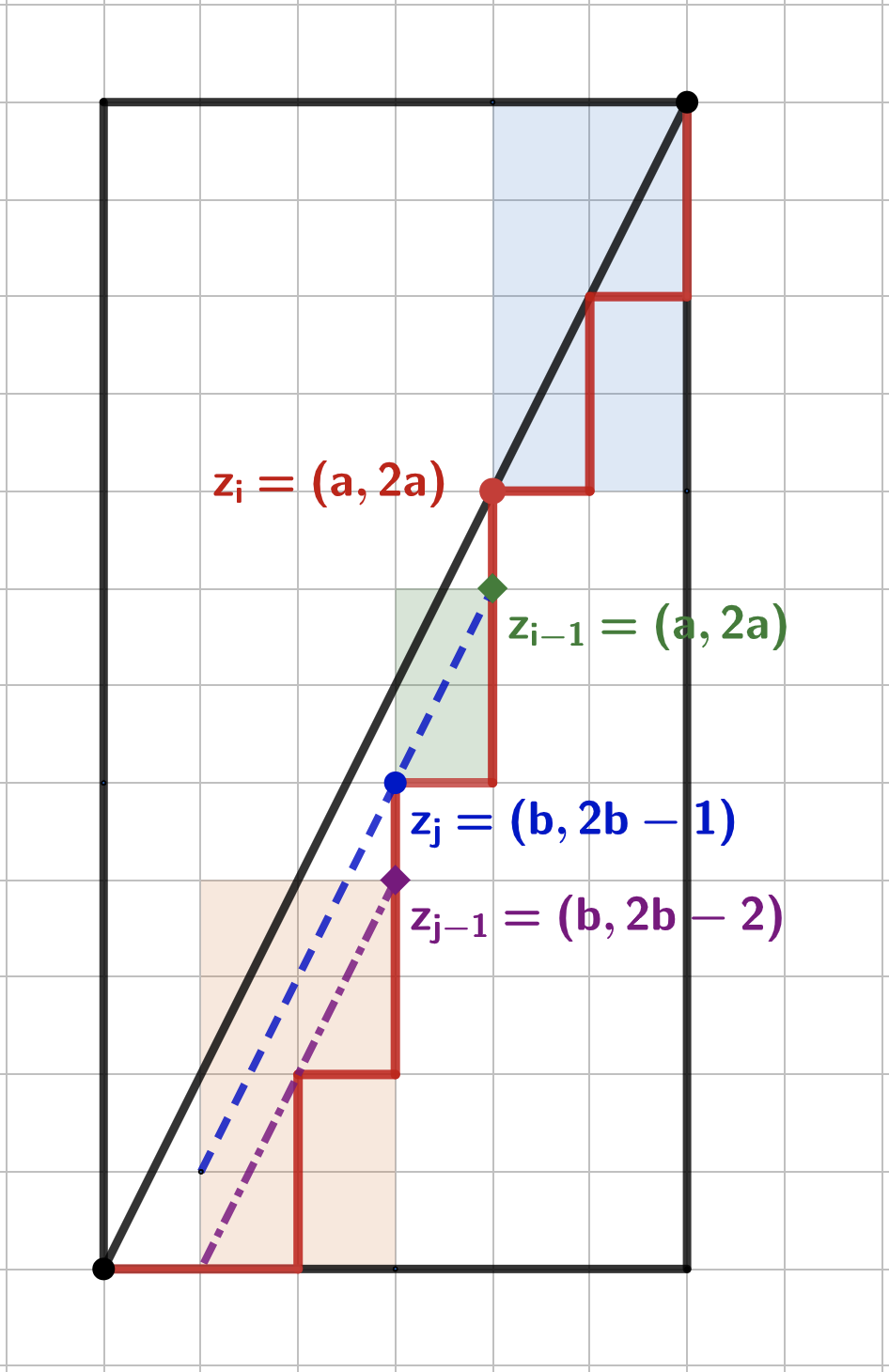}
            \caption{proof part 1}
            \label{fig:gridwalk 1}
        \end{subfigure}
        \begin{subfigure}{.66\textwidth}
            \includegraphics[width=1\linewidth]{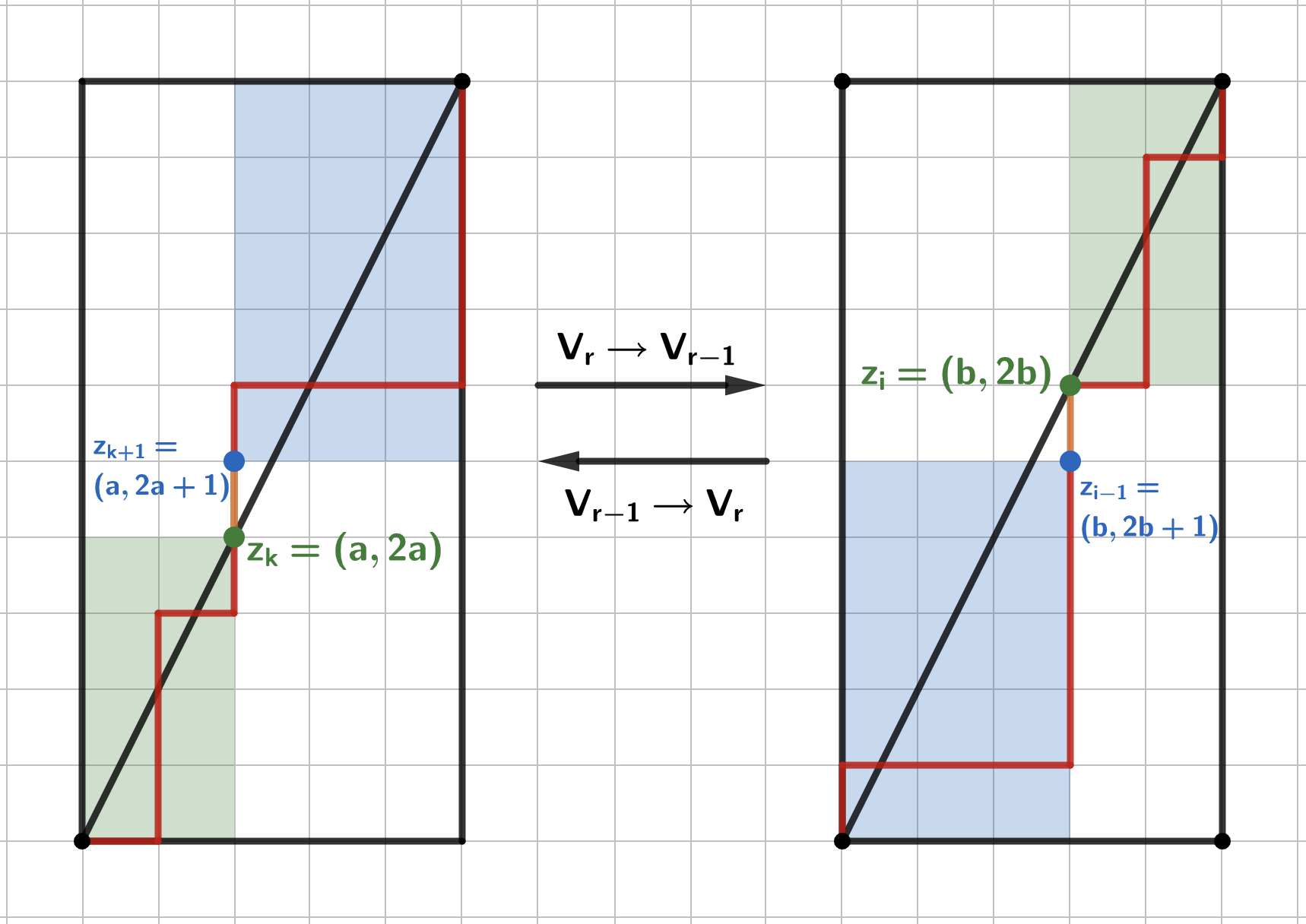}
            \caption{proof part 2}
            \label{fig:gridwalk 2}
        \end{subfigure}
        \caption{}
        \label{fig:grid walk}
    \end{figure}

\end{proof}
\begin{lemma}\label{lem:gridwalkcounting}
    $\displaystyle{d_k = C'_k = \frac{1}{2k+1}\binom{3k}{k}}$. 
\end{lemma}
\begin{proof}
    For $r\in\{0,1,\dots,2k\}$, let $V_r$ be the set of grid walks from $(0,0)$ to $(k,2k)$ that have $r$ vertical steps above the diagonal. In other words, for each $w=(z_0,z_1,\dots, z_{3k})\in V_r$, writing $z_j=(x_j,y_j)$, there are $r$ indices $j \in [3k]$ such that $y_j >2x_j$ and $z_j - z_{j-1} = (0,1)$. Let $G_k$ be the set of all grid walk from $(0,0)$ to $(k,2k)$. We have that $\displaystyle \abs{G_n}=\binom{3k}{k}$. Note that $V_0=W_k$ and $\displaystyle \bigcup_{r=0}^{2k} V_r=G_k$. We will prove that $|V_r|=|V_{r-1}|$ for all $r\in [2k]$. This implies that $\displaystyle \abs{V_0}=d_k=\frac{1}{2k+1}\binom{3k}{k}$ as needed.
    
    \begin{lemma}
        $|V_{r}| = |V_{r-1}|$ for all $r\in [2k]$.
    \end{lemma}
    \begin{proof}
        For each $r \in [2k]$, we have the following bijection between $V_{r}$ and $V_{r-1}$. 
    \begin{enumerate}
        \item Given $w=(z_0,\dots, z_{3k}) \in V_r$, let $z_{j}$ be the \textbf{last} point where the walk is on the diagonal and then takes a step upwards. In other words, $j\in [3k-1] \cup \{0\}$ is the largest index such that $z_j=(a,2a)$ for some $a\in[k]$ and $z_{j+1}-z_{j}=(0,1)$. 
        
        Let $w_1=(z_0,\dots, z_{j})$ and let $w_2=(z_{j+1},\dots, z_{3n})$. Let $w'=(w'_1,w'_2)$ where $w'_1 = w_2-(a,2a+1)$ and $w'_2 = w_1+(k-a,2k-2a))$ (see \Cref{fig:grid walk}, $w'$ exchanges the green and blue part of $w$). Observe that $w_2'$ has the same number of steps above the diagonal as $w_1$ does. Since $z_j$ is the last point such that $w$ is on the diagonal and then takes a step upwards, $w'_1$ has the same number of vertical steps above the diagonal as $w_2$. However, the vertical step from $z_j = (a,2a)$ to $z_{j+1} = (a,2a+1)$ in $w$ (i.e., the step between the end of $w_1$ and the start of $w_2$) is above the diagonal while the vertical step from $(k-a,2(k-a)-1)$ to $(k-a,2(k-a))$ in $w'$ (i.e., the step between the end of $w'_1$ and the start of $w'_2$) is below the diagonal. Thus, $w'$ has one less vertical step above the diagonal than $w$ so $w'\in V_{r-1}$.
        
        \item Let $w = (z_0,\dots, z_{3k})\in V_{r-1}$. Let $z_j$ be the first point such that $w$ touches the diagonal from below. In other words, $j \in [3k]$ is the first index such that $z_j=(b,2b)$ for some $b \in [k]$ and $z_j-z_{j-1}=(0,1)$. Let $w_1=(z_0,\dots,z_{j-1})$ and $w_2=(z_{j},\dots,z_{3k})$. Let $w'=(w_1',w_2')$ where $w_1' = w_2-(b,2b)$ and $w_2' = w_1 + (k-b,2k-2b+1)$. Observe that $w_1'$ has the same number of steps above the diagonal as $w_2$ does. Since $z_j$ is the first point where $w$ arrives at the diagonal after a vertical step, $w'_2$ has the same number of vertical steps above the diagonal as $w_1$. However, the vertical step $z_{j-1} = (b,2b-1)$ to $z_{j} = (b,2b)$ in $w$ (i.e., the step between the end of $w_1$ and the start of $w_2$) is below the diagonal while the vertical step from $((k-b,2(k-b))$ to $((k-b,2(k-b) + 1)$ in $w'$ (i.e., the step between the end of $w'_1$ and the start of $w'_2$) is below the diagonal. Thus, $w'$ has one less vertical step above the diagonal than $w$ so $w'\in V_{r-1}$.
    \end{enumerate}
    It is not hard to check that this gives a bijection.
    \end{proof}
\end{proof}
    \Cref{thm:catalan2-recurrence-relation} follows directly from \Cref{lem:gridwalkrecurrencerelation} and \Cref{lem:gridwalkcounting}. By \Cref{lem:gridwalkcounting}, $\displaystyle{d_k = C'_k}$. By \Cref{lem:gridwalkrecurrencerelation}, $\displaystyle d_{k} = \sum_{\substack{a,b,c\geq 0:\\ a+b+c=k-1}} {d_a}{d_b}{d_c} $. Thus, $\displaystyle C'_{k+1} = \sum_{\substack{a,b,c\geq 0:\\ a+b+c=k}} {C'_a}{C'_b}{C'_c}
               =\sum_{a=0}^{k}{ {C'_a} \(\sum_{b=0}^{k-a} {C'_b}{C'_{k-a-b}}\)}$, as needed.  
\end{proof}
\subsection{Dominant constraint graphs on a cycle}\label{sec:dominantcycleconstraintgraphs}
In order to count the number of dominant constraint graphs in $\mathcal{C}_{\(\a_Z,2k\)}$, we need a few properties of these constraint graphs. As a warm-up, we first analyze dominant constraint graphs on a cycle of length $2k$. This is essentially the trace power analysis for Wigner matrices (\cite{wigner1993characteristic}).

Recall the \emph{line shape} $\a_0$ and its corresponding $H\(\a_0, 2k\)$ from \Cref{defn:line-shape} and \Cref{defn:line-shape-H-label}.

\subsubsection{Isolated vertices in constraint graphs}

When analyzing constraint graphs, it is often useful to consider vertices which are not incident to any constraint edges. We call such vertices \emph{isolated}. 
\begin{defn}\label{defn:isolated vertex in constr graph}
Given a multi-graph $H$ and a constraint graph $C$ on $H$, we say that a vertex $v \in V(C) = V(H)$ is \emph{isolated} if $v$ is not incident to any edge in $E(C)$.
\end{defn}
\begin{prop}\label{prop:line shape isolated vertex}
If $C$ is a nonzero-valued constraint graph on $H(\a_0,2k)$ and $C$ has an isolated vertex $x_j$, then $x_{j-1}\=x_{j+1}$. In the cases when $j=1$ or $j=2k$, we take $x_{0}=x_{2k}$ and $x_{2k+1} = x_1$ respectively. 
\end{prop}

\begin{proof}
    Recall that by \Cref{prop:nonzero iff even multi-edges}, $C$ is nonzero-valued if and only if each edge in $H(\a_0,2k)/C$ appears an even number of times. Since $x_j$ is isolated, the only way this can happen is if $x_{j-1}\=x_{j+1}$.
\end{proof}

\subsubsection{Properties of dominant constraint graphs on a cycle}
We now describe several properties of dominant constraint graphs in $\C_{(\a_0,2k)}$.
\begin{defn}\label{defn:cyclewellbehaved}
We say that a constraint graph $C \in \C_{(\a_0,2k)}$ is \emph{well-behaved} if all constraint edges in $C$ are either between two copies of $u$ (i.e., two vertices in $\{u_1,\ldots,u_k\}$) or between two copies of $v$ (i.e., two vertices in $\{v_1,\ldots,v_k\}$). Equivalently, $C$ is well-behaved if whenever $x_i \= x_j$, $j-i$ is even.
\end{defn}

\begin{defn}\label{defn:cyclenoncrossing}
We say that a constraint graph $C \in \C_{(\a_0,2k)}$ is \emph{non-crossing} if there are no indices $i_1,i_2,i_3,i_4 \in [k]$ such that $i_1 < i_2 < i_3 < i_4$, $x_{i_1} \= x_{i_3}$, $x_{i_2} \= x_{i_4}$, and there is no path from $x_{i_2}$ to $x_{i_3}$ in $C$. Equivalently, $C$ is \emph{non-crossing} if there is a constraint graph $C'$ such that $C' \equiv C$ and if we draw $C'$ so that the vertices $x_1,\ldots,x_{2k}$ are on a circle and the edges of $C'$ are chords of this circle, no two edges of $C'$ cross.
\end{defn}

\begin{lemma}\label{lem:line-shape-catalan-lemma}
For all $k \in \mathbb{N}$, all dominant constraint graphs in $\mathcal{C}_{(\a_0,2k)}$ are well-behaved, non-crossing, and have $k-1$ edges.
\end{lemma}
\begin{proof}
    We first show that all dominant constraint graphs in $\mathcal{C}_{(\a_0,2k)}$ have $k-1$ edges. To show this, we make the following observations:
    \begin{enumerate}
        \item The constraint graph $C_0$ with constraint edges $E(C_0) = \left\{\{x_{2j-1},x_{2j+1}\}: j \in [k-1]\right\}$ has value $1$ so dominant constraint graphs in $\mathcal{C}_{(\a_0,2k)}$ have at most $k-1$ edges.
        \item If $C \in \C_{(\a_0,2k)}$ and $\val(C) \neq 0$ then every edge in $H(\a_0,2k)/C$ has even multiplicity so $H(\a_0,2k)/C$ has at most $k$ distinct edges. Since $H(\a_0,2k)/C$ is connected, this implies that $H(\a_0,2k)/C$ has at most $k+1$ vertices so $|E(C)| = \abs{V\(H(\a_0,2k)\)} - \abs{V\(H(\a_0,2k)/C\)} \geq k-1$.
    \end{enumerate}
    We prove that all dominant constraint graphs in $\mathcal{C}_{(\a_0,2k)}$ are well-behaved and non-crossing by induction. The base case $k = 1$ is trivial.

    For the inductive step, assume that all dominant constraint graphs in $\mathcal{C}_{(\a_0,2(k-1))}$ are well-behaved and non-crossing. Let $C \in \mathcal{C}_{(\a_0,2k)}$ be a dominant constraint graph. There are two cases to consider.
    \begin{enumerate}
        \item $C$ does not have any isolated vertices. In this case, $C$ must have at least $k$ edges as $|V(C)| = |V(H(\a_0,k))| = 2k$. This is a contradiction as dominant constraint graphs in $\mathcal{C}_{(\a_0,2(k-1))}$ have $k-1$ edges.
        \item $C$ has an isolated vertex $x_j$. In this case, by \Cref{prop:line shape isolated vertex}, $x_{j-1} \= x_{j+1}$ (where $x_0 = x_{2k}$ and $x_{2k+1} = x_1$). Without loss of generality, we can assume that $\{x_{j-1},x_{j+1}\} \in E(C)$.
        
        Let $H'$ be the multi-graph obtained from $H$ by deleting the vertex $x_j$, deleting the edges $\{x_{j-1},x_j\}$ and $\{x_j,x_{j+1}\}$, and merging the vertices $x_{j-1}$ and $x_{j+1}$. Similarly, let $C'$ be the constraint graph obtained from $C$ by deleting the vertex $x_j$, deleting the constraint edge $\{x_{j-1},x_{j+1}\}$, and merging the vertices $x_{j-1}$ and $x_{j+1}$. 
        
        Observe that $H'$ is isomorphic to $H(\a_0,k-1)$ so we can view $C'$ as a constraint graph on $H(\a_0,k-1)$. Moreover, $\val(C') = \val(C) \neq 0$ and $|E(C')| = |E(C)| - 1 = k-2$, so $C'$ is a dominant constraint graph on $H(\a_0,k-1)$. By the inductive hypothesis, $C'$ is well-behaved and non-crossing and it is not hard to check that $C$ is well-behaved and non-crossing as well.
    \end{enumerate}
\end{proof}

\subsection{Recurrence relation for dominant constraint graphs on \texorpdfstring{$H(\a_0,2k)$}{H(a0, 2k)}}\label{sec:cyclerecurrencerelation}
\begin{figure}[hbt!]
    \centering
    \includegraphics[scale=0.32]{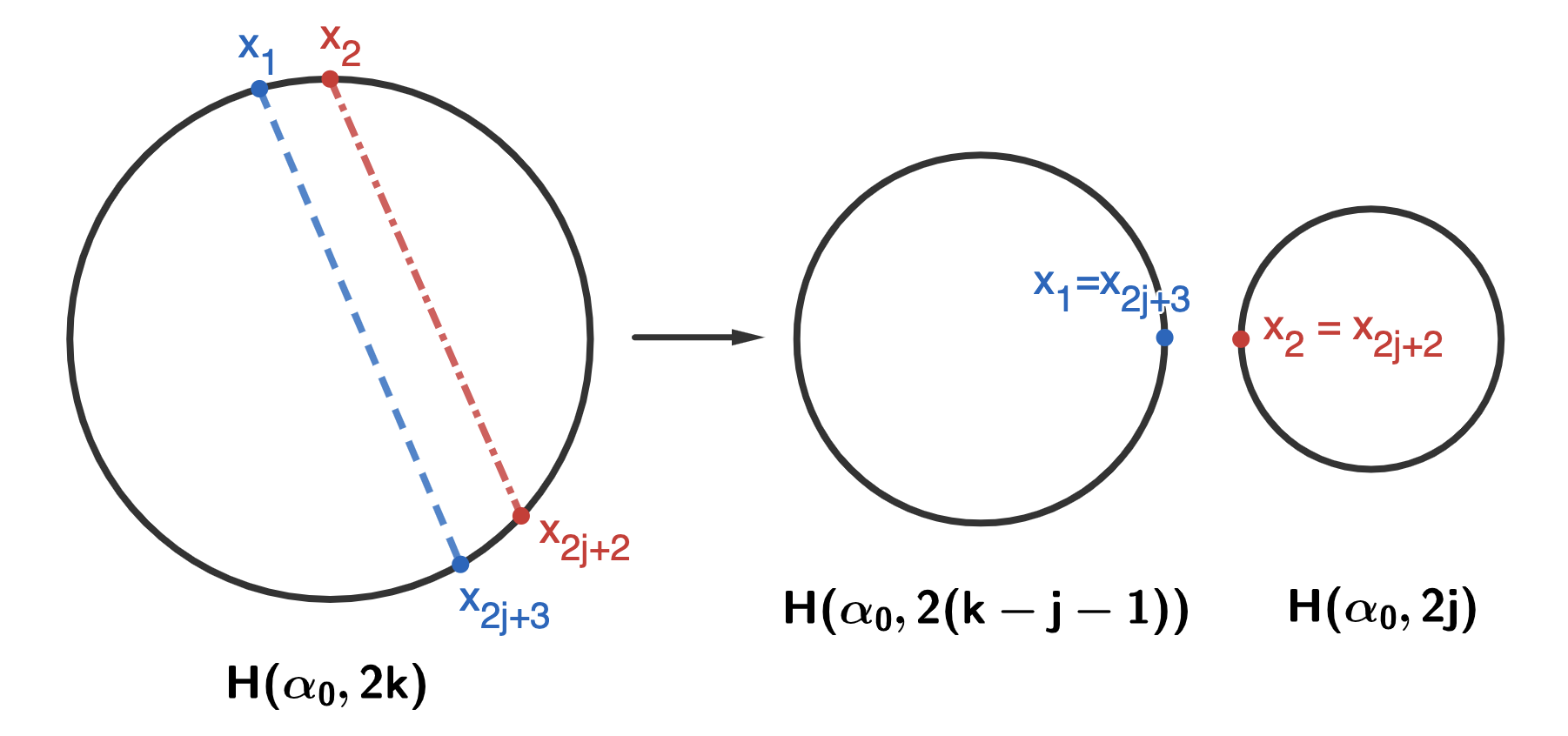}
    \caption{Illustration of \Cref{lem:lineshapepartitioning} case 2: $E = \{\{x_1,x_{2j+3}\},\{x_2,x_{2j+2}\}\}$. This is also an illustration of \Cref{cor:line-shape-recurrence}, case 2: $j\in[1,k-2]$ is the smallest index such that $x_1 \longleftrightarrow_C x_{2j+3}$.}
    \label{fig:line-shape-split}
\end{figure}

We now show that if $C_k = \abs{C \in \mathcal{C}_{(\a_0,2k)}: C \text{ is dominant}}$ then $\displaystyle C_k = \sum_{j=0}^{k-1}\,{C_{j}C_{k-1-j}}$. Since $C_0 = C_1 = 1$, this implies that $C_k = \frac{1}{k+1}\binom{2k}{k}$ is the $k$th Catalan number.
\begin{lemma}\label{lem:lineshapepartitioning}
The following sets of constraint edges partition $H(\a_0,2k)$ into subgraphs $H_1,H_2$.

\begin{enumerate}
\item If $E = \{\{x_1,x_3\}\}$ then $E$ partitions $H(\a_0,2k)$ into $H_1$ and $H_2$ where $V(H_1) = \emptyset$ and $V(H_2) = \{x_1\} \cup \{x_4,\ldots,x_{2k}\}$. Observe that $H_1 \simeq H(\a_0,0)$ and $H_2 \simeq H(\a_0,k-1)$

\item If $E = \{\{x_1,x_{2j+3}\},\{x_2,x_{2j+2}\}\}$ for some $j \in [k-2]$ then $E$ partitions $H(\a_0,2k)$ into $H_1$ and $H_2$ where $V(H_1) = \{x_2,\ldots,x_{2j+1}\}$ and $V(H_2) = \{x_1\} \cup \{x_{2j+4},\ldots,x_{2k}\}$. Observe that $H_1 \simeq H(\a_0,j)$ and $H_2 \simeq H(\a_0,k - 1 - j)$.

\item If $E = \{\{x_2,x_{2k}\}\}$ then $E$ partitions $H$ into $H_1$ and $H_2$ where $V(H_1) = \{x_2,\ldots,x_{2k-1}\}$ and $V(H_2) = \emptyset$. Observe that $H_1 \simeq H(\a_0,k-1)$ and $H_2 \simeq H(\a_0,0)$.
\end{enumerate}

See \Cref{fig:line-shape-split} for an illustration.
\end{lemma}
\begin{proof}
This follows from \Cref{lem:Hpartitionconditions}. For each of these sets of edges $E$, the first condition of \Cref{defn:Hpartitions} and the first condition of \Cref{lem:Hpartitionconditions} can be seen directly. The second and third conditions of \Cref{lem:Hpartitionconditions} are trivial as $V(H_1) \cap V(H_2) = \emptyset$. The fourth condition of \Cref{lem:Hpartitionconditions} is trivial in cases $1$ and $3$. In case $2$, the fourth condition of \Cref{lem:Hpartitionconditions} follows from \Cref{lem:line-shape-catalan-lemma} which says that all dominant constraint graphs $C$ on $H(\a_0,2k)$ are non-crossing and well-behaved. To see this, let $C$ be a dominant constraint graph on $H(\a_0,2k)$ and observe that for each pair of edges $e_1,e_2 \in E(H)$ such that $e_1,e_2$ appear with odd multiplicity in $H/{E_C}$, $e_1$ appears in $H_1$, and $e_2$ appears in $H_2$, one of the endpoints of $e_1$ must be $x_{i}$ for some odd $i \in [3,2j+1]$ and one of the endpoints of $e_2$ must be $x_{i'}$ for some odd $i' \in [2j+3,2k] \cup \{1\}$. Since $C$ contains the edge $\{x_2,x_{2j+1}\}$, is well-behaved, and is non-crossing, we cannot have that $x_{i} \=_{C} x_{i'}$ so we cannot have that $e_1 \=_{C} e_2$.
\end{proof}

\begin{cor}\label{cor:line-shape-recurrence}
Letting $C_k = \abs{C \in \mathcal{C}_{(\a_0,2k)}: C \text{ is dominant}}$, $\displaystyle C_k = \sum_{j=0}^{k-1}\,{C_{j}C_{k-1-j}}$
\end{cor}
\begin{proof}
We partition the dominant constraint graphs $C \in \mathcal{C}_{(\a_0,2k)}$ based on the first index $j' \in [2k]$ such that $j' > 1$ and $x_{1} \=_C x_{j'}$. If no such $j'$ exists then we set $j' = 2k+1$.

We must have that $j'$ is odd and $j' \geq 3$ so we can write $j' = 2j + 3$ for some $j \in \{0,1,\ldots,j-1\}$. We have the following cases:
\begin{enumerate}
\item If $j = 0$ then this gives case $1$ of \Cref{lem:lineshapepartitioning}. There are $C_{0}C_{k-1}$ such dominant constraint graphs $C \in \mathcal{C}_{(\a_0,2k)}$. 
\item If $0 < j < k-1$ then since $C$ is non-crossing, we must have that $x_2 \=_C x_{2j + 2}$ as otherwise the edges $\{x_1,x_2\}$ and $\{x_{2j+2},x_{2j+3}\}$ would appear with odd multiplicity in $H/C$ and we would have that $\val(C) = 0$. This gives case $2$ of \Cref{lem:lineshapepartitioning} and there are ${C_j}C_{k-1-j}$ such dominant constraint graphs $C \in \mathcal{C}_{(\a_0,2k)}$.
\item If $j = k-1$ then we must have that $x_2 \=_C x_{2k}$ as otherwise the edges $\{x_1,x_2\}$ and $\{x_{2k},x_{1}\}$ would appear with odd multiplicity in $H/C$ and we would have that $\val(C) = 0$. This gives case $3$ of \Cref{lem:lineshapepartitioning} and there are $C_{k-1}C_0$ such dominant constraint graphs $C \in \mathcal{C}_{(\a_0,2k)}$.
\end{enumerate}

See \Cref{fig:line-shape-split} for an illustration.
\end{proof}

\subsection{Properties of dominant constraint graphs on \texorpdfstring{$H(\a_{Z},2k)$}{H(aZ, 2k)}}

Now that we have analyzed dominant constraint graphs in $\mathcal{C}_{(\a_0,2k)}$, we can analyze dominant constraint graphs in $\mathcal{C}_{\(\a_Z,2k\)}$.

\begin{defn}\label{defn:z-shape constraint graph}
    Let $\a_Z$ be the Z-shape as defined in \Cref{defn:z-shape} and let $H\(\a_Z,2k\)$ be the multi-graph as defined in \Cref{def:copies}. We label the vertices of $V_{\(\a_Z\)_i}$ as $\{a_{i1},a_{i2},b_{i1},b_{i2}\}$ and the vertices of $V_{\(\a_Z^T\)_i}$ as $\{b_{i1},b_{i2},a_{(i+1)1},a_{(i+1)2}\}$. We call the induced subgraph of $H\(\a_Z,2k\)$ on vertices $\left\{a_{i1},b_{i1}:i\in[k]\right\}$ the \textit{outer wheel $W_1$} and the induced subgraph on vertices $\{a_{i2},b_{i2}:i\in[k]\}$ the \textit{inner wheel $W_2$}. We denote the vertices of $W_i$ as $V_i$ and edges as $E_i$. 
    
    We label the ``middle edges" of $H(\a,2k)$ in the following way: let $e_{2i-1}=\{a_{i2},b_{i1}\}$ and $e_{2i}=\{b_{i1},a_{(i+1)2}\}$ for $i \in [k]$. We call the edges $\{e_i: i \in [2k]\}$ the \textit{spokes} of $H\(\a_Z,2k\)$. See \Cref{fig:zshape-H-2k} for an illustration.
    
    \begin{figure}[hbt!]
        \centering
        \includegraphics[scale=0.32]{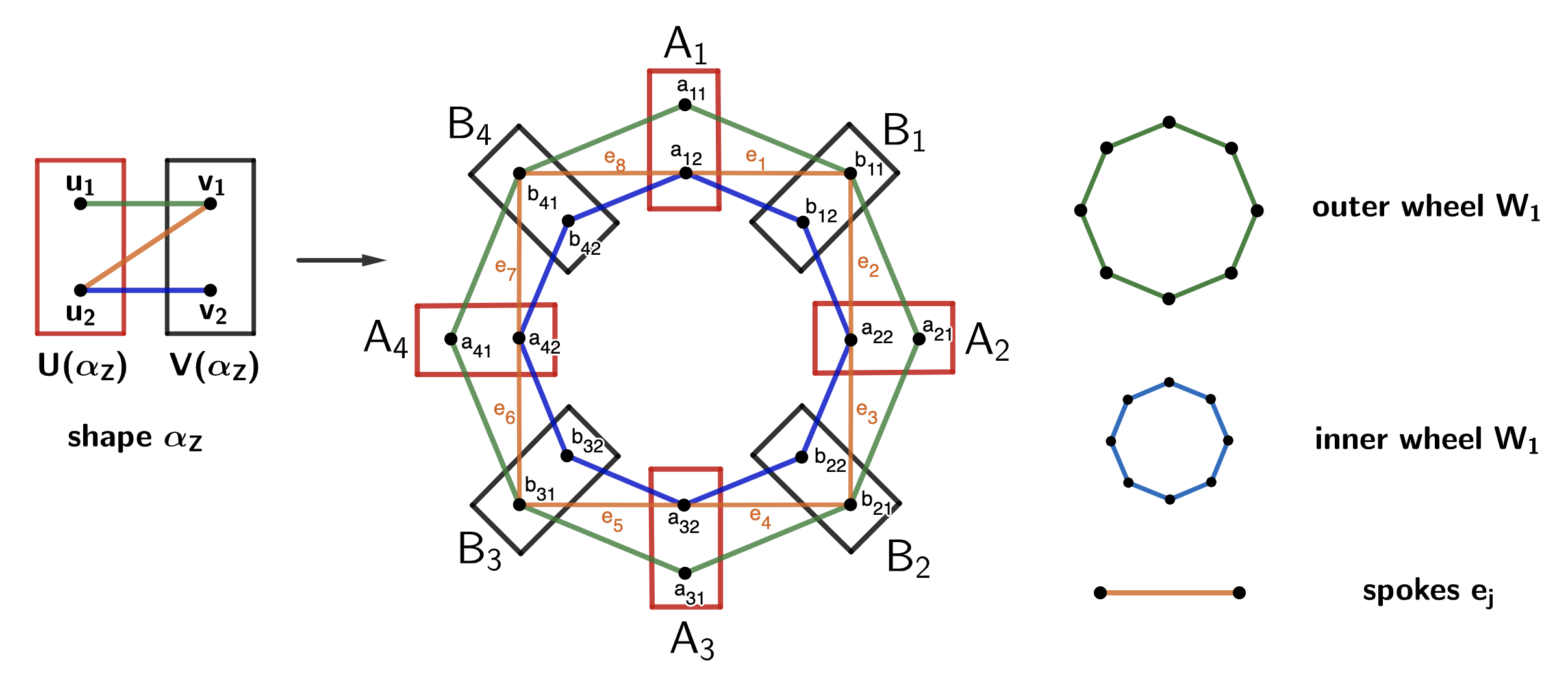}
        \caption{$H(\a_{Z},2k)$ where $k=4$.}
        \label{fig:zshape-H-2k}
    \end{figure}
\end{defn}

\begin{defn}\label{defn:z-shape induced constraint graph}
Let $C$ be a constraint graph in $\mathcal{C}_{\(\a_Z,2k\)}$. For $i=1,2$, we take $C_i$ to be the induced constraint graph of $C$ on the vertices $V_i$.
\end{defn}

\begin{rmk}
Each wheel $W_i$ can be viewed as $H(\a_0,2k)$. The induced constraint graph $C_i$ can be viewed as a constraint graph on $H(\a_0,2k)$.
\end{rmk}

The key property that we need about dominant constraint graphs in $\mathcal{C}_{\(\a_Z,2k\)}$ is that they are well-behaved. This implies that the induced constraint graphs $C_1,C_2$ are dominant constraint graphs in $\mathcal{C}_{(\a_0,2k)}$.

\begin{defn}\label{defn:wellbehaved}
Given a shape $\a$, we say that a constraint graph $C \in \mathcal{C}_{(\a,2k)}$ is \emph{well-behaved} if whenever $u \=_C v$, $u$ and $v$ are copies of the same vertex in $\a$ or $\a^T$.
\end{defn}

\begin{thm}\label{thm:dominantiswellbehaved}
All dominant constraint graphs in $\mathcal{C}_{\(\a_Z,2k\)}$ are well-behaved.
\end{thm}
This theorem is surprisingly tricky to prove, so we defer its proof to the appendix.
\begin{rmk}
This theorem is not true for all shapes $\a$. In particular, this theorem is false for the bipartite shape $\a$ with $U_{\a} = (u_1,u_2)$, $V_{\a} = (v_1,v_2)$, $V(\a) = U_{\a} \cup V_{\a}$, and $E(\a) = \left\{\{u_1,v_1\}, \{u_1,v_2\}, \{u_2,v_1\}, \{u_2,v_2\}\right\}$.
\end{rmk}
\begin{prop}
If $C \in \mathcal{C}_{\(\a_Z,2k\)}$ is a dominant constraint graph then $E(C) = E(C_1) \cup E(C_2)$ and the induced constraint graphs $C_1$ and $C_2$ are dominant constraint graphs on $W_1$ and $W_2$.
\end{prop}
\begin{proof}
Since $C$ is well-behaved, $E(C) = E(C_1) \cup E(C_2)$. We must have that $C_1$ and $C_2$ are well-behaved as otherwise we would have that $|E(C)| > 2k-2$, which is too large.
\end{proof}
We now observe that if $C \in \mathcal{C}_{\(\a_Z,2k\)}$ is a dominant constraint graph then the corresponding constraints for which spokes are equal to each other are non-crossing (see \Cref{prop:noncrossingspokes} below). We will then use this observation to prove a key fact about dominant constraint graphs in $\mathcal{C}_{\(\a_Z,2k\)}$.
\begin{prop}\label{prop:noncrossingspokes}
If $C \in \mathcal{C}_{\(\a_Z,2k\)}$ is a dominant constraint graph then for all $i_1,i_2,i_3,i_4 \in [2k]$ such that $i_1 < i_2 < i_3 < i_4$, if $e_{i_1} \=_C e_{i_3}$ and $e_{i_2} \=_C e_{i_4}$ then $e_{i_1} \=_C  e_{i_2} \=_C e_{i_3} \=_C e_{i_4}$
\end{prop}
\begin{proof}
Observe that $e_i$ goes between $a_{(\lfloor{\frac{i}{2}}\rfloor + 1)2}$ and $b_{{\lceil{\frac{i}{2}}\rceil}1}$. Thus, if $e_{i_1} \= e_{i_3}$ and $e_{i_2} \= e_{i_4}$ then letting $C_1$ and $C_2$ be the restrictions of $C$ to $W_1$ and $W_2$ respectively, we have that 
\begin{enumerate}
\item $a_{(\lfloor{\frac{i_1}{2}}\rfloor + 1)2} \= a_{(\lfloor{\frac{i_3}{2}}\rfloor + 1)2}$ and $a_{(\lfloor{\frac{i_2}{2}}\rfloor + 1)2} \= a_{(\lfloor{\frac{i_4}{2}}\rfloor + 1)2}$. Since $C_2$ is non-crossing, 
\begin{equation*}
a_{(\lfloor{\frac{i_1}{2}}\rfloor + 1)2} \= a_{(\lfloor{\frac{i_2}{2}}\rfloor + 1)2} \= a_{(\lfloor{\frac{i_3}{2}}\rfloor + 1)2} \= a_{(\lfloor{\frac{i_4}{2}}\rfloor + 1)2}.
\end{equation*}
\item $b_{{\lceil{\frac{i_1}{2}}\rceil}1} \= b_{{\lceil{\frac{i_3}{2}}\rceil}1}$ and $b_{{\lceil{\frac{i_2}{2}}\rceil}1} \= b_{{\lceil{\frac{i_4}{2}}\rceil}1}$. Since $C_1$ is non-crossing, 
\begin{equation*}
b_{{\lceil{\frac{i_1}{2}}\rceil}1} \= b_{{\lceil{\frac{i_2}{2}}\rceil}1} \= b_{{\lceil{\frac{i_3}{2}}\rceil}1} \= b_{{\lceil{\frac{i_4}{2}}\rceil}1}.
\end{equation*}
\end{enumerate}
This implies that $e_{i_1} \=_C  e_{i_2} \=_C e_{i_3} \=_C e_{i_4}$.
\end{proof}
\begin{lemma}\label{lem:zshape-spoke-constraint}
If $C \in \mathcal{C}_{\(\a_Z,2k\)}$ is a dominant constraint graph then 
\begin{enumerate}
\item For all $i < j \in [k]$, if $a_{i1} \=_C a_{j1}$ then $a_{i2} \=_C a_{j2}$ and the spokes $\{e_{i'}: i' \in [2i-1,2j-2]\}$ can only be made equal to each other. 
\item Similarly, for all $i < j \in [k]$, if $b_{i2} \=_C b_{j2}$ then $b_{i1} \=_C b_{j1}$ and the spokes $\{e_{i'}: i' \in [2i,2j-1]\}$ can only be made equal to each other.

\begin{figure}[hbt!]
    \centering
    \begin{subfigure}[t]{.38\textwidth}
      \includegraphics[width=1\linewidth]{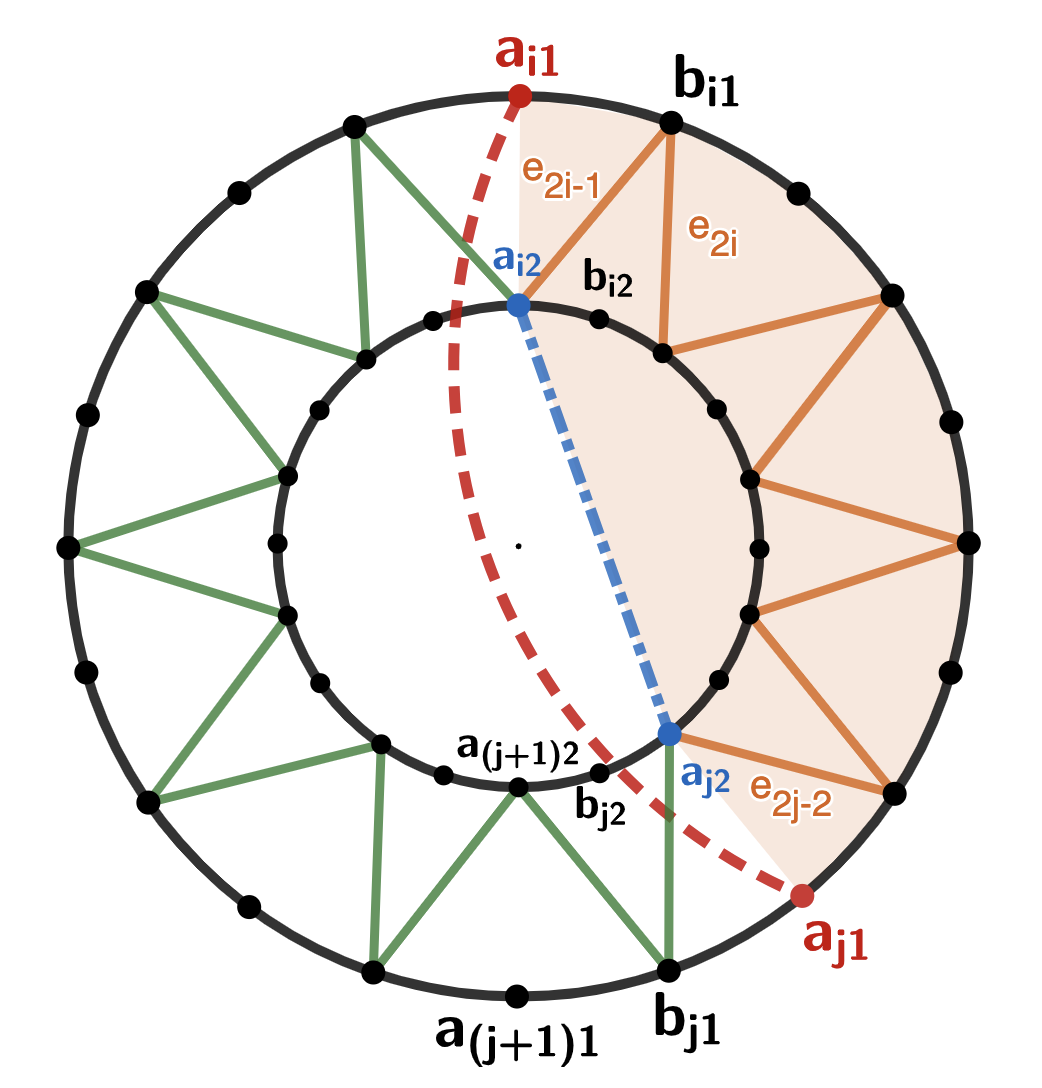}
      \caption{$a_{i1} \longleftrightarrow_C a_{j1} \implies a_{i2} \longleftrightarrow_C a_{j2}$}
      \label{fig:spoke-constraint-a}
    \end{subfigure}\hspace{1cm}
    \begin{subfigure}[t]{.38\textwidth}
      \includegraphics[width=1\linewidth]{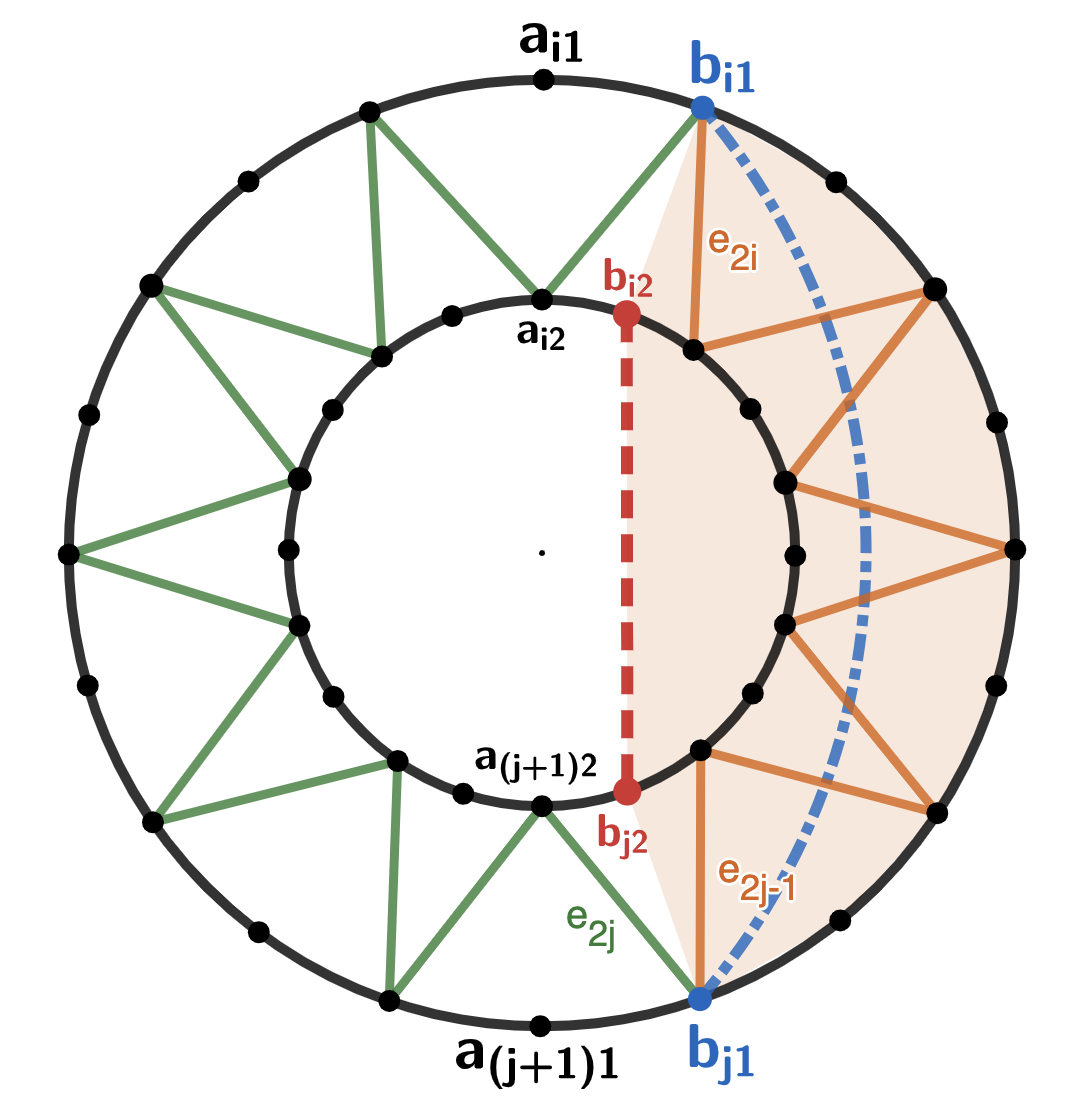}
      \caption{$b_{i2} \longleftrightarrow_C b_{j2} \implies b_{i1} \longleftrightarrow_C b_{j1}$}
      \label{fig:spoke-constraint-b}
    \end{subfigure}
    \caption{Illustration of \Cref{lem:zshape-spoke-constraint}: a spoke can only be made equal to spokes of the same color.}
    \label{fig:spoke-constraint}
\end{figure}
\end{enumerate}
\end{lemma}
\begin{proof}
We only prove the first statement as the proof of the second statement is similar. Let $C_1$ and $C_2$ be the restrictions of $C$ to $W_1$ and $W_2$ respectively.

To see that the spokes $\{e_{i'}: i' \in [2i-1,2j-2]\}$ can only be made equal to each other, observe that since $a_{i1} \=_C a_{j1}$ and $C_1$ is non-crossing, we cannot have that $b_{{\lceil{\frac{i'}{2}}\rceil}1} \= b_{{\lceil{\frac{i''}{2}}\rceil}1}$ for any $i',i''$ such that $i' \in [2i-1,2j-2]$ and $i'' \in [2k] \setminus [2i-1,2j-2]$.

To show that $a_{i2} \=_C a_{j2}$, we apply the following lemma to $\{e_{i'}: i' \in [2i-1,2j-2]\}$. 
\begin{lemma}\label{lem:matchinglemma}
For all $m \in \mathbb{N}$, if $M$ is a perfect matching on the indices $[2m]$ such that no two edges of $M$ cross (i.e. there is no pair of edges $\{i_1,i_2\},\{i_3,i_4\} \in M$ such that $i_1 < i_3 < i_2 < i_4$) then either $\{1,2m\} \in M$ or there is a sequence of indices $i_1 < \ldots < i_l$ such that $i_1 \geq 1$, $i_l < m$, and 
\begin{enumerate}
    \item $\{1,2{i_1}\} \in M$ and $\{2{i_l}+1,2m\} \in M$.
    \item For all $j \in [l-1]$, $\{2i_{j} + 1,2{i_{j+1}}\} \in M$.
\end{enumerate} 
See Figure \ref{fig:perfect-matching} for an illustration.
\end{lemma}
\begin{proof}
We prove this by induction on $m$. The base case $m = 1$ is trivial. For the inductive step, assume the result is true for $m$ and consider a matching $M$ on $[2m+2]$ such that no edges of $M$ cross. If $\{1,2m+2\} \in M$ then we are done so we can assume that $\{1,2m+2\} \notin M$.

Choose $s < t \in [2m+2]$ such that $\{s,t\} \in M$ and $t - s$ is minimized. We claim that $t = s+1$. To see this, assume that $t > s+1$. Since $M$ is a perfect matching, $\{s+1,x\} \in M$ for some $x \in [2m+2]$. Since no two edges of $M$ cross, we must have that $s+1 < x < t$. However, this implies that $x - (s+1) < t-s$, contradicting our choice of $s$ and $t$.

Now consider the matching $M'$ obtained from $M$ by deleting the indices $s,s+1$ and decreasing all indices greater than $s+1$ by $2$. By the inductive hypothesis, either $\{1,2m\}\in M$, or there is a sequence $i'_1 < \ldots < i'_l$ such that $i'_1 \geq 1$, $i'_l < m$, $\{1,2i'_1\} \in M'$, $\{2{i'_l}+1,2m\} \in M'$, and for all $j \in [l-1]$, $\{2i'_{j} + 1,2i'_{j+1}\} \in M'$. In the latter case we can modify this sequence for $M'$ as follows to obtain the desired sequence for $M$:
\begin{enumerate}
    \item Increase all indices in this sequence which are greater than or equal to $s$ by $2$. 
    \item If $s$ is odd and $i'_j = \frac{s-1}{2}$ for some $j \in [l]$, insert $i'_j + 1 = \frac{s+1}{2}$ after $i'_j$. This accounts for the fact that $2i'_j + 1 = s$ was shifted to $s+2$ so we need to add the edge $\{s,s+1\}$ between $2i'_j = s-1$ and $2i'_j + 1 + 2 = s+2$.
\end{enumerate}

If $\{1,2m\} \in M'$, then there are three cases:
\begin{enumerate}
    \item If $1<s<2m+1$ then $\{1,2m+2\} \in M$.
    \item If $s=1$ then $\{1,2\}, \{3,2m+2\} \in M$ so we can take the sequence with the single element $i_1 = 1$.
    \item If $s = 2m+1$ then $\{1,2m\}, \{2m+1,2m+2\} \in M$ so we can take the sequence with the single element $i_1 = m$.
\end{enumerate}
Applying \Cref{lem:matchinglemma} to $\{e_{i'}: i' \in [2i-1,2j-2]\}$, we have that either $e_{2i-1} \=_C e_{2j-2}$ or there is a sequence of indices $i_1 < \ldots < i_l$ such that $i \leq i_1$, $i_l < j-1$, and 
\begin{enumerate}
    \item $e_{2i-1} \=_C e_{2i_{1}}$ and $e_{2i_{l} + 1} \=_C e_{2j-2}$
    \item For all $j' \in [l-1]$, $e_{2i_{j'} + 1} \=_C e_{2i_{(j'+1)}}$
\end{enumerate}
As before, we observe that $e_{i'}$ goes between $a_{(\lfloor{\frac{i'}{2}}\rfloor + 1)2}$ and $b_{{\lceil{\frac{i'}{2}}\rceil}1}$. Thus, if $e_{2i-1} \=_C e_{2j-2}$ then $a_{i2} \= a_{j2}$. If we have a sequence of $i_1 < \ldots < i_l$ satisfying the above properties then we have that 
\begin{enumerate}
\item $a_{i2} \=_C a_{({i_1}+1)2}$ and $a_{(i_l + 1)2} \=_C a_{j2}$
\item For all $j' \in [l-1]$, $a_{(i_{j'} + 1)2} \=_C a_{(i_{j'+1} + 1)2}$
\end{enumerate}
Together, these statements imply that $a_{i2} \=_C a_{j2}$, as needed.
\end{proof}

\begin{figure}[hbt!]
    \centering
    \begin{subfigure}[b]{.42\textwidth}
      \includegraphics[width=1\linewidth]{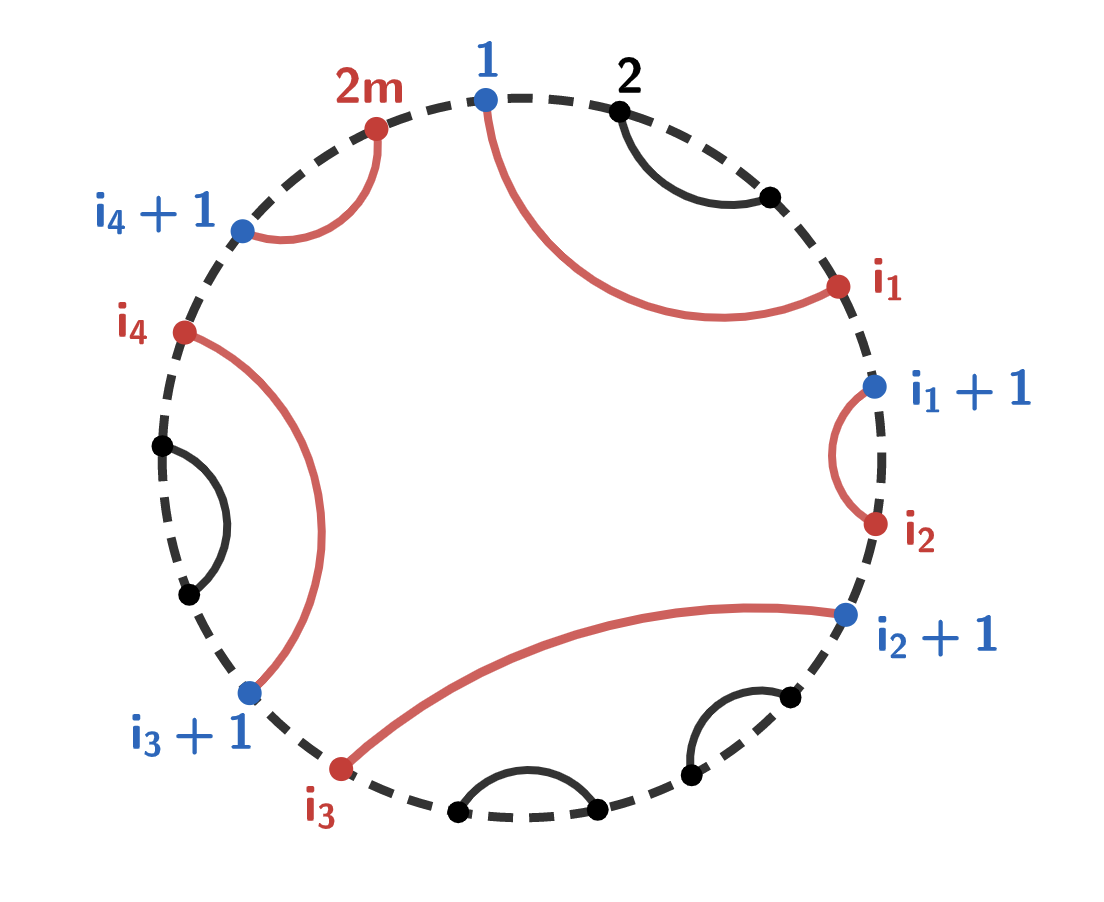}
      \caption{Illustration of Lemma \ref{lem:matchinglemma}: solid lines are matchings.}
      \label{fig:perfect-matching}
    \end{subfigure}
    \begin{subfigure}[b]{.35\textwidth}
      \includegraphics[width=1\linewidth]{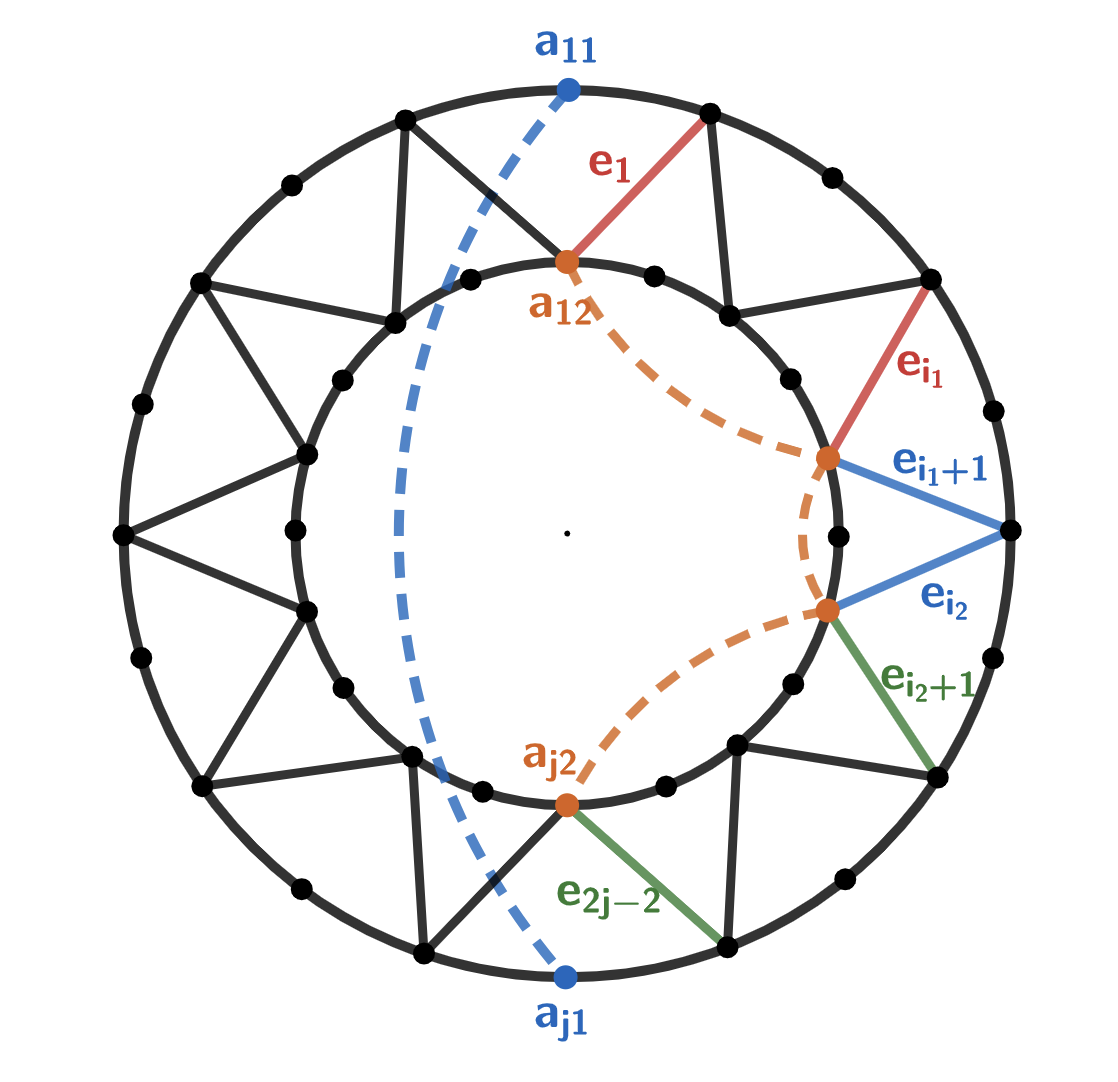}
      \caption{Illustration of Lemma \ref{lem:zshape-spoke-constraint}: $a_{11}\longleftrightarrow a_{j1}$ implies $a_{12}\longleftrightarrow a_{j2}$.}
      \label{fig:z-shape-split2}
    \end{subfigure}
    \label{fig:matching-z-split}
    \caption{Illustration of \Cref{lem:matchinglemma} and \Cref{lem:zshape-spoke-constraint}}
\end{figure}
\end{proof}

\begin{rmk}\label{rmk:zshape-spoke-constraint-rmk}

\begin{figure}[hbt!]
    \centering
    \includegraphics[scale=0.3]{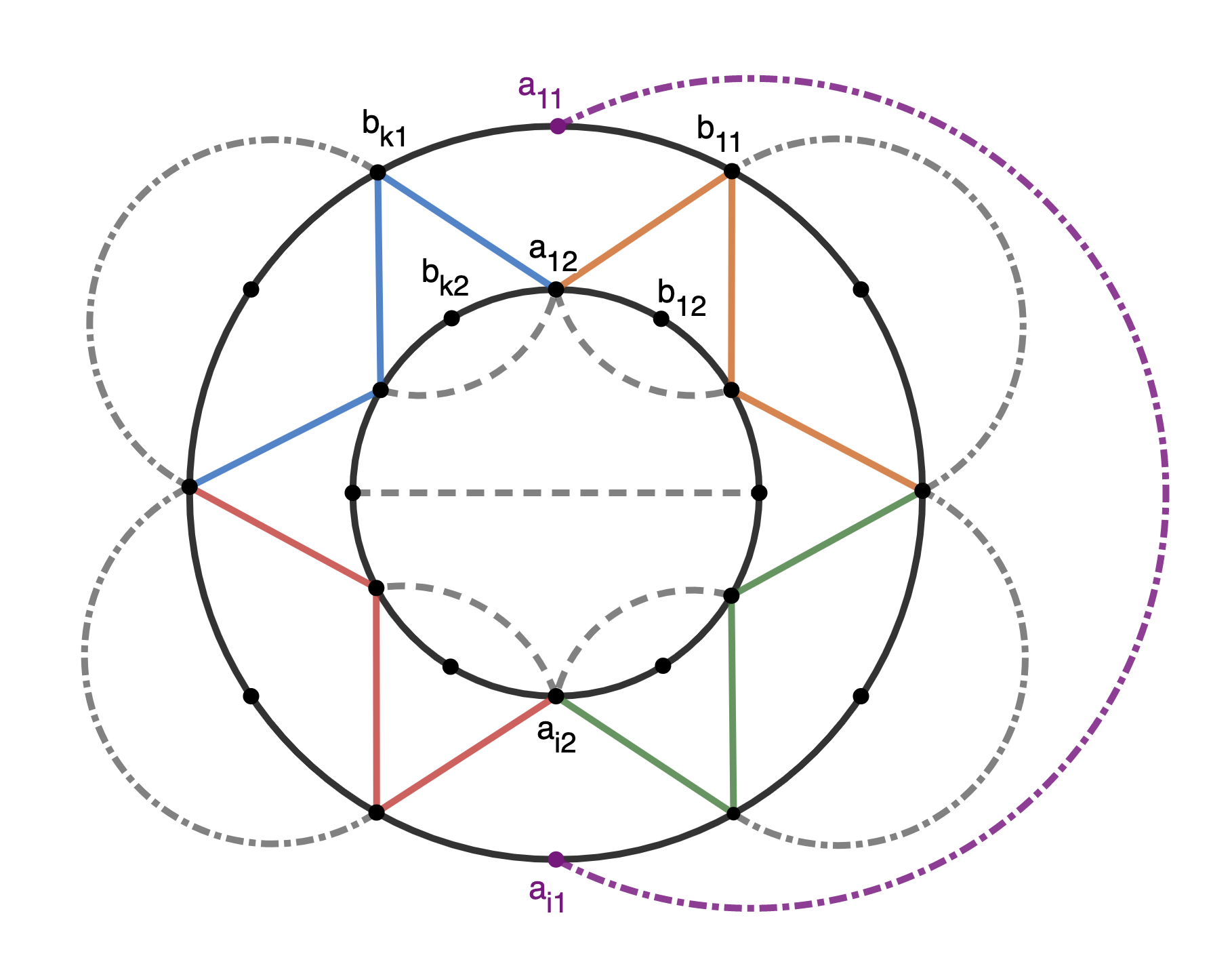}
    \caption{Illustration of \Cref{rmk:zshape-spoke-constraint-rmk}. If we only know that each edge of $H(\a_Z,2k)/C$ has multiplicity at least two rather than the stronger statement that each edge of $H(\a_Z,2k)/C$ has even multiplicity then $\a_{11}\longleftrightarrow_C a_{i1}$ does not imply $a_{12}\longleftrightarrow_C a_{i2}$.}
    \label{fig:zshape-spoke-constraint-rmk}
\end{figure}

\Cref{lem:zshape-spoke-constraint} relies on the fact that $\val(C) = 0$ whenever there is an edge in $H(\a_Z,2k)/C$ which has odd multiplicity. There are constaint graphs $C \in \mathcal{C}_{\(\a_Z,2k\)}$ which do not satisfy the conclusion \Cref{lem:zshape-spoke-constraint} where each edge in $H(\a_Z,2k)/C$ appears with multiplicity at least $2$. For an example of this, see \Cref{fig:zshape-spoke-constraint-rmk}.
\end{rmk}

\subsection{Proof of Theorem \ref{thm:num-zshape-constraint-graph}}
We are now ready to prove the main result of this section, \Cref{thm:num-zshape-constraint-graph}.

\begin{defn}\label{defn:H_j_Z}
For all $k \in \mathbb{N}$ and $j \in \{0,1,2\}$, we define $H_j\(\a_Z,2k\)$ to be the graph obtained by starting with $H\(\a_Z,2k\)$, merging the vertices $b_{1j'}$ and $b_{kj'}$ for all $j' \in [j]$, deleting the vertices $\left\{a_{1j'}: j' \in [j]\right\}$ and all edges incident to these vertices, and deleting the spokes incident to $a_{1(j+1)}$.

Equivalently, letting $E_j = \left\{\{b_{1j'},b_{kj'}\}: j' \in [j]\right\}$, $H_j\(\a_Z,2k\)$ is the graph obtained by taking $H\(\a_Z,2k\)/{E_j}$, deleting all multi-edges with even multiplicity, and then deleting all isolated vertices.

\begin{figure}[hbt!]
    \centering
    \begin{subfigure}[c]{.35\textwidth}
      \includegraphics[width=1\linewidth]{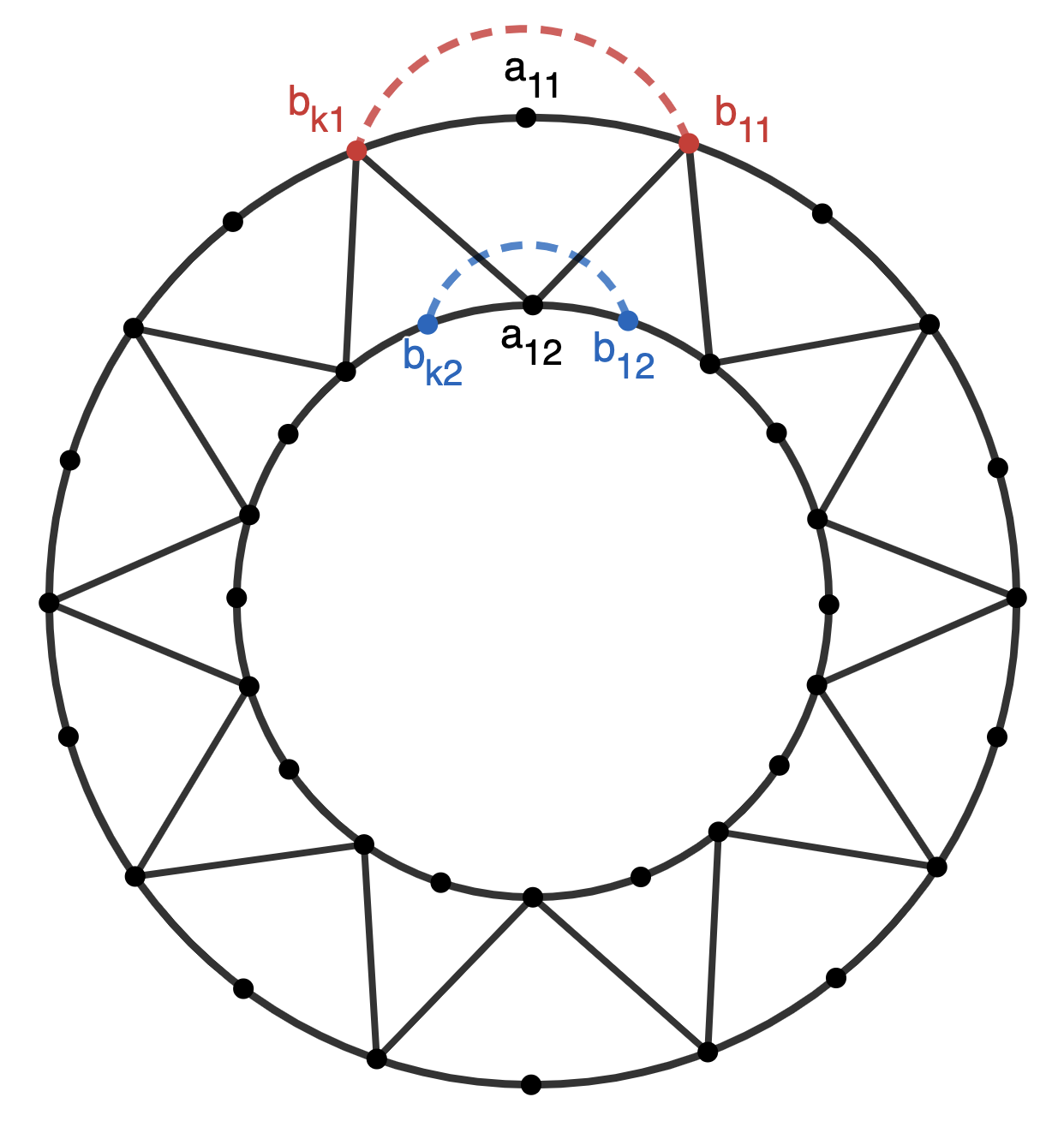}
    \end{subfigure}\hspace{1cm}
    \begin{subfigure}[c]{.4\textwidth}
      \includegraphics[width=1\linewidth]{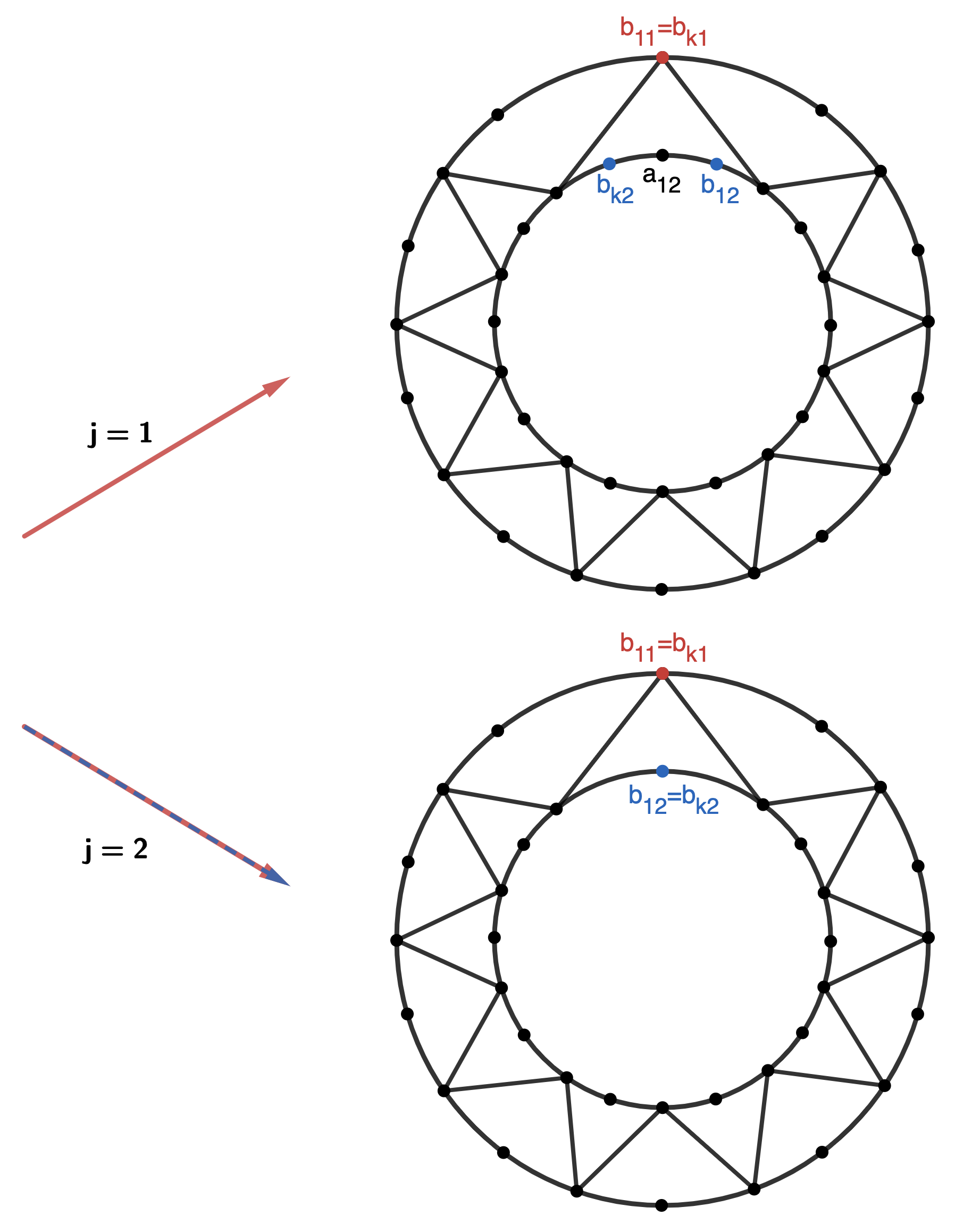}
    \end{subfigure}
    \caption{Illustration of \Cref{defn:H_j_Z}. Note that for $j = 1$, we have the red constraint edge between $b_{11}$ and $b_{k1}$ but we do not have the blue constraint edge beween $b_{12}$ and $b_{k2}$.}
    \label{fig:H_j_z}
\end{figure}

\end{defn}

\begin{defn}\label{defn:D-k-j}
For all $k \in \mathbb{N}$ and all $j \in \{0,1,2\}$, define $D(k,j)$ to be the number of dominant constraint graphs on $H_j\(\a_Z,2k\)$. Equivalently, $D(k,j)$ is the number of dominant constraint graphs $C \in \mathcal{C}_{\(\a_Z,2k\)}$ such that for all $j' \in [j]$, $b_{1j'} \=_{C} b_{kj'}$. We set $D(0,0) = 1$.

Note that $D(k,0) = \abs{C \in \mathcal{C}_{\(\a_Z,2k\)}: C \text{ is dominant}}$ and $D(k,2) = D(k-1,0)$.
\end{defn}

\begin{lemma}\label{lem:z-shape-split-1}
For all $k \in \mathbb{N}$ and $j \in \{0,1\}$, $\displaystyle D(k,j) = \sum_{i=0}^{k}{D(i,0)D(k-i,j+1)}$
\end{lemma}
\begin{proof}
Given a dominant constraint graph $C$ on $H_j\(\a_Z,2k\)$, let $i \in [k]$ be the last index such that $i=1$ or $a_{1(j+1)} \=_{C} a_{i(j+1)}$. 

We have that $i = 1$ if and only if $b_{1(j+1)} \=_{C} b_{k(j+1)}$. There are $D(k,j+1)$ dominant constraint graphs on $H_j\(\a_Z,2k\)$ for which this is true.

If $i > 1$ then we make the following observations:
\begin{enumerate}
\item By \Cref{lem:zshape-spoke-constraint}, for all $j' > j+1$, $a_{1j'} \=_{C} a_{ij'}$.
\item $b_{i(j+1)} \=_C b_{k(j+1)}$ as otherwise the edges $\left\{a_{i(j+1)},b_{i(j+1)}\right\}$ and $\left\{b_{k(j+1)},a_{1(j+1)}\right\}$ would only appear once in $H_j\(\a_Z,2k\)/C$. By \Cref{lem:zshape-spoke-constraint}, for all $j' \leq j$, $b_{ij'} \=_C b_{kj'}$.
\end{enumerate}

\begin{figure}[hbt!]
    \centering
    \begin{subfigure}[c]{.36\textwidth}
      \includegraphics[width=1\linewidth]{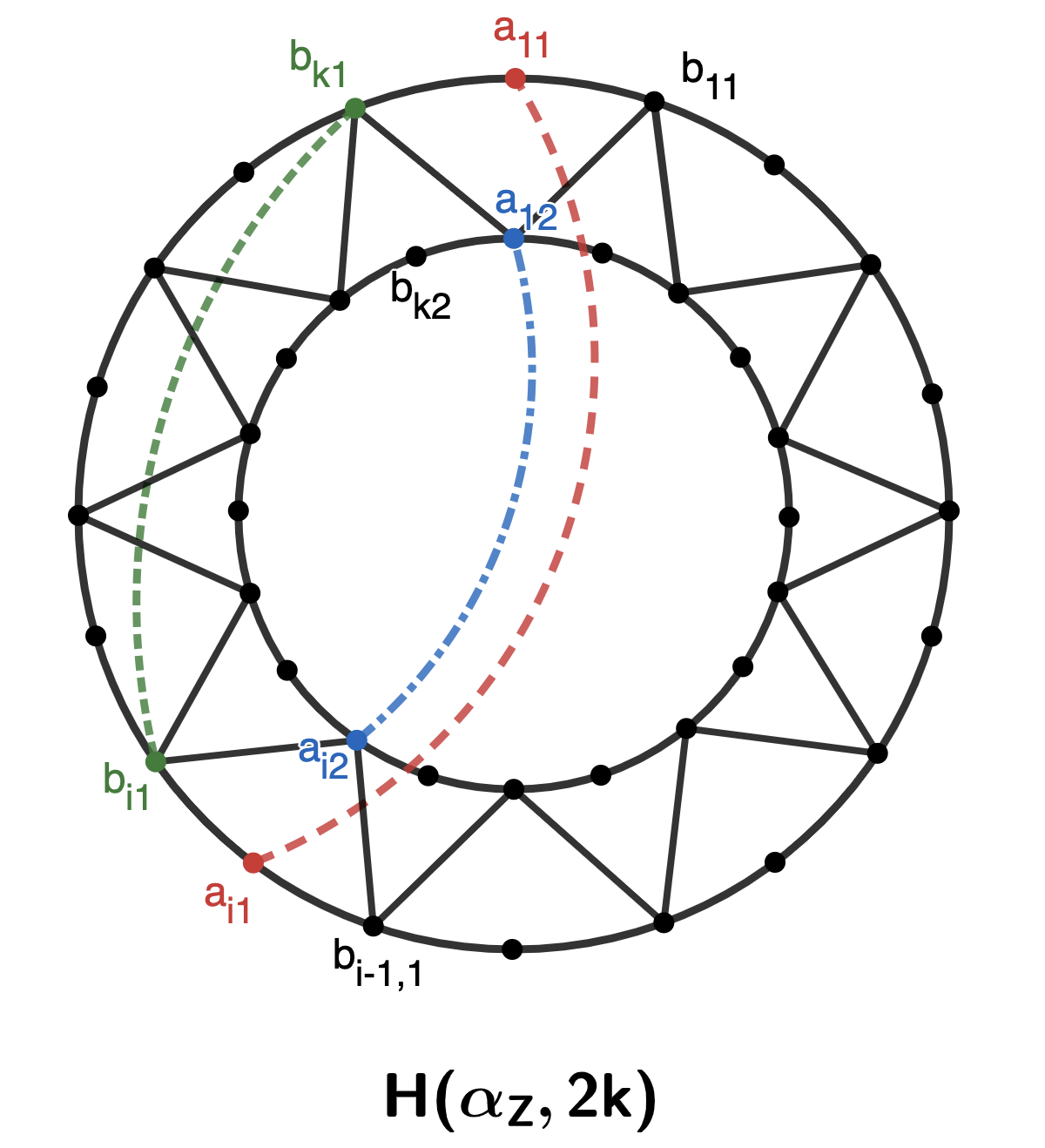}
    \end{subfigure}
    \begin{subfigure}[c]{.63\textwidth}
      \includegraphics[width=1\linewidth]{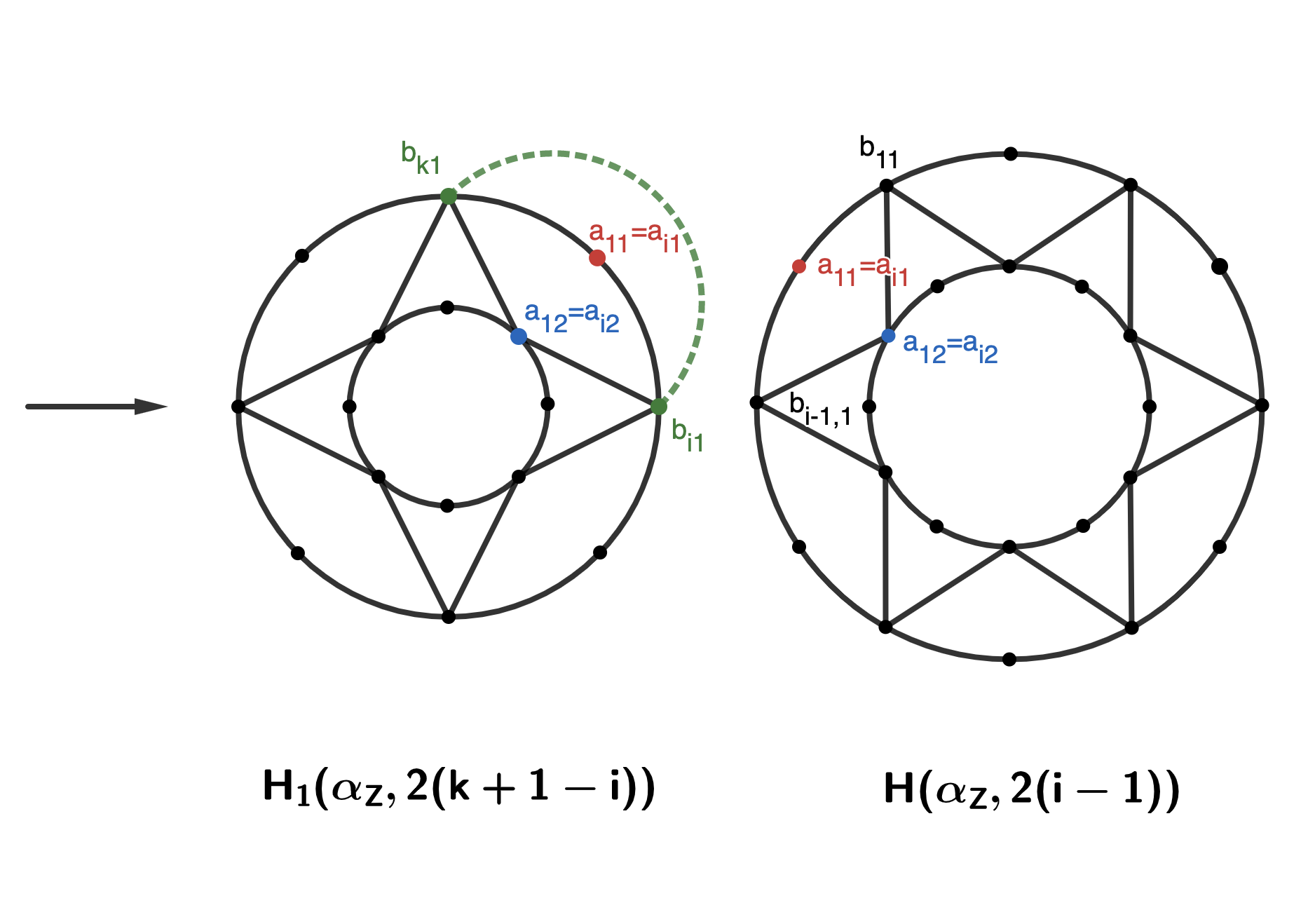}
    \end{subfigure}
    \caption{Illustration of \Cref{lem:z-shape-split-1} and \Cref{lem:z-shape-split-2}: case $j=0$}
    \label{fig:z-shape-split}
\end{figure}

\begin{lemma}\label{lem:z-shape-split-2}
$E = \left\{\{a_{1j'},a_{ij'}\}: j' \geq j+1\right\} \bigcup \left\{\{b_{ij'},b_{kj'}\}: j' \leq j+1\right\}$
partitions $H_j\(\a_Z,2k\)$ into two parts $H_1$ and $H_2$ where $H_1 \cong H\(\a_Z,2(i-1)\)$ and $H_2 \cong H_{j+1}\(\a_Z,2(k+1-i)\)$.
\end{lemma}
\begin{proof}
Let $m = 2$. We take $H_1$ to be the induced subgraph of $H_j\(\a_Z,2k\)$ on the vertices 
\begin{equation*}
V(H_1) = \left\{a_{i'j'}: i' \in [2,i], j' \in [m]\right\} \cup \left\{b_{i'j'}: i' \in [i-1], j' \in [m]\right\}
\end{equation*}
and we take $H_2$ to be the induced subgraph of $H_j\(\a_Z,2k\)$ on the vertices 
\begin{align*}
V(H_2) = &\left\{a_{i'j'}: i' \in [i+1,k], j' \in [m]\right\} \cup \left\{a_{ij'}: j' \in [j+2,m]\right\} \\
&\cup \left\{b_{i'j'}: i' \in [i+1,k], j' \in [m]\right\} \cup \left\{b_{ij'}: j' \in [j+2,m]\right\}
\end{align*}

To see that $H_2 \cong H_{j+1}\(\a_Z,2(k+1-i)\)$, note that by the above observations, for all $j' \leq j+1$, $b_{ij'} \=_C b_{kj'}$ and for all $j' \in [j+1,m]$, $a_{1j'} \=_C a_{ij'}$. Subtracting $i-1$ from all of the $i'$ indices of $H_2$ gives $H_{j+1}\(\a_Z,2(k+1-i)\)$.

To see that $H_1 \cong H\(\a_Z,2(i-1)\)$, observe that for all $j' \in [m]$, all edges and spokes in the range between $b_{1j'}$ and $a_{ij'}$ are present. For the edges $\left\{\{a_{ij'},b_{1j'}\}: j' \in [m]\right\}$ and the spokes $\left\{\{a_{i(j')}, b_{1(j'-1)}\}: j' \in [2,m]\right\}$, observe that

\begin{enumerate}
\item If $j' \leq j$ then the edge $\left\{a_{ij'},b_{1j'}\right\}$ comes from the edge between $a_{ij'}$ and $b_{ij'}$ as $b_{1j'} \=_C b_{kj'} \=_C b_{ij'}$. Similarly, if $j' \leq j+1$ then the spoke $\left\{a_{i(j')}, b_{1(j'-1)}\right\}$ comes from the spoke $\left\{a_{i(j')}, b_{i(j'-1)}\right\}$
\item If $j' \geq j+1$ then the edge $\left\{a_{ij'},b_{1j'}\right\}$ comes from the edge between $a_{1j'}$ and $b_{1j'}$ as $a_{1j'} \=_C a_{ij'}$. Similarly, if $j' \geq j+2$ then the spoke $\left\{a_{ij'},b_{1(j'-1)}\right\}$ comes from the spoke $\left\{a_{1(j')}, b_{1(j'-1)}\right\}$.
\end{enumerate}

\begin{rmk}
Note that the spoke $\left\{a_{i(j+1)}, b_{1(j)}\right\}$ cannot come from the spoke $\left\{a_{1(j+1)}, b_{1j}\right\}$ as this spoke was deleted in $H_j\(\a_Z,2k\)$.
\end{rmk}
Since all of the needed edges are present and no extra edges are present, $H_1 \cong H\(\a_Z,2(i-1)\)$.

Finally, we use \Cref{lem:Hpartitionconditions} to confirm that $E$ partitions $H_j\(\a_Z,2k\)$ into the two parts $H_1$ and $H_2$. The first condition of \Cref{defn:Hpartitions} and the first condition of \Cref{lem:Hpartitionconditions} can be seen directly. For the second and third conditions of \Cref{lem:Hpartitionconditions}, we observe that  $V(H_1) \cap V(H_2) = \left\{a_{1j'}: j' \in [j+2,m]\right\} \cup \left\{b_{1j'}: j' \in [j]\right\}$. Since all dominant constraint graphs on $H_j\(\a_Z,2k\)$, $H_1$, and $H_2$ are well-behaved, none of these constraint graphs can have a path between any two vertices in $V(H_1) \cap V(H_2) = \left\{a_{1j'}: j' \in [j+2,m]\right\} \cup \left\{b_{1j'}: j' \in [j]\right\}$.

To see that the fourth condition of \Cref{lem:Hpartitionconditions} holds, assume that there are edges $e_1,e_2 \in E(H)$ which appear with odd multiplicity in $H/{C_E}$ such that $e_1$ appears in $H_1$, $e_2$ appears in $H_2$, and $e_1 \=_{C} e_2$. Since $C$ is well-behaved, we must either have that $e_1$ and $e_2$ are edges in $W_i$ for some $i \in [m]$ or $e_1$ and $e_2$ are spokes between $W_i$ and $W_{i+1}$ for some $i \in [m-1]$. Using similar logic as we used to prove \Cref{lem:lineshapepartitioning}, we can rule out the possibility that $e_1$ and $e_2$ are both in $W_i$ for some $i \in [m]$. The possibility that $e_1$ and $e_2$ are spokes between wheels $W_i$ and $W_{i+1}$ for some $i \in [m-1]$ is ruled out by \Cref{lem:zshape-spoke-constraint} which implies that spokes in $H_1$ are not made equal to spokes in $H_2$.
\end{proof}
\end{proof}

\begin{cor}\label{cor:z-shape-recurrence-relation}
\begin{equation*}
    D(k,0) = \sum_{\substack{i_1, i_2, i_3\geq 0: \\ i_1+i_2+i_3 = k-1}} D\(i_1,0\)D\(i_2,0\)D\(i_3, 0\).
\end{equation*}
\end{cor}

\begin{proof}
By \Cref{lem:z-shape-split-1} and \Cref{defn:D-k-j},
    \begin{align*}
        D(k,0) 
        &= \sum_{i_1=0}^{k}D(i_1,0)D(k-i_1,1) \\
        &= \sum_{i_1=0}^{k}D(i_1,0)\(\sum_{i_2=0}^{k-i_1} D(i_2,0)D(k-i_1-i_2, 2)\) \\
        &= \sum_{i_1=0}^{k} D(i_1,0) \(\sum_{i_2=0}^{k-i_1} D(i_2,0)D\(k-i_1-i_2-1, 0\)\) \\
        &= \sum_{\substack{i_1, i_2, i_3\geq 0: \\ i_1+i_2+i_3 = k-1}} D\(i_1,0\)D\(i_2,0\)D\(i_3, 0\).
    \end{align*}
\end{proof}

By \Cref{cor:z-shape-recurrence-relation} and \Cref{thm:catalan2-recurrence-relation}, $D(k,0)$ has the same recurrence relation as $\displaystyle C'_k = \dfrac{1}{2k+1}\binom{3k}{k}$, thus \Cref{thm:num-zshape-constraint-graph} is proved.

\subsection{Analyzing \texorpdfstring{$\displaystyle \var \Big(\trace\(\(M_{\a_Z}M_{\a_Z}^T\)^k\)\Big)$}{Var(tr(Ma MaT)k)}}\label{subsection:z-shape-variance-trace-bound}
We now sketch how to show that $\displaystyle \var \Big(\trace\(\(M_{\a_Z}M_{\a_Z}^T\)^k\)\Big)$ is $O\(n^{4k+2}\)$.
    \begin{claim}
        For each $k$, $\displaystyle \var \Big(\trace\(\(M_{\a_Z}M_{\a_Z}^T\)^k\)\Big)$ is $O\(n^{4k+2}\)$.
    \end{claim}

    \begin{proof}
        Observe that 
        \small
        \begin{align*}
            &\var\(\trace\(\(M_{\alpha_Z}M_{\alpha_Z}^T\)^k\)\) 
            = \Eb\[\trace\(\(M_{\alpha_Z}M_{\alpha_Z}^T\)^k\)^2\] - \Eb\[\trace\(\(M_{\alpha_Z}M_{\alpha_Z}^T\)^k\)\]^2 \\
            =& \sum_{A_1,B_1,A'_1,B'_1,\ldots,A_k,B_k,A'_k,B'_k \subseteq \binom{[n]}{2}}\, \Bigg(\Eb\[\prod_{i=1}^{k}\,  M_{\alpha_Z}\(A_i,B_i\)M_{\alpha_Z}^T\(B_i,A_{i+1}\) M_{\alpha_Z}\(A_i',B_i'\) M_{\alpha_Z}^T\(B_i',A_{i+1}'\)\] \\
            &\hspace{3.5cm} - \Eb\[\prod_{i=1}^{k}\,  M_{\alpha_Z}\(A_i,B_i\)M_{\alpha_Z}^T\(B_i,A_{i+1}\)\] \cdot \Eb\[\prod_{i=1}^{k}\, M_{\alpha_Z}\(A_i',B_i'\) M_{\alpha_Z}^T\(B_i',A_{i+1}'\)\]\Bigg)
        \end{align*}
        \normalsize
        
    To analyze the terms in the summation, we consider two copies $H_1$ and $H_2$ of $\(\alpha_Z, 2k\)$ with vertices $\{A_i: i \in [k]\} \cup \{B_i: i \in [k]\}$ and $\{A'_i: i \in [k]\} \cup \{B'_i: i \in [k]\}$ respectively, where the constraint edges of a constraint graph $C$ can now have edges between $H_1$ and $H_2$ as well as edges within $H_1$ and edges within $H_2$. Denote \begin{align*}
    \displaystyle \val(C) &= \Eb\bigg[\prod_{i=1}^{k}\,  M_{\alpha_Z}\(A_i,B_i\)M_{\alpha_Z}^T\(B_i,A_{i+1}\) M_{\alpha_Z}\(A_i',B_i'\) M_{\alpha_Z}^T\(B_i',A_{i+1}'\)\bigg] - \\
    &\Eb\bigg[\prod_{i=1}^{k}\,  M_{\alpha_Z}\(A_i,B_i\)M_{\alpha_Z}^T\(B_i,A_{i+1}\)\bigg] \Eb\bigg[\prod_{i=1}^{k}\, M_{\alpha_Z}\(A_i',B_i'\) M_{\alpha_Z}^T\(B_i',A_{i+1}'\)\bigg]. 
    \end{align*}
    Observe that $\val(C) = 1$ if under $C$, each edge appears an even number of times in $H_1 \cup H_2$, but there is at least one edge that appears an odd number of times when restricted to only $H_1$ or $H_2$ (call this special edge the ``odd" edge). $\val(C) = 0$ in all other situations. 
    
    It suffices to prove that for the constraint graphs $C$ such that $\val(C) = 1$, $(H_1 \cup H_2)/C$ contains at most $4k+2$ distinct vertices as this implies that $\displaystyle \sum_{C} N(C)\val(C) = O_k\(n^{4k+2}\)$.
    
    In this proof sketch, we will only prove the above claim for well-behaved constraint graphs $C$. For the full analysis, see \Cref{thm:variancewellbehaved}. For well-behaved constraint graphs $C$, there are two cases:
    \begin{enumerate}
        \item The odd edge $e$ is among one of the wheels of $H_1$ or $H_2$. Let $W_{i,j}$ denote the $j^{th}$ wheel of $H_i$. Observe that for each $j=1,2$, $W_{1,j}\cup W_{2,j}$ has $4k$ edges, thus $(W_{1,j}\cup W_{2,j})/C$ has at most $2k$ multi-edges (since each edge has multiplicity at least $2$). W.L.O.G. assume $e$ is in $W_{1,1}$, the outer wheel of $H_1$. In this case, the following statements are true.
        \begin{enumerate}[i.]
            \item[1.] $(W_{1,1}\cup W_{2,1})/C$ is a connected graph: since $e$ is an odd edge, it is constrained to some edge in $W_{2,1}$.
            \item[2.] $(W_{1,1}\cup W_{2,1})/C$ has at least one cycle: each vertex in $W_{1,1}/C$ has an even degree and the multi-edge containing $e$ has an odd multiplicity, thus there must be at least one cycle. 
        \end{enumerate}
        This implies that $(W_{1,1}\cup W_{2,1})/C$ has at most $2k$ vertices. On the other hand, $(W_{1,2}\cup W_{2,2})/C$ has at most $2k+2$ vertices because it has at most two connected components. Thus, all together there are at most $2k + (2k+2) = 4k+2$ distinct vertices.
        
        \item The odd edge $e$ is among the spokes. Without loss of generality, assume $e\in H_1$. Then $e$ is constrained to some $e'$ in $H_2$. This implies that for each $j=1,2$, $(W_{1,j}\cup W_{2,j})/C$ is connected, and thus has at most $2k+1$ vertices. Thus, all together there are at most $2(2k+1) = 4k+2$ distinct vertices.
    \end{enumerate}
    \end{proof}

\section{The spectrum of the singular values of \texorpdfstring{$\frac{M_{\a_Z}}{n}$}{MaZ/n}}
\label{section:zshape-spectrum}

\setlength{\parskip}{1.5mm}
\setlength{\baselineskip}{1.3em}

In this section, we show how to use \Cref{thm:zshape-trace-convergence} to derive the limiting distribution of the singular values of $\frac{M_{\a_Z}}{n}$ as $n \to \infty$.
\subsection{The spectrum of Wigner matrices}

\setlength{\parskip}{1.5mm}
\setlength{\baselineskip}{1.3em}

Before we derive the limiting distribution of the spectrum of the singular values of  Z-shaped graph matrices, we will show a simpler but similar analysis for deriving Wigner's Semicircle Law from the trace power moments. 

\begin{defn}\label{defn:wigner-matrix}
    A \emph{Wigner matrix $M$} is a random Hermitian matrix where entries $M_{ij}$ for $i<j$ are i.i.d. complex random variables with mean $0$ and variance $1$, $M_{ij} = \overline{M_{ji}}$ for $i>j$, and the diagonal entries are i.i.d real random variables with bounded mean and variance. 
\end{defn}

\begin{defn}\label{defn:semicircle-law}
    The \emph{semi-circle law $F_{sc}(x)$} is a distribution function whose density function is 
    \begin{equation}
        f_{sc}(x) = F_{sc}'(x) = 
        \begin{cases}
            \;\frac{1}{2\pi} \sqrt{4-x^2} & \text{ if } |x|\leq 2\\
            \;0 &\text{ otherwise}
        \end{cases}
    \end{equation}
\end{defn}

\begin{thm}[Wigner's Semicircle Law]\label{thm:wigner}
    If $M_n$ is a sequence of $n\times n$ Wigner matrices whose entries have bounded moments then the distribution of eigenvalues of $\frac{M_n}{\sqrt{n}}$ weakly converges almost surely to $F_{sc}(x)$.
\end{thm}
The analysis in Sections \ref{sec:dominantcycleconstraintgraphs} and \ref{sec:cyclerecurrencerelation} for the $\pm{1}$ random matrix $M_{\a_0}$ can be easily generalized to Wigner matrices so we have the following theorem. 
\begin{thm}\label{thm:trace-power-wigner}
    If $M_n$ is a sequence of $n\times n$ Wigner matrices whose entries have bounded moments then $\lim_{n\to\infty} \dfrac{\,1\,}{n} \Eb\, \[\trace\(\(\frac{M_n}{\sqrt{n}}\)^{2k}\)\] = C_k $
    where $\displaystyle C_k = \dfrac{1}{k+1}\binom{2k}{k}$ is the $k^{th}$ Catalan number.
\end{thm}

One can verify that $\displaystyle \int_{-2}^{2}\, \dfrac{1}{2\pi}\sqrt{4-x^2}\cdot x^{2k} \,dx = C_k$ and conclude \Cref{thm:wigner} after checking the other conditions in \Cref{lem:moment-convergence}. Here, we will instead show how to derive $f_{sc}(x)$ by deriving a differential equation for $f_{sc}(x)$ from the recurrence relation for $C_k$. 

\begin{lemma}
    Assume $f(x)$ is a function satisfying that for all $k \in \mathbb{N}$,
    \begin{equation}
        \int_{-2}^{2} x^{2k}\cdot f(x)\,dx = C_k \label{eqn:moments-f-Ck}
    \end{equation}
    and moreover, 
    \begin{enumerate}
        \item $f(x)$ is continuously differentiable on $(-2,2)$.
        \item $\displaystyle \lim_{x\to 2^{-}} f(x) = \lim_{x\to -2^{+}} f(x) = 0$.
    \end{enumerate}
    Then $f(x)$ satisfies the following differential equation on $(-2,2)$:
    \begin{equation}
        (4-x^2)f'(x) + xf(x) = 0. \label{eqn:ODE-wigner}
    \end{equation}
\end{lemma}
\begin{proof}
We use the following recurrence relation for the Catalan numbers.
\begin{prop}
    \begin{equation}
        \dfrac{C_{k+1}}{C_{k}} = \dfrac{2(2k+1)}{k+2}. \label{eqn:rel-C_k} 
    \end{equation}
\end{prop}
\begin{proof}
\begin{align*}
    \dfrac{C_{k+1}}{C_{k}} = \dfrac{\dfrac{1}{k+2}\cdot \dfrac{(2k+2)!}{(k+1)!(k+1)!}}{\dfrac{1}{k+1}\cdot \dfrac{(2k)!}{k!k!}} = \dfrac{(k+1)\cdot (2k+2)(2k+1)}{(k+2)\cdot (k+1)^2} = \dfrac{2(2k+1)}{k+2}.
\end{align*}
\end{proof}

For any $k\geq 1$, using integration by parts and the boundary condition that $\displaystyle \lim_{x\to 2^{-}} f(x) = \lim_{x\to -2^{+}} f(x) = 0$, 
\begin{align*}
    \int_{-2}^2 x^{2k+1}\cdot f'(x) \,dx 
    &= \[x^{2k+1}\cdot f(x)\]_{-2}^2 - (2k+1)\cdot \int_{-2}^2 f(x)\cdot x^{2k} \,dx =  -(2k+1) C_k
\end{align*}

We have the following relations together with \Cref{eqn:rel-C_k}.
\begin{enumerate}[i.]
    \item $\displaystyle \int_{-2}^2 x^{2k}\cdot f(x) \,dx = C_k$
    \item $\displaystyle \int_{-2}^2 x^{2k+2}\cdot f(x) \,dx = C_{k+1}$
    \item $\displaystyle \int_{-2}^2 x^{2k+1}\cdot f'(x) \,dx = -(2k+1)C_k $
    \item $\displaystyle \int_{-2}^2 x^{2k+3}\cdot f'(x) \,dx = -(2k+3)C_{k+1}$
    \item $\displaystyle (k+2) C_{k+1} = 2(2k+1) C_k$
\end{enumerate}

Multiplying v. by $-2$ and rewriting using i. to iv., we get
\begin{align*}
    & -2(k+2)C_{k+1} = -(2k+3)C_{k+1} - C_{k+1} = -4(2k+1)C_k \\
    \implies & \int_{-2}^2 x^{2k+3}\cdot f'(x) \,dx - \int_{-2}^2 x^{2k+2}\cdot f(x) \,dx = 4 \int_{-2}^2 x^{2k+1}\cdot f'(x) \,dx\\
    \implies & \int_{-2}^2 x^{2k+1}\cdot \((4-x^2)f'(x) + xf(x)\) \, dx = 0.
\end{align*}

Using a similar argument as in \Cref{lem:zero even moments}, we can conclude that 
\begin{equation*}
    (4-x^2)f'(x) + xf(x) = 0.
\end{equation*}
\end{proof}

Finally we can solve the ODE for the solution of $f(x)$.
\begin{prop}
    Assume $f(x)$ satisfies \Cref{eqn:ODE-wigner} and $\displaystyle \int_{-2}^2 f(x)\,dx = 1 $. Then $f(x) = \dfrac{1}{2\pi}\sqrt{4-x^2}$ on $(-2,2)$.
\end{prop}

\begin{proof}
    \begin{align*}
        & (4-x^2)f'(x) + xf(x) = 0 
        \implies \dfrac{f'(x)}{f(x)} = -\dfrac{x}{4-x^2} \\
        \implies &  \int \dfrac{f'(x)}{f(x)} \, dx = \int -\dfrac{x}{4-x^2} \, dx
        \implies \ln\abs{f(x)} = \dfrac{\,1\,}{2}\ln\abs{4-x^2} + C \\
        \implies & f(x) = A\cdot \sqrt{4-x^2} \text{ for some constant } A \\
        \implies & f(x) = \dfrac{1}{2\pi} \sqrt{4-x^2} \text{ since } \int_{-2}^2 f(x)\,dx = 1 
    \end{align*}
\end{proof}

\subsection{The spectrum of the singular values of Z-shaped graph matrices}

We now find the limiting distribution of the spectrum of the singular values of $\frac{M_{\alpha_Z}}{n}$ as $n \to \infty$.

\begin{defn}
Let $a = \dfrac{3\sqrt{3}}{2}$ and let $g_{\alpha_Z}: (0,\infty) \to \mathbb{R}$ be the function such that 
\begin{equation*} 
    g_{\alpha_Z}(x) = \dfrac{i}{\pi}\cdot \(\sqrt{3}\cdot \sin\(\dfrac{1}{3}\cdot\arctan\(\dfrac{3}{\sqrt{4x^2/3-9}}\)\)+ \cos\(\dfrac{1}{3}\cdot\arctan\(\dfrac{3}{\sqrt{4x^2/3-9}}\)\)\)
\end{equation*}
if $x \in (0,a]$ and $g_{\alpha_Z}(x) = 0$ if $x > a$ where we take $\arctan(ix) = \dfrac{i}{2}\ln\(\dfrac{1+x}{1-x}\)$ for all real $x$, we take $ln(-x) = {\pi}i + ln(x)$ for all $x > 0$, and we take $sin(x) = \frac{e^{ix} - e^{-ix}}{2i}$ and $cos(x) = \frac{e^{ix} + e^{-ix}}{2}$ for all complex $x$.
\end{defn}

\begin{thm}\label{Zspectrum}
With probability $1$, the spectrum of the singular values of $\frac{1}{n}M_{\alpha_Z}$ approaches $g_{\alpha_Z}$ weakly as $n\to\infty$.
\end{thm}

\begin{proof}
    By \Cref{thm:zshape-trace-convergence}, it suffices to prove that for all $k \in \mathbb{N}$,
    \begin{equation*}
        \Eb_{x\sim g_{\a_Z}(x)} \big[x^{2k}\big] = \int_{x = 0}^{a}{x^{2k}g_{\a_{Z}}(x) \,dx} = C'_k.
    \end{equation*}

    To prove that $\displaystyle\int_{x = 0}^{a}{x^{2k}g_{\a_{Z}}(x) \,dx} = C'_k$, we proceed as follows:
\begin{enumerate}
    \item We derive a differential equation for $g_{\a_{Z}}$ based on a recurrence relation for $C'_k$ (see Theorem \ref{thm:ODE for f(x) of Z-shape}).
    \item We prove that if $g_{\a_{Z}}$ satisfies this differential equation and some conditions at $x = 0$ and $x = a$ then $\displaystyle\int_{x = 0}^{a}{x^{2k}g_{\a_{Z}}(x)dx} = C'_k$ (see Theorem \ref{thm:verify ODE for Z-shape spec}).
    \item We verify that $g_{\a_{Z}}$ satisfies the required conditions (see Theorem \ref{thm:solution to ODE for Z-shape}).
\end{enumerate}

\begin{rmk}
Technically, only steps 2 and 3 are needed. We include the first step because it gives better intuition for where the differential equation comes form.
\end{rmk}

\begin{thm}\label{thm:ODE for f(x) of Z-shape}
    Let $\displaystyle C'_k = \dfrac{1}{2k+1}\binom{3k}{k}$ and $a=\lim_{k\to\infty} C'_{k+1}/C'_k = 3\sqrt{3}/2$. Assume $f(x)$ is an function satisfying that for all $k \in \mathbb{N}$,
    \begin{equation}
        \int_{0}^{a} x^{2k}\cdot f(x)\dx = C'_k \label{eqn:moments-f-Dk}
    \end{equation}
    and moreover, 
    \begin{enumerate}
        \item $f(x)$ is twice continuously differentiable on $(0,a)$.
        \item $\lim_{x\to 0^{+}}xf(x)=0$ and $\lim_{x\to 0^{+}}x^2f'(x)=0$.
        \item $\lim_{x\to a^{-}}f(x)=0$ and $\lim_{x\to a^{-}}f'(x)(4x^2-27)=\lim_{x \to a^{-}}8af'(x)(x-a)=0$.
        \item $\lim_{x \to 0^{+}}{x^3f''(x)} = 0$ and $\lim_{x\to a^{-}}(a-x)^{2}f''(x) = 0$.
    \end{enumerate}
    Then $f(x)$ satisfies the following differential equation on $(0,a)$:
    \begin{equation}
        (4x^4-27x^2)f''(x)+(8x^3-27x)f'(x)+3f(x)=0. \label{eqn:ODE for f(x) of Z-shape}
    \end{equation}
\end{thm}

\begin{proof}
To prove this, we use the following recurrence relation for $\displaystyle C'_k = \dfrac{1}{2k+1}\binom{3k}{k}$.
\begin{prop}\label{prop:Drecurrence}
\begin{equation}
\dfrac{{C'_k}}{{C'_{k-1}}} = \dfrac{3(3k-1)(3k-2)}{2k(2k+1)}. \label{eqn:rel between D_k} 
\end{equation}
\end{prop}

\begin{proof}
Observe that
\begin{align*}
\dfrac{C'_k}{C'_{k-1}} &= \frac{2k-1}{2k+1} \cdot \frac{(3k)!(2k-2)!(k-1)!}{(3k-3)!(2k)!k!} \\
&= \frac{2k-1}{2k+1} \cdot \frac{(3k)(3k-1)(3k-2)}{(2k)(2k-1)k} = \dfrac{3(3k-1)(3k-2)}{2k(2k+1)}.
\end{align*}
\end{proof}

We also need the following relationship between the moments of $f$ and the moments of its derivatives.
\begin{defn}
    For all $j \in \{0,1,2\}$ and $k \in \mathbb{Z}$ such that $k \geq j$, we define $A(j,k)$ to be $\displaystyle A(j,k):=\int_{0}^{a} f^{(j)}(x)\cdot x^{k}\dx$ where $f^{(j)}(x)$ denotes the $j^{th}$ derivative of $f$. Notice that $A(0,2k) = C'_k$.
\end{defn}
\begin{lemma}\label{lemma:relationsbetweenA}
For all $j \in \{1,2\}$ and $k \in \mathbb{Z}$ such that $k \geq j$, 
\begin{equation*}
A(j,k) = \[f^{(j-1)}(x) \cdot x^k\]_0^{a} - kA(j-1,k-1).
\end{equation*}
\end{lemma}
\begin{proof}
Using integration by parts, we have that 
 \begin{align*}
         A(j,k) = \int_{0}^{a}{f^{(j)}(x)\cdot x^{k}\dx} &= \[f^{(j-1)}(x) \cdot x^k\]_0^{a} - \int_{0}^{a} kf^{(j-1)}(x)\cdot x^{k-1}\dx \\&
         = \[f^{(j-1)}(x) \cdot x^k\]_0^{a} - kA(j-1,k-1).
\end{align*}
\end{proof}

\begin{cor}\label{cor:relationsbetweenA}
If $\lim_{x\to 0^{+}}xf(x)=0$ and $\lim_{x\to a^{-}}f(x)=0$ then
\begin{enumerate}
\item For all $k \in \mathbb{N}$, $A(1,k) = kA(0,k-1)$
\item For all $k \in \mathbb{N}$ such that $k \geq 2$, 
\begin{equation*}
A(2,k) = \[f'(x) \cdot x^k\]_0^{a} - kA(1,k-1) = \[f'(x) \cdot x^k\]_0^{a} + k(k-1)A(0,k-2).
\end{equation*}
\end{enumerate}
\end{cor}

    Using Corollary \ref{cor:relationsbetweenA}, Proposition \ref{prop:Drecurrence}, and the fact that $A(0,2k) = {C'_k}$, we have that for all $k \in \mathbb{N}$,
    \begin{align*}
        A(2,2k+2) 
        &= \[f'(x)\cdot x^{2k+2}\]_{0}^{a} - 2kA(1,2k+1) - 2A(1,2k+1)\\
        &= \[f'(x)\cdot x^{2k+2}\]_{0}^{a} + (2k)(2k+1)A(0,2k) - 2A(1,2k+1)\\
        & = \[f'(x)\cdot x^{2k+2}\]_{0}^{a} + 3(3k-1)(3k-2)A(0,2k-2) - 2A(1,2k+1).
    \end{align*}
    
    Multiplying both sides by $4$ and repeatedly applying Corollary \ref{cor:relationsbetweenA}, we get that
    \begin{align*}
        &4A(2,2k+2) \\
        & = 4\[f'(x)\cdot x^{2k+2}\]_{0}^{a} + 27(2k)(2k-1)A(0,2k-2) + (-54k+24)A(0,2k-2) - 8A(1,2k+1)\\
        & = 4\[f'(x)\cdot x^{2k+2}\]_{0}^{a} + 27\cdot\(-\[f'(x)\cdot x^{2k}\]_{0}^{a} + A(2,2k)\) - 27(2k-1)A(0,2k-2)\\
        &- 3A(0,2k-2) - 8A(1,2k+1)\\
        & = \[x^{2k}f'(x)\cdot\(4x^2-27\)\]_0^{a} + 27A(2,2k) + 27A(1,2k-1) - 3A(0,2k-2) - 8A(1,2k+1) \\
        & = 27A(2,2k) + 27A(1,2k-1) - 3A(0,2k-2) - 8A(1,2k+1).
    \end{align*}
    where the last equality holds because $\lim_{x \to a^{-}}(4x^2-27)f'(x)=0$ and $\lim_{x \to 0^{+}}x^2f'(x)=0$ by assumption. 
    
    Writing the $A(j,k)$'s above as integrals, we get that for all $k \in \mathbb{N}$
    \begin{align*}
        \int_{0}^{a} \(4f''(x)\cdot x^4 - 27 f''(x)\cdot x^{2} - 27 f'(x) \cdot x + 8 f'(x) \cdot x^{3} + 3 f(x) \)\cdot x^{2k-2} \dx = 0.
    \end{align*}
    
    One way for this equation to be true is if $(4x^4-27x^2)f''(x)+(8x^3-27x)f'(x)+3f(x)=0$ on $(0,a)$. As shown by the following lemma and corollary, this is the only way for this equation to be true for all $k \in \mathbb{Z}$, which completes the proof of Theorem \ref{thm:ODE for f(x) of Z-shape}.
    
\begin{lemma}\label{lem:zero even moments}
    Let $a$ be some positive constant. If $f$ is continuous on $[0,a]$ and $\displaystyle\int_{0}^{a} f(x)x^{2k}\dx=0$ for all nonnegative integers $k$, then $f=0$ on $(0,a)$.
\end{lemma}\label{lem:zeromomentsimplyzero}
\begin{proof}
    Let $M>0$ be an upper bound of $f$ on $[0,a]$. For an arbitrary $\epsilon>0$, let $p(x)$ be a polynomial such that $\abs{p(x)-f(\sqrt{x})}<\dfrac{\epsilon}{M\cdot a}$ for all $x\in(0,a^2)$. Taking $p_\epsilon(x)=p(x^2)$,  $p_\epsilon$ is a linear combination of monomials of even power and $\abs{p_\epsilon(x)-f(x)}<\dfrac{\epsilon}{M\cdot a}$ for all $x\in(0,a)$. Thus $\displaystyle\int_{0}^a \(f(x)-p_\epsilon(x)\)\cdot f(x)\dx <\epsilon$. On the other hand, since all even moments of $f(x)$ are zero, 
    \begin{equation*}
        \int_{0}^a \(f(x)-p_\epsilon(x)\)\cdot f(x)\dx
        = \int_{0}^a f(x)^2 - p_\epsilon(x)f(x) \dx = \int_{0}^a f(x)^2\dx.
    \end{equation*}
    
    Thus $\displaystyle\int_0^a f(x)^2\dx < \epsilon$ for all $\epsilon>0$ and we conclude that $f(x)=0$ on $(0,a)$. 
\end{proof}
\begin{cor}
Let $a$ be some positive constant. If $f$ is continuous on $(0,a)$, $\lim_{x\to 0^{+}}x^{2}f(x)=0$, $\lim_{x \to a^{-}}{(a-x)^2f(x)} = 0$, and $\displaystyle\int_{0}^{a} f(x)x^{2k}\dx=0$ for all nonnegative integers $k$, then $f=0$ on $(0,a)$.
\end{cor}
\begin{proof}
This follows by applying Lemma \ref{lem:zero even moments} to the function $f(x)x^2(a - x)^2$.
\end{proof}
\end{proof}
We now confirm that if $f$ satisfies the differential equation $(4x^4-27x^2)f''(x)+(8x^3-27x)f'(x)+3f(x)=0$, the conditions of Theorem \ref{thm:ODE for f(x) of Z-shape}, and the condition that $\displaystyle\int_{0}^a f(x)\dx =1$, then $\displaystyle\int_{0}^{a} x^{2k}\cdot f(x)\dx = {C'_k}$.
\begin{thm}\label{thm:verify ODE for Z-shape spec}
    Let $a=\lim_{k\to\infty} {C'_k}/{C'_{k-1}} = 3\sqrt{3}/2$. Let $f$ be a function satisfying the following ODE on $(0,a)$
    \begin{equation}
        (4x^4-27x^2)f''(x)+(8x^3-27x)f'(x)+3f(x)=0
        \label{eqn: ODE for spectrum of Z-shape}
    \end{equation}
    and the first three conditions in Theorem \ref{thm:ODE for f(x) of Z-shape}, i.e.
    \begin{enumerate}
        \item $f(x)$ is twice continuously differentiable on $(0,a)$.
        \item $\lim_{x\to 0^{+}}xf(x)=0$ and $\lim_{x\to 0^{+}}x^2f'(x)=0$.
        \item $\lim_{x\to a^{-}}f(x)=0$ and $\lim_{x\to a^{-}}f'(x)(4x^2-27)=\lim_{x \to a^{-}}8af'(x)(x-a)=0$.
    \end{enumerate} 
    Moreover, assume that $\displaystyle\int_{0}^a f(x)\dx =1$. Then for all $k \in \mathbb{N} \cup \{0\}$,
    \begin{equation}
        \int_{0}^{a} x^{2k}\cdot f(x)\dx = C'_k \,.
    \end{equation}
\end{thm}

\begin{proof}
    Notice that $\displaystyle\int_{0}^{a}f(x)\dx=1 = {C_0}'$ by assumption. We aim to prove that for all $k \in \mathbb{N} \cup \{0\}$, 
 \begin{equation*}
 (2k+3)(2k+2)\int_{0}^{a} x^{2k+2}\cdot f(x)\dx = 3(3k+2)(3k+1)\int_{0}^{a} x^{2k}\cdot f(x)\dx.
 \end{equation*}
 If so, then since $(2k+3)(2k+2){C'_{k+1}} = 3(3k+2)(3k+1){C'_{k}}$, we can prove Theorem \ref{thm:verify ODE for Z-shape spec} by induction on $k$. 
    
    We multiply \eqref{eqn: ODE for spectrum of Z-shape} by $x^{2k}$ and integrate from $0$ to $a$:
    \begin{align*}
        0 
        &=\int_{0}^{a}(4x^4-27x^2)f''(x)\cdot x^{2k} + (8x^3-27x)f'(x)\cdot x^{2k} + 3x^{2k}f(x) \dx \\
        & = \(\[f'(x)(4x^4-27x^2)x^{2k}\]_{0}^a - \int_{0}^a f'(x) \(4(2k+4)x^{2k+3} - 27(2k+2)x^{2k+1}\) \dx \) \\
        &\hspace{2cm}+ \int_{0}^a (8x^{2k+3}-27x^{2k+1})f'(x)\dx + \int_{0}^a 3x^{2k}f(x)\dx \\
        & = -\int_{0}^a \(8(k+1)x^{2k+3} - 27(2k+1)x^{2k+1}\)f'(x)\dx + \int_{0}^a 3x^{2k}f(x)\dx\\
        & = -\[f(x)\(8(k+1)x^{2k+3} - 27(2k+1)x^{2k+1}\)\]_0^a \\
        &\hspace{2cm}+ \int_{0}^a \(8(k+1)(2k+3)x^{2k+2} - 27(2k+1)^2 x^{2k}\)f(x)\dx +  \int_{0}^a 3x^{2k}f(x)\dx \\
        & = \int_{0}^a \(8(k+1)(2k+3)x^{2k+2} - 3(36k^2+36k+8) x^{2k}\)f(x)\dx\\
        & = \int_{0}^a \(8(k+1)(2k+3)x^{2k+2} - 12(3k+1)(3k+2) x^{2k}\)f(x)\dx 
    \end{align*}
    as $\displaystyle \[f'(x)(4x^4-27x^2)x^{2k}\]_{0}^a$ and $\displaystyle \[f(x)\(8(k+1)x^{2k+3} - 27(2k+1)x^{2k+1}\)\]_0^a$ are zero by the assumed conditions on $f$.
    
    Rearranging the last step we get
    \begin{equation*}
        (2k+2)(2k+3)\int_0^a f(x)x^{2k+2} \dx = 3(3k+1)(3k+2)\int_{0}^a f(x)x^{2k}\dx
    \end{equation*}
    as needed. 
\end{proof}
Using WolframAlpha to solve the above ODE and analyzing the constant coefficient by the imposed initial conditions, we can get an explicit solution for $f(x)$. We verify the solution below.

\begin{thm}\label{thm:solution to ODE for Z-shape}
    Let $a=3\sqrt{3}/2$ and $f: (0,a) \to \mathbb{R}$ be the function such that
    \begin{equation}
        f(x)= \dfrac{i}{\pi}\cdot \(\sqrt{3}\cdot \sin\(\dfrac{1}{3}\cdot\arctan\(\dfrac{3}{\sqrt{4x^2/3-9}}\)\)+ \cos\(\dfrac{1}{3}\cdot\arctan\(\dfrac{3}{\sqrt{4x^2/3-9}}\)\)\) 
        \label{eqn:explicit soln to f(x) for Z-shape}
    \end{equation}
    where we take $\arctan(ix) = \dfrac{i}{2}\ln\(\dfrac{1+x}{1-x}\)$ for all real $x$, we take $ln(-x) = {\pi}i + ln(x)$ for all $x > 0$, and we take $sin(x) = \frac{e^{ix} - e^{-ix}}{2i}$ and $cos(x) = \frac{e^{ix} + e^{-ix}}{2}$ for all complex $x$. Then $f(x)$ is an solution to the ODE \eqref{eqn:ODE for f(x) of Z-shape} on $(0,a)$. Moreover, $f$ satisfies the conditions listed in Theorem \ref{thm:verify ODE for Z-shape spec}. 
\end{thm}

\begin{proof}
We first verify that this $f(x)$ satisfies the ODE \eqref{eqn:ODE for f(x) of Z-shape}
\begin{equation*}
    (4x^4-27x^2)f''(x)+(8x^3-27x)f'(x)+3f(x)=0.
\end{equation*}
on $(0,a)$.

For simplicity, let $\displaystyle g(x)=\dfrac{1}{3}\cdot\arctan\(\dfrac{3}{\sqrt{4x^2/3-9}}\)$. Then 
\begin{enumerate}

\item $\displaystyle f(x) = \dfrac{i}{\pi}\(\sqrt{3}\sin(g(x))+\cos(g(x))\)$.
\item $f'(x)=\dfrac{i}{\pi}\(\sqrt{3}\cos(g(x))-\sin(g(x))\)\cdot g'(x) = \dfrac{i}{\pi}\(\sqrt{3}\cos(g(x))-\sin(g(x))\)\dfrac{-1}{x\sqrt{4x^2/3-9}}$.
\item
\begin{align*}
    f''(x)
    &= \dfrac{i}{\pi}\(-\sqrt{3}\sin g(x)-\cos g(x)\)\cdot (g'(x))^2 + \dfrac{i}{\pi}\(\sqrt{3}\cos g(x)-\sin g(x)\)\cdot g''(x)\\
    & = \dfrac{-i}{\pi}\(\sqrt{3}\sin g(x)+\cos g(x)\)\cdot \dfrac{1}{x^2(4x^2/3-9)} \\
    &+\dfrac{i}{\pi} \(\sqrt{3}\cos g(x)-\sin g(x)\)\cdot\dfrac{8x^2/3-9}{x^2(4x^2/3-9)^{3/2}}.
\end{align*}

\end{enumerate}
Plugging the above into the LHS of \eqref{eqn:ODE for f(x) of Z-shape}, one can verify that $(4x^4-27)f''(x)+(8x^3-27x)f'(x)+3f(x)=0$.

We now check the conditions listed in Theorem \ref{thm:verify ODE for Z-shape spec}. For this purpose, it is more convenient to re-write $f(x)$ as a function of all real terms. Let $\displaystyle y = \dfrac{3}{\sqrt{9-4x^2/3}}$ and let $z=\dfrac{y-1}{y+1}=\dfrac{27-2x^2-9\sqrt{9-4x^2/3}}{2x^2}$. Note that when $0 < x < \frac{3\sqrt{3}}{2}$, $y$ and $z$ are real variables, $y \geq 1$, and $0 \leq z \leq 1$. We make the following observations:
\begin{enumerate}
\item $\displaystyle g(x)=\dfrac{1}{3}\cdot\arctan(-iy)=\dfrac{i}{6}\ln\(\dfrac{1-y}{1+y}\) = \dfrac{i}{6}ln(-z) = \dfrac{i}{6}ln(z) - \frac{\pi}{6}$
\item $sin(g(x)) = \frac{i}{2}\left(e^{-i(iln(z)/6 - {\pi}/6)} - e^{i(iln(z)/6 - {\pi}/6)}\right) = \frac{i}{2}\left(\left(\frac{\sqrt{3} + i}{2}\right)z^{1/6} - \left(\frac{\sqrt{3} - i}{2}\right)z^{-1/6}\right)$
\item $cos(g(x)) = \frac{1}{2}\left(e^{-i(iln(z)/6 - {\pi}/6)} + e^{i(iln(z)/6 - {\pi}/6)}\right) = \frac{1}{2}\left(\left(\frac{\sqrt{3} + i}{2}\right)z^{1/6} + \left(\frac{\sqrt{3} - i}{2}\right)z^{-1/6}\right)$
\end{enumerate}
Plugging in the above equations to $f(x)$ and simplifying, we get that
    \begin{align*}
        &f(x)=\dfrac{i}{\pi}\(\sqrt{3}\sin(g(x))+\cos(g(x))\)=\dfrac{-1}{\pi}\cdot\(z^{1/6}-z^{-1/6}\), \\
        &f'(x)=\dfrac{i}{\pi}\(\sqrt{3}\cos(g(x))-\sin(g(x))\)\cdot\dfrac{-1}{x\sqrt{4x^2/3-9}} = \dfrac{-1}{\pi}\cdot\(z^{1/6}+z^{-1/6}\)\cdot\dfrac{1}{x\sqrt{9-4x^2/3}}.
    \end{align*}
    
Recall that $y = \dfrac{3}{\sqrt{9-4x^2/3}}$ and $z=\dfrac{y-1}{y+1}$. Observe that 
\begin{enumerate}
\item As $x \to 0^{+}$, $y \approx 1 + \dfrac{2x^2}{27}$. Thus, $\lim_{x \to 0^{+}}{\dfrac{z}{x^2}} = \dfrac{1}{27}$.
\item As $x \to a^{-}$, $y \to \infty$. Thus, $\lim_{x \to a^{-}}{z} = 1$.
\end{enumerate}
Thus,
\begin{enumerate}
    \item $f$ is twice differentiable on $(0,a)$.
    \item $\lim_{x\to 0^{+}}xf(x)= \lim_{x\to 0}x\cdot\(\dfrac{z^{-1/6}}{\pi}\) = 0$.
    \item $\lim_{x\to 0^{+}}x^{2}f'(x)= \lim_{x\to 0}x\cdot\(\dfrac{-z^{-1/6}}{3\pi}\) = 0$.
    \item $\lim_{x\to a^{-}} f(x)=\dfrac{-1}{\pi}(1-1)=0$.
    \item $\lim_{x\to a^{-}}(4x^2-27)f'(x)=\lim_{x\to a^{-}}\dfrac{1}{\pi}(z^{1/6}+z^{-1/6})\cdot\(\dfrac{\sqrt{3(27-4x^2)}}{x}\) = 0$.
\end{enumerate}
    
We now prove the last piece of this Theorem:  $\displaystyle\int_{0}^a f(x)\dx=1$. We have that $a=3\sqrt{3}/2$, 
$z=\dfrac{y-1}{y+1}=\dfrac{27-2x^2-9\(9-4x^2/3\)^{1/2}}{2x^2}=\dfrac{27-27\(1-x^2/a^2\)^{1/2}}{2x^2}-1$, and $\displaystyle f(x)=\dfrac{-1}{\pi}\cdot\(z^{1/6}-z^{-1/6}\)$.
    
Let $x=a\sin\theta$. Then $z=\dfrac{27-27\cos\theta}{2a^2\sin^2\theta}-1=\dfrac{2(1-\cos\theta)}{\sin^2\theta}-1=\dfrac{1-\cos\theta}{1+\cos\theta}$. Expressing $\cos\theta$ in terms of $z$, we get $\cos\theta = \dfrac{1-z}{1+z}$, thus $\sin\theta=\dfrac{2\sqrt{z}}{1+z}$. Moreover,
    \begin{align*}
        \dz 
        & = \(\dfrac{1-\cos\theta}{1+\cos\theta}\)'d\theta 
        = \dfrac{2\sin\theta}{(1+\cos\theta)^2}d\theta 
        = \dfrac{2\sin\theta(1-\cos\theta)}{\sin^2\theta(1+\cos\theta)}
        = \dfrac{2z}{\sin\theta}d\theta
        \implies d\theta = \dfrac{\sqrt{z}}{z(1+z)}\dz.
    \end{align*}
    
    Thus 
    \begin{align*}
        \int_{0}^a f(x)\dx
        &= \dfrac{-1}{\pi}\int_{0}^{\pi/2} \(\(\dfrac{1-\cos\theta}{1+\cos\theta}\)^{1/6} - \(\dfrac{1-\cos\theta}{1+\cos\theta}\)^{-1/6}\)a\cos\theta d\theta \\
        &= \dfrac{-a}{\pi}\int_{0}^{1} \(z^{1/6}-z^{-1/6}\)\cdot\(\dfrac{1-z}{1+z}\)\cdot\dfrac{\sqrt{z}}{z(1+z)} \dz\\
        &= \dfrac{-a}{\pi}\int_{0}^1 \dfrac{(1-z)(z^{2/3}-z^{1/3})}{z(1+z)^2}\dz.
    \end{align*}
    
    Let $w=z^3$, then
    \begin{align*}
        \int_{0}^a f(x)\dx
        &=\dfrac{-a}{\pi}\int_{0}^1 \dfrac{(1-w^3)(w-1)}{(1+w^3)^2}\dw \\
        &= \dfrac{-a}{\pi}\int_{0}^1 \dfrac{-4/3}{(1+w)^2} + \dfrac{2w}{(w^2-w+1)^2} +\dfrac{-5/3}{w^2-w+1}\dw \\
        &= \dfrac{-a}{\pi}\(\dfrac{4}{3}\[\dfrac{1}{1+w}\]_0^1 + \int_{0}^1 \dfrac{2w-1}{(w^2-w+1)^2}\dw + \int_{0}^1 \dfrac{1}{(w^2-w+1)^2}\dw + \int_0^1\dfrac{-5/3}{w^2-w+1}\dw\) \\
        &= \dfrac{-a}{\pi}\(-\dfrac{2}{3} + \[\dfrac{-1}{(w^2-w+1)}\]_0^1 +  \int_{0}^1 \dfrac{1}{(w^2-w+1)^2}\dw + \int_0^1\dfrac{-5/3}{w^2-w+1}\dw\) \\
        &= \dfrac{-a}{\pi}\(-\dfrac{2}{3} + \int_{0}^1 \dfrac{1}{\((w-\frac{1}{2})^2+\frac{3}{4}\)^2}\dw + \int_0^1\dfrac{-5/3}{w^2-w+1}\dw\). 
    \end{align*}
    
\begin{lemma}
For any $b\neq 0$,
    \begin{equation}
        \int\dfrac{1}{(x^2+b^2)^2}\dx = \dfrac{1}{2b^2}\(\int \dfrac{1}{x^2+b^2} \dx - \dfrac{x}{x^2+b^2}\).
    \end{equation}
\end{lemma}
\begin{proof}
    \begin{align*}
        \int\dfrac{1}{(x^2+b^2)^2}\dx
        & = \dfrac{1}{b^2}\int\dfrac{x^2+b^2}{(x^2+b^2)^2} + \dfrac{-x^2}{x^2+b^2}\dx \\
        & = \dfrac{1}{b^2}\(\int\dfrac{1}{x^2+b^2}\dx + \int \dfrac{-x}{2}d\(\dfrac{1}{x^2+b^2}\)\)\\
        & = \dfrac{1}{b^2}\(\int\dfrac{1}{x^2+b^2} \dx - \dfrac{x}{2(x^2+b^2)} + \int\dfrac{-\frac{1}{2}}{x^2+b^2}\dx\) \\
        & = \dfrac{1}{2b^2}\(\int\dfrac{1}{x^2+b^2} \dx - \dfrac{x}{x^2+b^2}\).
    \end{align*}
\end{proof}

Apply the lemma to the $\displaystyle\int_{0}^1 \dfrac{1}{\((w-1/2)^2+3/4\)^2}\dw$ term, we get that
\begin{align*}
    \int_{0}^a f(x)\dx
    &= \dfrac{-a}{\pi}\(-\dfrac{2}{3} +\dfrac{2}{3}\(\int_{0}^1 \dfrac{1}{(w-1/2)^2+3/4}\dw -\[\dfrac{w-1/2}{w^2-w+1}\]_{0}^1 \) + \int_0^1\dfrac{-5/3}{(w-1/2)^2+3/4} \dw\)\\
    &= \dfrac{-a}{\pi}\(-\dfrac{2}{3} +\dfrac{2}{3} + \int_0^1\dfrac{-1}{(w-1/2)^2+3/4} \dw\) \\
    &=\dfrac{a}{\pi}\[\dfrac{1}{\sqrt{3}/2}\arctan\(\dfrac{w-1/2}{\sqrt{3}/2}\)\]_{0}^1 
    = \dfrac{3\sqrt{3}/2}{\pi}\cdot\dfrac{2\pi}{3\sqrt{3}} 
    = 1.
\end{align*}
\end{proof}

\end{proof}

\begin{figure}[hbt!]
    \centering
    \begin{subfigure}{.49\textwidth}
      \includegraphics[width=1\linewidth]{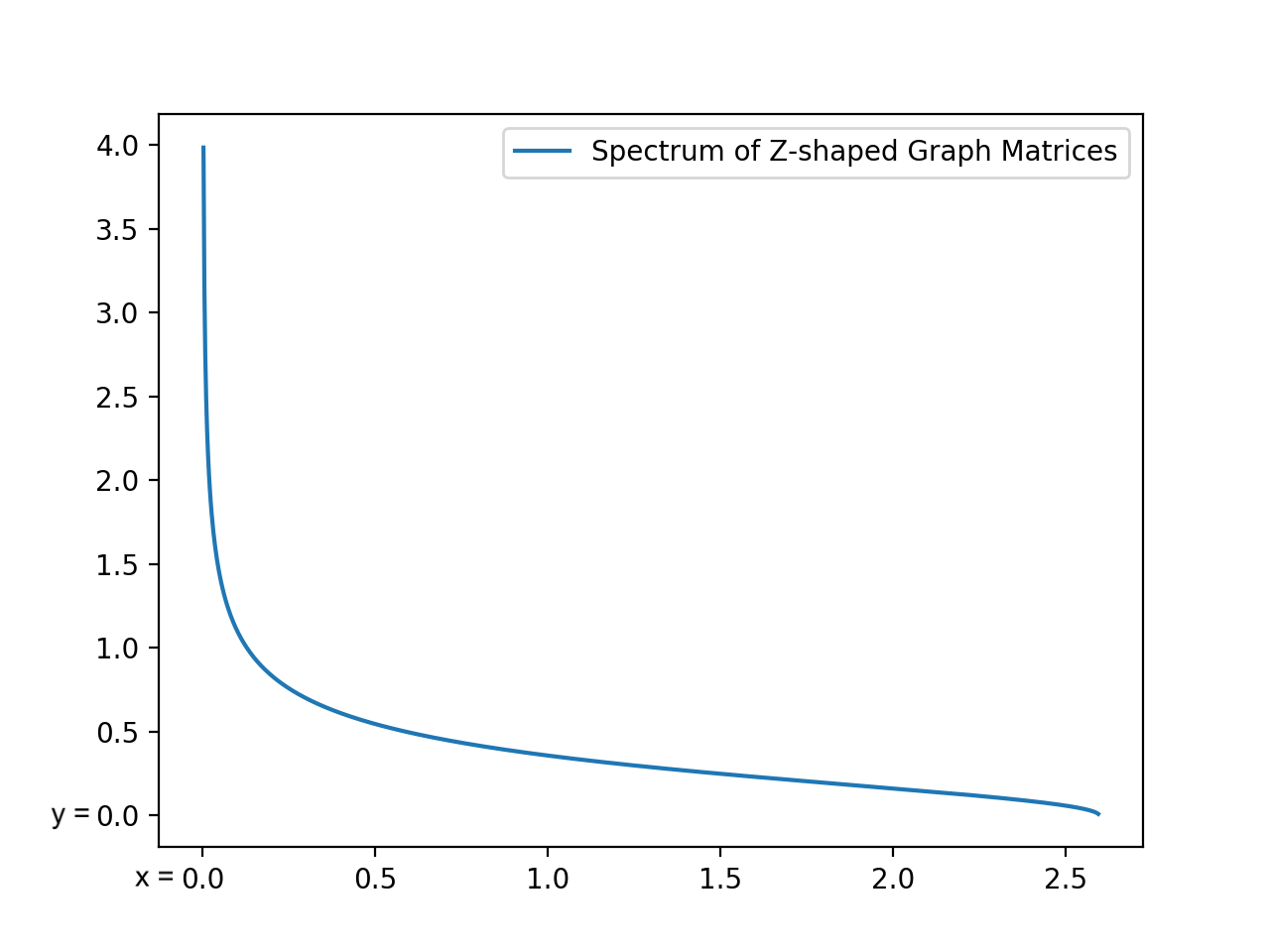}
      \caption{The spectrum of singular values}
      \label{fig:z-shape spec}
    \end{subfigure}
    \begin{subfigure}{.49\textwidth}
      \includegraphics[width=1\linewidth]{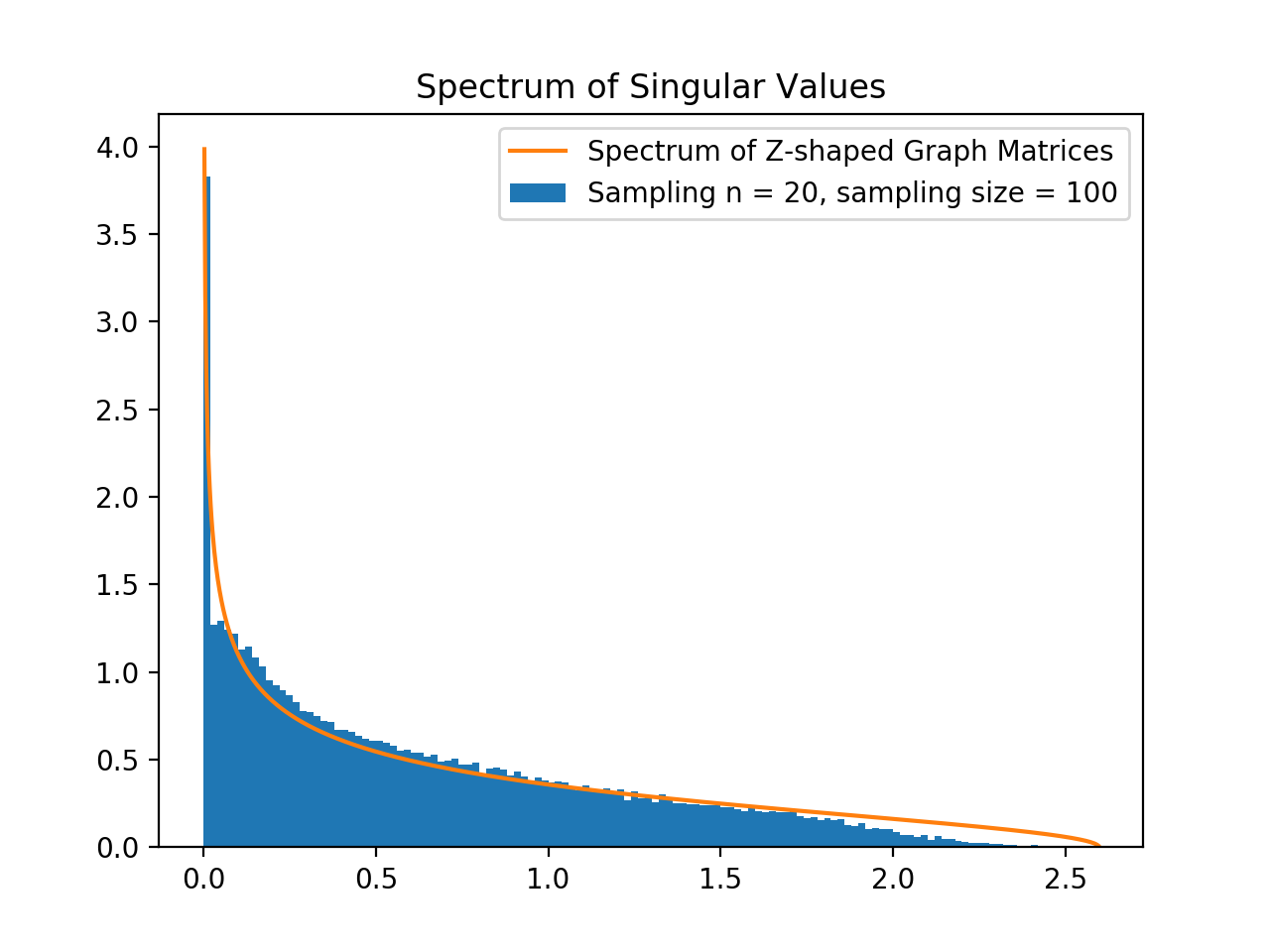}
      \caption{Sampling of singular values of $M_{\alpha_Z}$ where $n=20$}
      \label{fig:z-shape sam 20}
    \end{subfigure}
    
    \begin{subfigure}{.49\textwidth}
      \includegraphics[width=1\linewidth]{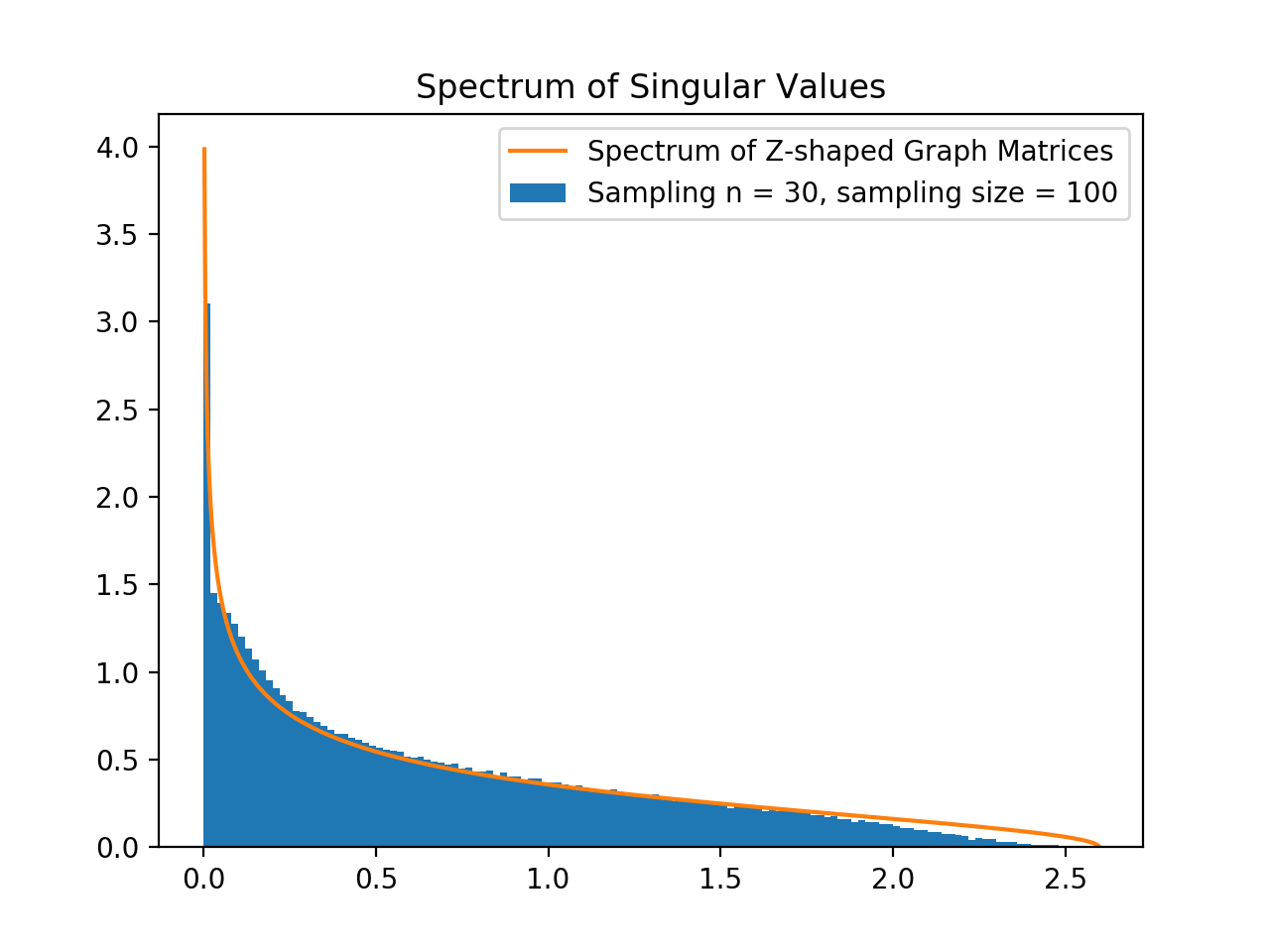}
      \caption{Sampling of singular values of $M_{\alpha_Z}$ where $n=30$}
      \label{fig:z-shape sam 30}
    \end{subfigure}
    \begin{subfigure}{.49\textwidth}
      \includegraphics[width=1\linewidth]{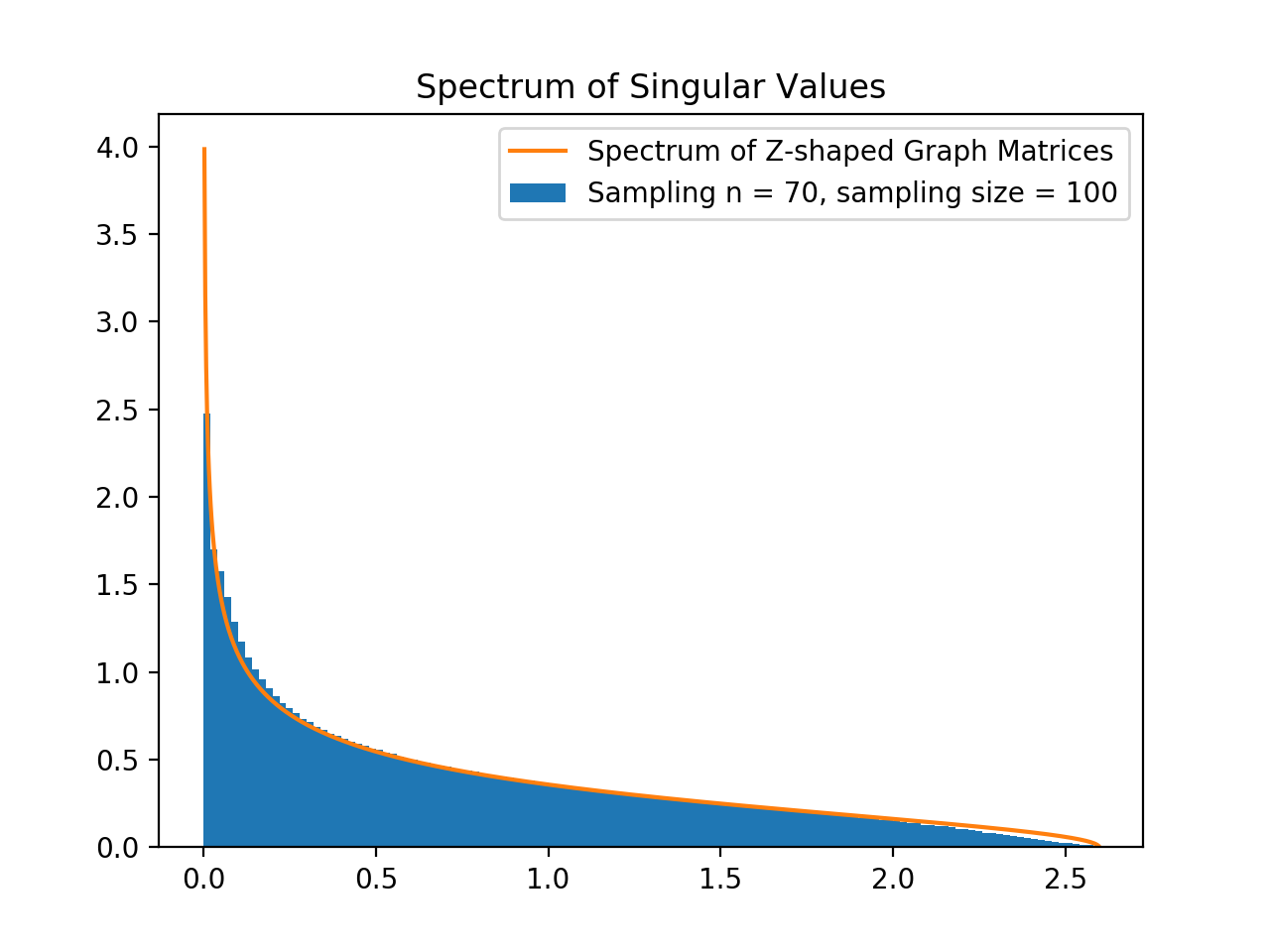}
      \caption{Sampling of singular values of $M_{\alpha_Z}$ where $n=70$}
      \label{fig:z-shape sam 70}
    \end{subfigure}
    \caption{The Spectrum of singular values of the Z-shape graph matrix and some samplings of the Z-shape graph matrices with random input graphs on $n$ vertices, for $n=20,30,70$.}
    \label{fig:z-shape spec with sampling}
\end{figure}

Figure \ref{fig:z-shape spec with sampling} shows some graphs of $g_{\alpha_Z}(x)$ and samplings of singular values of $M_{\alpha_Z}$ for $n=20,30,70$. We can see that the sampled distribution of the singular values of $M_{\alpha_Z}$ gets closer to $g_{\alpha_Z}(x)$ as $n$ gets bigger. 

\section{Counting dominant constraint graphs on \texorpdfstring{$H\(\mzshape,2k\)$}{H(aZ(m), 2k)}} \label{section:mzshape}

\setlength{\parskip}{1.5mm}
\setlength{\baselineskip}{1.3em}

In this section, we consider a generalization of the Z-shape $\a_Z$ which which we call the $m$-layer Z-shape.

\begin{defn}\label{defn:m-zshape}
\begin{figure}[hbt!]
    \centering
    \includegraphics[scale=0.35]{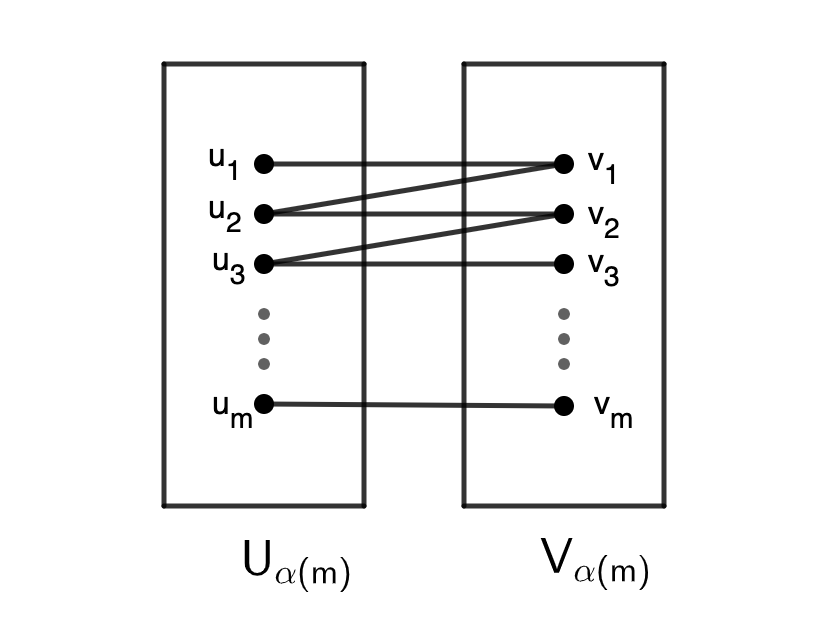}
    \caption{The $m-$layer Z-shape $\mzshape$}
    \label{fig:m-zshape}
    \end{figure}

    Let $\mzshape$ be the bipartite shape with vertices $V(\mzshape)=\{u_1,\dots,u_m,v_1,\dots, v_m\}$ and edges $E(\mzshape)=\left\{\{u_i,v_i\}: i\in[m]\right\}\cup \left\{\{u_{i+1},v_{i}\}: i\in[m-1]\right\}$ with distinguished tuples of vertices $U_{\mzshape}=(u_1,\dots,u_m)$ and $V_{\mzshape}=(v_1,\dots, v_m)$. See \Cref{fig:m-zshape} for an illustration.
    
    We refer to $\mzshape$ as the \emph{$m$-layer Z-shape}. Note that $\a_{Z(2)}$ is the Z-shape $\a_Z$ as in \Cref{defn:z-shape}.
\end{defn}

\begin{rmk}\label{rmk:mzshape-num-constr-edges}
    For the $m$-layer Z-shape $\mzshape$, $s_{\mzshape}=m$ where $s_{\mzshape}$ is the size of the minimum separator. By \Cref{lem:min constr edges for bipartite shape}, for any nonzero-valued constraint graph $C\in\C_{\(\mzshape,2k\)}$, $\abs{E(C)}\geq m\cdot(k-1)$. By \Cref{cor:dominantnumberofedges}, dominant constraint graphs $C\in\C_{\(\mzshape,2k\)}$ have $m\cdot(k-1)$ edges.
\end{rmk}

Applying $s_{\mzshape}$ to \Cref{cor:expected trace for bipartite shapes}, we get a similar result as \Cref{cor:Zshapeexpectedtrace}.

\begin{cor}\label{cor:mzshape-expected-trace}
\begin{equation*}
    \Eb\Big[\trace\(\(M_{\mzshape}M_{\mzshape}^T\)^k\)\Big] = \abs{\left\{C \in \mathcal{C}_{\(\mzshape,2k\)}: C \text{ is dominant}\right\}}\cdot n^{mk + m} \pm O\(n^{mk+m-1}\).
\end{equation*}
\end{cor}

Before stating our main result, we need one more definition, which is an extension of $\displaystyle C'_k = \dfrac{1}{2k+1}\binom{3k}{k}$ from \Cref{section:zshape}.

\begin{defn}\label{defn:D(m,n)}
For $m,n$ positive integers,
    \begin{equation}
        D(m,n)=\frac{1}{m\cdot n+1}\binom{(m+1)\cdot n}{n}\,.
    \end{equation}
\end{defn}

\begin{rmk}
     $C'_k$ from \Cref{section:zshape} is equal to $D(2,k)$ here. $D(m,n)$ is also $C_{m+1,1}(n)$ where the generalized Catalan number $C_{m,r}(k)$ is defined in \Cref{defn:generalized-calatan-number}.
\end{rmk}

An analogous theorem to \Cref{thm:zshape-trace-convergence} for $\mzshape$ is stated below.
\begin{thm}\label{thm:mzshape-trace-convergence}
For all $k \in \mathbb{N} \cup \{0\}$,
\begin{enumerate}
    \item $\Eb\Big[\trace\(\(M_{\mzshape}M_{\mzshape}^T\)^k\)\Big] = D(m,k)\cdot n^{mk + m} \pm O\(n^{mk+m-1}\)$.
    \item $\var \Big(\trace\(\(M_{\mzshape}M_{\mzshape}^T\)^k\)\Big)$ is $O\(n^{2(mk+m-1)}\)$.
\end{enumerate}
\end{thm}

For part 2, the variance trace bound, its proof is very similar to \Cref{subsection:z-shape-variance-trace-bound} and one can find the complete proof in \Cref{thm:variancewellbehaved}. We will instead just focus on the first part.

Below is the main result of this section.
\begin{thm}\label{thm:num-mzshape-constraint-graph}
    Let $\mzshape$ be the $m$-layer Z-shape as in \Cref{defn:m-zshape}. Then the number of dominant constraint graphs $C\in\C_{\(\mzshape,2k\)}$ is $D(m,k)$. 
\end{thm}

\begin{rmk}
    When $m=2$, $D(m,k)=D(2,k)=C'_k$, $\mzshape=\a_Z$, and this theorem is exactly \Cref{thm:num-zshape-constraint-graph}. 
\end{rmk}

Combined with \Cref{cor:mzshape-expected-trace}, this shows that $\Eb\Big[\trace\(\(M_{\mzshape}M_{\mzshape}^T\)^k\)\Big] = D(m,k)\cdot n^{mk + m} \pm O\(n^{mk+m-1}\)$, which is the first part of \Cref{thm:mzshape-trace-convergence}.

\subsection{Recurrence Relation for \texorpdfstring{$D(m,n)$}{D(m,n)}}

To prove the main result for this section, We need the following crucial recurrence relation for $D(m,n)$.
\begin{thm}\label{thm:generalized-calatan-recurrence-relation}
    \begin{equation}
        D(m,n+1)=\sum_{\substack{ i_0,\dots, i_{m}\geq 0: \\ i_0+\dots+i_{m}=n}} D(m,i_0)\dots D(m,i_{m})\,.
    \end{equation}
\end{thm}

\begin{proof}
    The proof is similar to the proof of \Cref{thm:catalan2-recurrence-relation}. Let \emph{$W_{m,n}:=$ the set of all grid walks from $(0,0)$ to $(n,mn)$ that are weakly below the diagonal} and $d_{m,n}=\abs{W_{m,n}}$. We will prove that $d_{m,n}=D(m,n)$ and that $d_{m,n}$ satisfies the recurrence relation in the theorem. 
    
    \begin{enumerate}[1.]
    
\setlength{\parindent}{1.5em}
\setlength{\parskip}{1.2mm}
\renewcommand{\baselinestretch}{1.2}

    \item $d_{m,n}=D(m,n)$: 
        
    For $r\in\{0,1,\dots,mn\}$, let \emph{$W_{m,n}(r):=$ the set of grid walks from $(0,0)$ to $(n,mn)$ that are $r$ steps above the diagonal}. Then $\displaystyle{\bigcup_{r=0}^{mn} W_{m,n}(r)}$ is the set of all grid walks from $(0,0)$ to $(n,mn)$, which has cardinality $\displaystyle \binom{(m+1)\cdot n}{n}$. Also $|W_{m,n}(0)|=|W_{m,n}|=d_{m,n}$. By the same proof as in the \Cref{thm:catalan2-recurrence-relation}, there is a bijection between $W_{m,n}(r)$ and $W_{m,n}(r-1)$ for each $r\in[mn]$. Thus $\displaystyle{d_{m,n}=\dfrac{1}{mn+1} \binom{(m+1)\cdot n}{n}}$. 
        
        \begin{figure}[hbt!]
            \centering
            \includegraphics[scale=0.3]{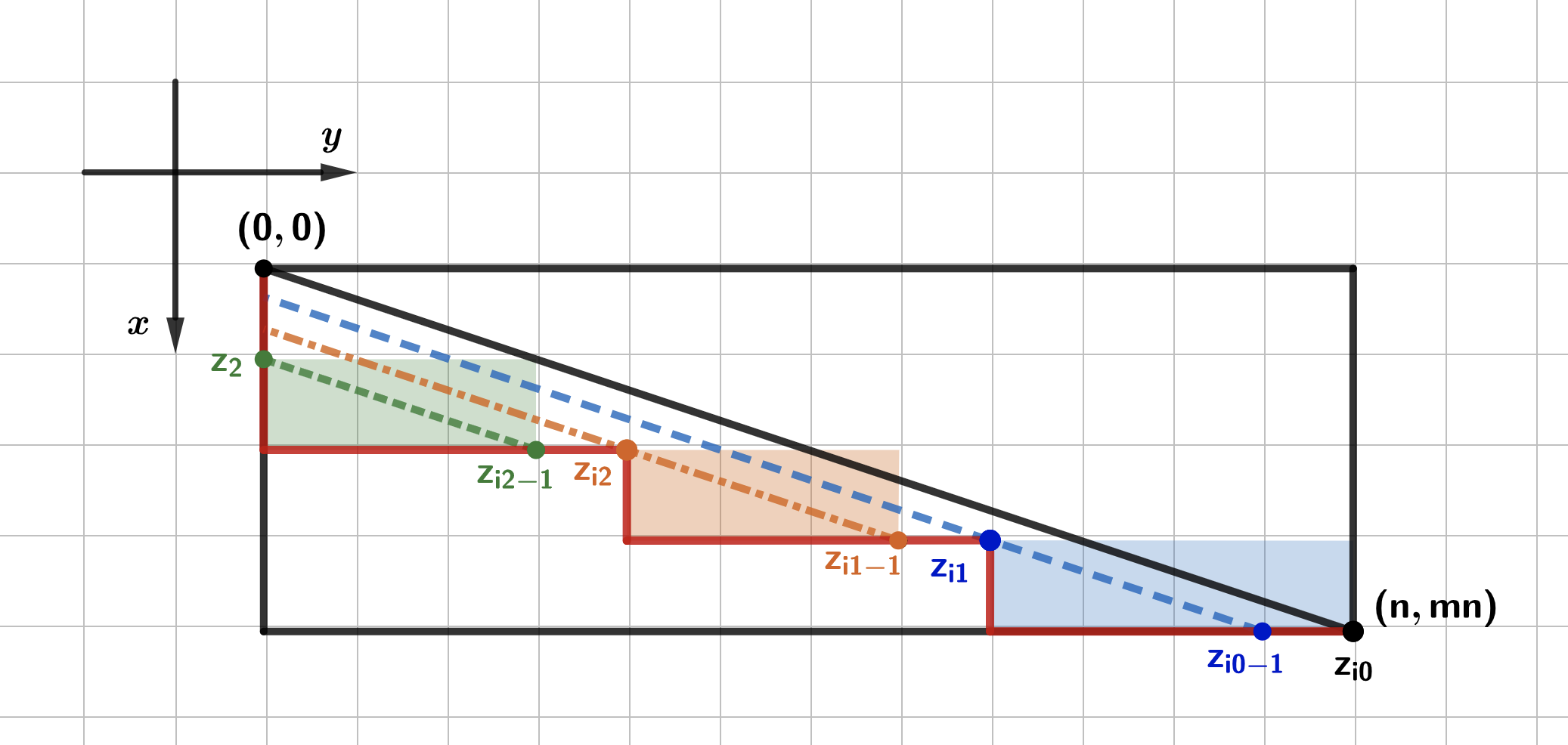}
            \caption{Illustration of the proof part 2 for \Cref{thm:generalized-calatan-recurrence-relation}.}
            \label{fig: D(m,n) recurrence}
        \end{figure}
    
    \item $\displaystyle{d_{m,n}=\sum_{\substack{ i_0,\dots, i_{m+1}\geq 0: \\ i_0+\dots+i_{m+1}=n-1}} d_{m,i_1}\dots d_{m,i_{m+1}}}$: 
        
    Similar to the proof of \Cref{thm:catalan2-recurrence-relation}, we now will find a bijection between $W_{m,n}$ and \\ $\displaystyle{\bigcup_{\substack{ i_0,\dots, i_{m}\geq 0: \\ i_0+\dots+i_{m}=n-1}} W_{m,i_0}\times\dots\times W_{m,i_{m}} }$.
    
    Let $d_k$ be the line that is shifted $k$ vertical grids down from the diagonal. i.e. $d_k$ is the line $y=m\cdot x - k$. Let $w=\(z_1,\dots,z_{n\cdot(m+1)}\)\in W_{m,n}$. 

    \begin{enumerate}[i.]
        \item Let $z_{i_0}=(a_0, m\cdot a_0)$ be the first point where $w$ touches the diagonal. Then $w_0:=\(z_{i_0},\dots,z_{n\cdot(m+1)}\)$ can be viewed as an element in $W_{m,n-a_0}$. Moreover, $w':=\(z_2,\dots,z_{i_0-1}\)$ is strictly below the diagonal $d_0$ and thus weakly below $d_1$. 
            
        \item Let $z_{i_1}=(a_1, m\cdot a_1-1)$ be the first point where $w'$ touches $d_1$. Then $w_1:=(z_{i_1},\dots,z_{i_0-1})$ can be viewed as an element in $W_{m,a_0-a_1}$. Moreover, $w':=(z_2,\dots,z_{i_1-1})$ is strictly below $d_1$ and thus weakly below $d_2$. 
        
        \item Continuing in this way, we get a sequence of points $z_{i_0},\dots,z_{i_{m-1}}$ and walks $w_0,\dots,w_{m-1}$ where each $z_{i_j}=(a_j,m\cdot a_j-j)$ is a point on $d_j$ and each $w_i$ can be viewed as an element in $W_{m,a_{i-1}-a_i}$. 
            
        \item Since $z_{i_{m-1}}$ is the first point touching $d_{m-1}$, $w_{m}=(z_2,\dots, z_{i_{m-1}-1})$ is strictly below $d_{m-1}$ and thus weakly below $d_m$. Since $d_m$ crosses $(1,0)=z_{2}$, $w_m$ can be viewed as an element in $W_{m,a_{m-1}-1}$.
    \end{enumerate}
        
    Since $(n-a_0)+(a_0-a_1)+\dots+(a_{m-1}-1)=n-1$, we conclude that $$\displaystyle{(w_0,\dots,w_m)\in\bigcup_{\substack{ i_0,\dots, i_{m}\geq 0: \\ i_0+\dots+i_{m}=n-1}} W_{m,i_0}\times\dots\times W_{m,i_{m}}}\,.$$
        
    Taking these steps in reverse gives the inverse map. it is not hard to confirm that this construction gives a bijection. 
        
    \end{enumerate}
    
    Combining statements 1 and 2, we conclude that $D(m,n)$ satisfies the recurrence relation. 
    
\end{proof}

\subsection{Properties of dominant constraint graphs on \texorpdfstring{$H\(\mzshape,2k\)$}{H(aZ(m), 2k)}}

\begin{defn}\label{defn:m-zshape-H}
    Let $\mzshape$ be the $m$-layer Z-shape as in \Cref{defn:m-zshape} and let $H\(\mzshape,2k\)$ be the multi-graph as in \Cref{def:copies}. We label the vertices of $V_{(\a_Z)_i}$ as $\{a_{ij}, b_{ij}:j\in[m]\}$ and the vertices of $V_{(\a_Z)_i}^T$ as $\{b_{ij}, a_{(i+1)j}\}$. Let $V_i=\{a_{ij},b_{ij}:i\in[q]\}$. For $j\in[m]$, we call the induced subgraph of $H\(\mzshape,2k\)$ on vertices $V_i$ the \textit{$j^{th}$ wheel $W_j$}. 
    
    We label the ``middle edges" of $H(\a,2k)$ in the following way: let $e_{2i-1,j}=\{a_{i(j+1)},b_{ij}\}$ and $e_{2i,j}=\{b_{ij},a_{(i+1)(j+1)}\}$ for $i=1,\dots,q$. For a fixed $j\in[m]$, we call the edges $e_{i,j}$'s the \textit{spokes between wheels $W_j$ and $W_{j+1}$} of $H\(\mzshape,2k\)$. See \Cref{fig:m-z-shape-H labeling} for an illustration.
    
    \begin{figure}[hbt!]
        \centering
        \includegraphics[scale=0.35]{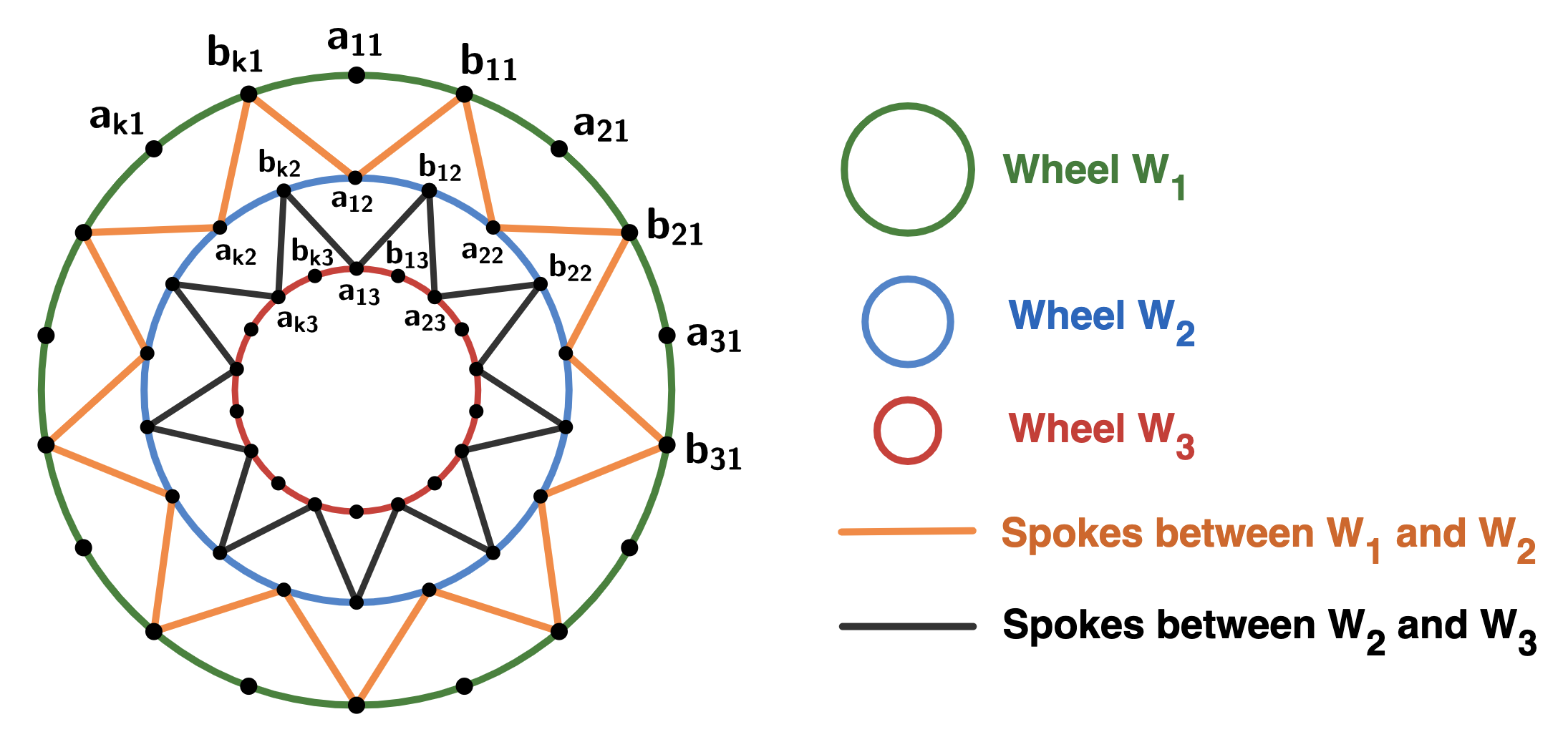}
        \caption{Illustration of \Cref{defn:m-zshape-H}: $H\(\mzshape, 2k\)$ where $m=3$.}
        \label{fig:m-z-shape-H labeling}
    \end{figure}
\end{defn}

\begin{defn}\label{defn:m-z-shape induced constraint graph}
   Let $\mzshape$ be the $m$-layer Z-shape. Let $C$ be a constraint graph on $H(\a_Z,2k)$. For $i=1,2$, We let $C_i$ denote the induced subgraph of $C$ on vertices $V_i$. We call $C_i$ the \textit{induced constraint graph} of $C$ on $V_i$. 
\end{defn}

Recall that a constraint graph $C\in\C_{(\a,2k)}$ is \emph{well-behaved} if whenever $u\=v$ in $C$, $u$ and $v$ are copies of the same vertex in $\a$ or $\a^T$. 

\begin{thm}\label{thm:mzshape dominant are well behaved}
    All dominant constraint graphs in $\C_{\(\mzshape,2k\)}$ are well-behaved.
\end{thm}
\begin{proof}
    See Appendix \ref{append:well-behaved}.
\end{proof}

\begin{prop}
If $C \in \mathcal{C}_{\(\mzshape,2k\)}$ is a dominant constraint graph, then $\displaystyle E(C) = \bigcup_{i=1}^{m} E(C_i)$ and the induced constraint graph $C_i$ is a dominant constraint graph on $W_i$ for all $i\in[m]$.
\end{prop}
\begin{proof}
Since $C$ is well-behaved, $E(C) = E(C_1) \cup \dots \cup E(C_m)$ and for all $i \in [m]$, $C_i$ is a nonzero-valued constraint graph on $W_i$. Moreover, each $C_i$ must have the minimum possible number of edges as otherwise $C$ would have too many edges.
\end{proof}

The proof of \Cref{lem:zshape-spoke-constraint} easily generalizes to $\mzshape$, yielding the following statement.
\begin{lemma}\label{lem:mzshape-spoke-constraint}
    Let $\mzshape$ be as in \Cref{defn:m-zshape} and let $C$ be a dominant constraint graph in $\C_{\(\mzshape,2k\)}$. 
    \begin{enumerate}
        \item If $a_{sj}\=a_{tj}$ for some $s < t \in [k]$ and $j\in[m-1]$, then $a_{s(j+1)}\=a_{t(j+1)}$. Moreover, the spokes $\left\{e_{x,j}:x \in [2s-1,2t-2]\right\}$ can only be made equal to each other.
        \item Similarly, if $b_{sj}\=b_{tj}$ for some $s< t \in [k]$ and $j\in \{2,3,\dots,m\}$, then $b_{s(j-1)}\=b_{t(j-1)}$. Moreover, the spokes $\left\{e_{x,j-1}:x \in [2s,2t-1]\right\}$ can only be made equal to each other.
    \end{enumerate}
\end{lemma}

\subsection{Proof of Theorem \ref{thm:num-mzshape-constraint-graph}}

We are now ready to prove the main result of this section, \Cref{thm:num-mzshape-constraint-graph}.

\begin{defn}\label{defn:H_j_mzshape}
For all $k \in \mathbb{N}$ and $j \in \{0,1,...,m\}$, we define $H_j\(\mzshape,2k\)$ to be the graph obtained by starting with $H\(\a_Z,2k\)$, merging the vertices $b_{1j'}$ and $b_{kj'}$ for all $j' \in [j]$, deleting the vertices $\left\{a_{1j'}: j' \in [j]\right\}$ and all edges incident to these vertices, and deleting the spokes incident to $a_{1j}$.

Equivalently, letting $E_j = \left\{\{b_{1j'},b_{kj'}\}: j' \in [j]\right\}$, $H_j\(\a_Z,2k\)$ is the graph obtained by taking $H\(\a_Z,2k\)/{E_j}$, deleting all multi-edges with even multiplicity, and then deleting all isolated vertices.
\end{defn}

\begin{defn}
For all $k \in \mathbb{N}$ and all $j \in \{0,1,...,m\}$, define $D_m(k,j)$ to be the number of dominant constraint graphs on $H_j\(\mzshape,2k\)$. Equivalently, $D_m(k,j)$ is the number of dominant constraint graphs $C \in \mathcal{C}_{\(\mzshape,2k\)}$ such that for all $j' \in [j]$, $b_{1j'} \=_{C} b_{kj'}$. We set $D_m(0,0) = 1$.

Note that $D_m(k,0) = \abs{C \in \mathcal{C}_{\(\mzshape,2k\)}: C \text{ is dominant}}$ and $D_m(k,m) = D_m(k-1,0)$.
\end{defn}

We can now generalize \Cref{lem:z-shape-split-1} and \Cref{lem:z-shape-split-2} to the case of a general $m$. Since the proofs are analogous to those in \Cref{section:zshape}, we will just state them.

\begin{lemma}\label{lem:mzshape-split-1}
For all $k \in \mathbb{N}$ and $j \in \{0,1,...,m-1\}$, $\displaystyle D_m(k,j) = \sum_{i=0}^{k}{D_m(i,0)D_m(k-i,j+1)}$
\end{lemma}

\begin{lemma}\label{lem:mzshape-split-2}
The constraint edges
$E = \left\{\{a_{1j'},a_{ij'}\}: j' \geq j+1\right\} \bigcup \left\{\{b_{ij'},b_{kj'}\}: j' \leq j+1\right\}$
partition $H_j\(\a_Z,2k\)$ into two parts $H_1$ and $H_2$ where $H_1 \cong H\(\a_Z,2(i-1)\)$ and $H_2 \cong H_{j+1}\(\a_Z,2(k+1-i)\)$. See \Cref{fig:mzshape-split} for an illustration.

\begin{figure}[hbt!]
    \centering
    \includegraphics[scale=0.34]{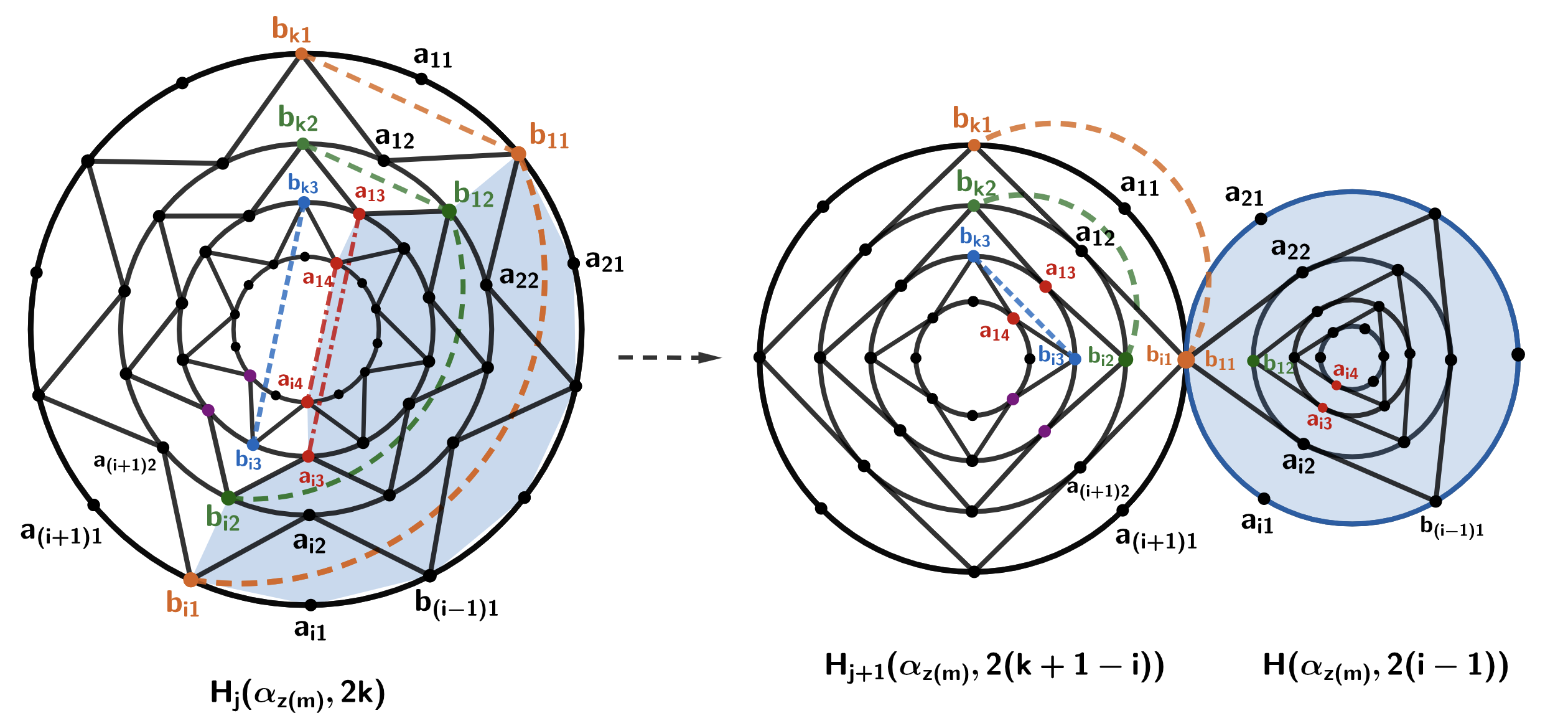}
    \caption{Illustration of \Cref{lem:mzshape-split-2}: here $j=2$.}
    \label{fig:mzshape-split}
\end{figure}
\end{lemma}

\begin{cor}\label{cor:mzshape-recurrence-relation}
For all $j\in[m]$ and $m\in\Nb$,
\begin{equation*}
    D_m(k,j) = \sum_{\substack{i_1, ..., i_{m-j+1}\geq 0: \\ i_1+\dots+i_{m-j+1} = k-1}} D_m\(i_1,0\)\cdot \dots \cdot D_m\(i_{m-j},0\)\cdot D_m\(i_{m-j+1}, 0\).
\end{equation*}

In particular, when $j=0$,
\begin{equation*}
    D_m(k,0) = \sum_{\substack{i_1, ..., i_{m+1}\geq 0: \\ i_1+\dots+i_{m+1} = k-1}} D_m\(i_1,0\)\cdot \dots \cdot D_m\(i_{m},0\)\cdot D_m\(i_{m+1}, 0\).
\end{equation*}
\end{cor}

\begin{proof}
By \Cref{lem:mzshape-split-1},
    \begin{align*}
        D_m(k,j) 
        &= \sum_{i_1=0}^{k}D_m(i_1,0)D_m(k-i_1,j+1) \\
        &= \sum_{i_1=0}^{k}D_m(i_1,0)\(\sum_{i_2=0}^{k-i_1} D_m(i_2,0)D_m(k-i_1-i_2, j+2)\) \\
        &= \ldots \\
        &= \sum_{\substack{i_1, ..., i_{m-j+1}\geq 0: \\ i_1+\dots+i_{m-j+1} = k}} D_m(i_1,0)D_m(i_2,0)\dots D_m\(i_{m-j},0\)D_m\(i_{m-j+1}, m\) \\
        &= \sum_{\substack{i_1, ..., i_{m-j+1}\geq 0: \\ i_1+\dots+i_{m-j+1} = k}} D_m(i_1,0)D_m(i_2,0)\dots D_m\(i_{m-j},0\)D_m\(i_{m-j+1}-1, 0\) \\
        &= \sum_{\substack{i_1, ..., i_{m-j+1}\geq 0: \\ i_1+\dots+i_{m-j+1} = k-1}} D_m(i_1,0)D_m(i_2,0)\dots D_m\(i_{m-j},0\)D_m\(i_{m-j+1}, 0\) \\
    \end{align*}
\end{proof}

By \Cref{cor:mzshape-recurrence-relation} and \Cref{thm:generalized-calatan-recurrence-relation}, $D_m(k,0)$ has the same recurrence relation as the generalized Calatan number $\displaystyle D(m,k) = \dfrac{1}{mk+1}\binom{(m+1)k}{k}$, thus \Cref{thm:num-mzshape-constraint-graph} is proved.

\section{The spectrum of an $m$-layer Z-shaped graph matrix} \label{section:mzshape-spectrum}

\setlength{\parskip}{1.5mm}
\setlength{\baselineskip}{1.3em}

In this section we aim to find the spectrum of the singular values for m-layer Z-shape graph matrices. 
By \Cref{thm:mzshape-trace-convergence}, we get the following corollary using similar analysis as in the proof for \Cref{lem:evenmomentsmatch}.

\begin{cor}\label{lem:mzshape-even-moments}
    Let $\{G_n:n \in \mathbb{N}\}$ be a sequence of $G\(n,\frac{1}{2}\)$ graphs and $M_{m,n} = \dfrac{1}{n^{m/2}}{M_{\mzshape}(G_n)}$. With probability $1$, for all $k \in \mathbb{N} \cup \{0\}$, $\displaystyle \lim_{n\to\infty}{\Eb_{X \sim D_{M_{m,n}}} [X^{2k}]} = D(m,k)$.
\end{cor}

Thus if we can find a function $g_m$ such that $\displaystyle{\int_{0}^{\infty} g_m(x)x^{2k} \dx = D(m,k)}$, then $g_m$ describes the limiting spectrum of singular values for the m-layer Z-shape graph matrix as $n$ goes to $\infty$. 

Recall that $\displaystyle{D(m,n)= \dfrac{1}{mn+1} \binom{(m+1)n}{n}}$ as defined in Definition \ref{defn:D(m,n)}. In this section we will generalize the arguments for the $m=2$ case in \Cref{section:zshape-spectrum} to the case $m=3$. The general steps will be:
\begin{enumerate}
    \item Assume $\displaystyle\int_{0}^{\infty} f(x)x^{2k}\dx = D(3,k)$ and derive a differential equation for $f(x)$.
    \item Prove that under this differential equation, the moments of $f(x)$ are indeed $D(3,k)$.
\end{enumerate}   

\begin{thm}\label{thm:ODE for f(x) of 3-Z-shape}
    Let $a=\lim_{n\to\infty} D(3,n+1)/D(3,n)=\dfrac{16}{3\sqrt{3}}$. If $f(x)$ is a function such that $\displaystyle\int_{0}^{a} f(x)x^{2k} dx=D(3,k)$ for all nonnegative integers $k$ and 
    \begin{enumerate}
        \item $f$ is three times continuously differentiable on $(0,a)$ and $\lim_{x\to a^-} f(x)=0$,
        \item $\lim_{x\to a^-}(x-a)f'(x)=0$,
        \item $\lim_{x\to a^-}2(x-a)f''(x)+f'(x)=0$,
        \item $\lim_{x\to 0^+}xf(x)=0$, $\lim_{x\to 0^+}x^2f'(x) = 0$, and $\lim_{x\to 0^+}x^3f''(x)=0$.
    \end{enumerate}
    then $f$ satisfies the following ODE on $(0,a)$:
    \begin{equation}
    (27x^4-256x^2)f'''(x) + (162x^3-768x)f''(x)+ (177x^2-192x)f'(x)+15xf(x)=0.
    \label{eqn:ODE for 3-Z-shape f(x)}
    \end{equation}
\end{thm}

\begin{lemma}\label{lem:zero odd moments}
    Let $a$ be some positive constant. If $f$ is continuous on $[0,a]$ and \\ $\displaystyle\int_{0}^{a} f(x)x^{2n+1}\dx=0$ for all nonnegative integers $n$, then $f=0$ on $(0,a)$.
\end{lemma}
\begin{proof}
    The proof is similar to the proof for \Cref{lem:zero even moments} for the case of all zero even moments. Here we approximate $f(\sqrt{x})\cdot\sqrt{x}$ by a polynomial $p(x)$, so that the odd polynomial $p(x^2)/x$ ($p(x)$ has $0$ constant term) approximates $f(x)$. 
\end{proof}

\begin{proof}[Proof of \Cref{thm:ODE for f(x) of 3-Z-shape}]
    Denote the LHS of the ODE by $G(x)$. Let $\displaystyle A(k,n)=\int_{0}^a f^{(k)}(x)\cdot x^n\dx$ and $\displaystyle B(k,n)=\[f^{k}(x)\cdot x^n\]_0^a$. Repeatedly doing integration by parts we get that
    \begin{align*}
        A(m,n) 
        &= B(m-1,n) - n\cdot A(m-1,n-1)\\
        &= B(m-1,n)-nB(m-2,n-1) + n(n-1)A(m-2,n-2) \\
        &= B(m-1,n)-nB(m-2,n-1) + n(n-1)B(m-3,n-2) \\ 
        &\hspace{0.5cm} -n(n-1)(n-2)A(m-3,n-3).
        \label{eqn:3-layer-z-A(m,n)}
    \end{align*}
    
    Also, we have that
    \begin{equation}
        \dfrac{D(3,n)}{D(3,n-1)}=\dfrac{A(0,2n)}{A(0,2n-2)}=\dfrac{4(4n-3)(4n-2)(4n-1)}{(3n+1)(3n)(3n-1)}.
        \label{eqn:D(3,n) ratio}
    \end{equation}
    
    The steps for deducing the ODE for $f(x)$ are very similar to the steps used in the proof of \Cref{thm:ODE for f(x) of Z-shape} to deduce the ODE for the Z-shaped graph matrix. We first apply $m=3$ and $n=2n+3$ to the first equation and rewrite the term $n(n-1)(n-2)A(m-3,n-3)$ using the second equation \eqref{eqn:D(3,n) ratio}. We then gradually eliminate the non-constant coefficients in front of $A(m,n)$'s using the first equation.
    
    Plugging in $m=3,n=2n+3$ into the first equation, we get 
    \begin{align*}
        27A(3,2n+3) 
        &= 27B(2,2n+3)-27(2n+3)B(1,2n+2)+\\
        &27(2n+3)(2n+2)B(0,2n+1) - 27(2n+3)(2n+2)(2n+1)A(0,2n).
    \end{align*}
    
    Rewriting the last term on the RHS above and applying the second equation \eqref{eqn:D(3,n) ratio}, we get
    \begin{align*}
        27(2n+3)&(2n+2)(2n+1)A(0,2n) \\
        & = 8(3n+1)(3n)(3n-1)A(0,2n) -\\
        & \hspace{1cm} 81\cdot 2(2n+2)(2n+1)A(0,2n) + 59\cdot 3(2n-1)A(0,2n) -15A(0,2n)\\
        & = 32(4n-3)(4n-2)(4n-1)A(0,2n-2) - \\
        & \hspace{1cm} 81\cdot 2(2n+2)(2n+1)A(0,2n) + 59\cdot 3(2n-1)A(0,2n) -15A(0,2n)
    \end{align*}
    
    We can rewrite the first term on the RHS as
    \begin{align*}
        32(4n-3)&(4n-2)(4n-1)A(0,2n-2)\\
        & = 256(2n+1)(2n)(2n-1)A(0,2n-2) - \\
        & \hspace{1cm} 256\cdot 3(2n)(2n-1)A(0,2n-2) + 128\cdot 3(2n-1)A(0,2n-2).
    \end{align*}
    
    Now apply the first equation to all the $(n+1)A(m,n)$, $(n+2)(n+1)A(m,n)$ and $(n+3)(n+2)(n+1)A(m,n)$ above, group together the $A(m,n)$ terms and $B(m,n)$ terms separately, and rewrite the $B(m,n)$ term using the definition of $B(m,n)$. We get that for all $n\geq 1$,
    \begin{align*}
        &27A(3,2n+3)-256A(3,2n+1) + 2\cdot 81A(2,2n+2)-3\cdot 256A(2,2n)  \\
        &\hspace{1cm} + 3\cdot 59 A(1,2n+1)-3\cdot 128A(1,2n-1) + 15A(0,2n)\\
        &= \[\((27x^2-256)f''(x)+27xf'(x)\)x^{2n+1}\]_0^a \\
        &\hspace{1cm}+ \[\((2n-2)(-27x^2+256)\)f'(x)x^{2n}\]_0^{a} + \[p(n,x)\cdot f(x)x^{2n-1}\]_0^a
    \end{align*}   
    where $p(n,x)$ is some polynomial in terms of $n$ and $x$. 
    
    Observe that the last term on the RHS is $0$ since $\displaystyle \lim_{x\to 0^+}xf(x)=0$ and $\displaystyle \lim_{x\to a^-}f(x)=0$ by assumption. The second last term is $0$ since $27x^2-256=27(x+a)(x-a)$, $\displaystyle \lim_{x\to a^-} (x-a)f'(x)=0$ and $\displaystyle \lim_{x\to 0^+}x^2f'(x)=0$. The first term top part is $\displaystyle \lim_{x\to a^-}27\((x+a)(x-a)f''(x)+xf'(x)\)x^{2n+1} \\= \lim_{x\to a^-} 27\(2a(x-a)f''(x)+af'(x)\)a^{2n+1} = 0 $ since $\displaystyle \lim_{x\to a^-} 2(x-a)f''(x)+f'(x)=0$ by assumption. The bottom part is $0$ since $\lim_{x\to 0^+}x^3f''(x)=0$ by assumption. Thus the RHS is $0$.
       
    Expanding out each $A(m,n)$ term by using the definition of $A(m,n)$, we get that \\ $\displaystyle\int_{0}^a G(x)x^{2n+1}=0$ for all $n\geq 1$. By \Cref{lem:zero odd moments}, $G(x)=0$ on $(0,a)$, which proves that $f$ satisfies the ODE.
    
\end{proof}

\begin{thm}\label{thm:verify ODE for 3-Z-shape spec}
    Let $a=\lim_{n\to\infty} D(3,n)/D(3,n-1) = \dfrac{16}{3\sqrt{3}}$ and let $f$ be a function satisfying the following ODE
    \begin{equation}
        (27x^4-256x^2)f'''(x) + (162x^3-768x)f''(x)+ (177x^2-192)f'(x)+15xf(x)=0
        \label{eqn: ODE for spectrum of 3-Z-shape}
    \end{equation}
    and the conditions listed in \Cref{thm:ODE for f(x) of 3-Z-shape}. Moreover, assume $\displaystyle\int_{0}^a f(x)\dx = 1$. Then for any nonnegative integer $k$,
    \begin{equation}
        A(0,2k)=\int_{0}^{a} x^{2k}\cdot f(x)\dx = D(3,k).
        \label{eqn:A(0,2k)=D(3,k)}
    \end{equation}
\end{thm}

The proof is very similar to the proof for \Cref{thm:verify ODE for Z-shape spec}. We integrate the ODE from $0$ to $a$, do integration by parts and use the conditions for $f$ to eliminate the redundant terms and finally arrive at the ratio between $A(0,2(k-1))$ and $A(0,2k)$ which matches the ratio between $D(3,k-1)$ and $D(3,k)$. By induction on $k$ we conclude \eqref{eqn:A(0,2k)=D(3,k)}.

\begin{cor}
Let $M_n = \dfrac{1}{n^{3/2}}\cdot M_{\a_{Z(3)}}$ where $M_{\a_{Z(3)}}$ is the graph matrix with random input graph $G\sim G(n,1/2)$ and $\a_{Z(3)}$ is the multi-Z-shape defined in \Cref{defn:m-zshape}. Let $g(x)$ be $f(x)$ as in \Cref{thm:verify ODE for 3-Z-shape spec} on $(0,a)$ and $0$ for $x\geq a$. With probability $1$, the spectrum of the singular values of $M_n$ approaches $g(x)$ weakly as $n\to\infty$.
\end{cor}

For this ODE \eqref{eqn:ODE for 3-Z-shape f(x)}, WolframAlpha fails to give us an explicit solution. Instead, we solve the ODE numerically by approximating the tail segment of $f(x)$ by $c\cdot(a-x)^{r}$ for some constants $c$ and $r$.

\begin{description}
\item[Step 1:] We analyze the behaviour of the ODE when $x$ is very close to $a$. Notice that $a=\dfrac{16}{3\sqrt{3}} \iff 27a^2-256=0$.
\begin{enumerate}
    \item $27x^4-256x^2 = 27x^2(x-a)(x+a) \sim 52a^3(x-a)$.
    \item $162x^3-768x = 81x^3 + 3x(27x^2-256) \sim 81a^3$.
    \item $177x^2-192 = \dfrac{3}{4}(209x^2+27x^2-256) \sim \dfrac{3\cdot 209a^2}{4}$.
\end{enumerate}

Thus when $x$ is very close to $a$, the ODE is
\begin{equation*}
    52a^3(x-a)f'''(x) + 81a^3f''(x) + \(\dfrac{3\cdot 209a^2}{4}\)f'(x)=0.
\end{equation*}

One can check that $f'(x)=C'\((a-x)^{-1/2}+\dfrac{209}{64\sqrt{3}}(a-x)^{1/2}\)$ is a solution to the above ODE. Thus
\begin{equation}
    f(x)\sim g(x)= C\((a-x)^{1/2} + \dfrac{209}{64\cdot3\sqrt{3}}(a-x)^{3/2}\) 
    \label{eqn:approx tail seg of f(x) for 3-Z-shape}
\end{equation}
for some constant $C$ when $x$ is very close to $a$.

\item[Step 2:] We approximate the solution of $f(x)$ by approximating the tail segment of $f(x)$ (where $|x-a|<\epsilon$) by $g(x)$ and use this approximation to obtain initial conditions for the ODE for $f$. In particular, we choose a small $\epsilon>0$ and set the initial conditions for the ODE as follows: 
$$f(a-\epsilon)= g(a-\epsilon), f'(a-\epsilon)=g'(a-\epsilon), f''(a-\epsilon)=g''(a-\epsilon).$$ 

We calculate the constant $C$ in $g$ by noticing that the integration of $f$ over $(-a,a)$ should be $1$. 

Setting $\epsilon=0.01$ and solving the ODE in python, we get the following plot for the solution to the ODE (concatenated with a tail segment where we use $g$ as an approximation):
\begin{figure}[hbt!]
    \centering
    \begin{subfigure}[t]{.45\textwidth}
      \includegraphics[width=1\linewidth]{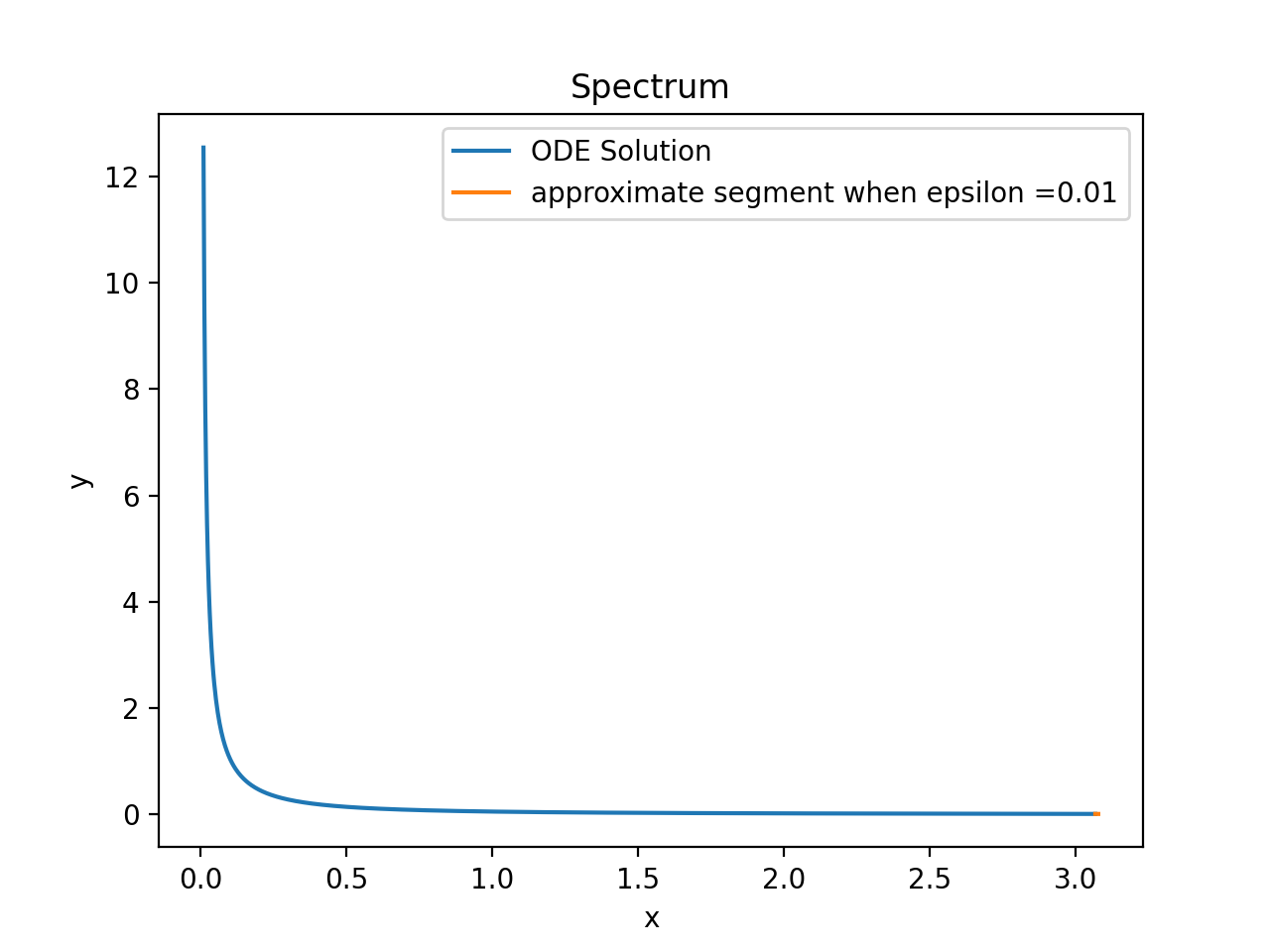}
      \caption{Plot of the ODE solution with the approximated tail segment.}
      \label{fig:eps=0.01}
    \end{subfigure}
    \begin{subfigure}[t]{.45\textwidth}
      \includegraphics[width=1\linewidth]{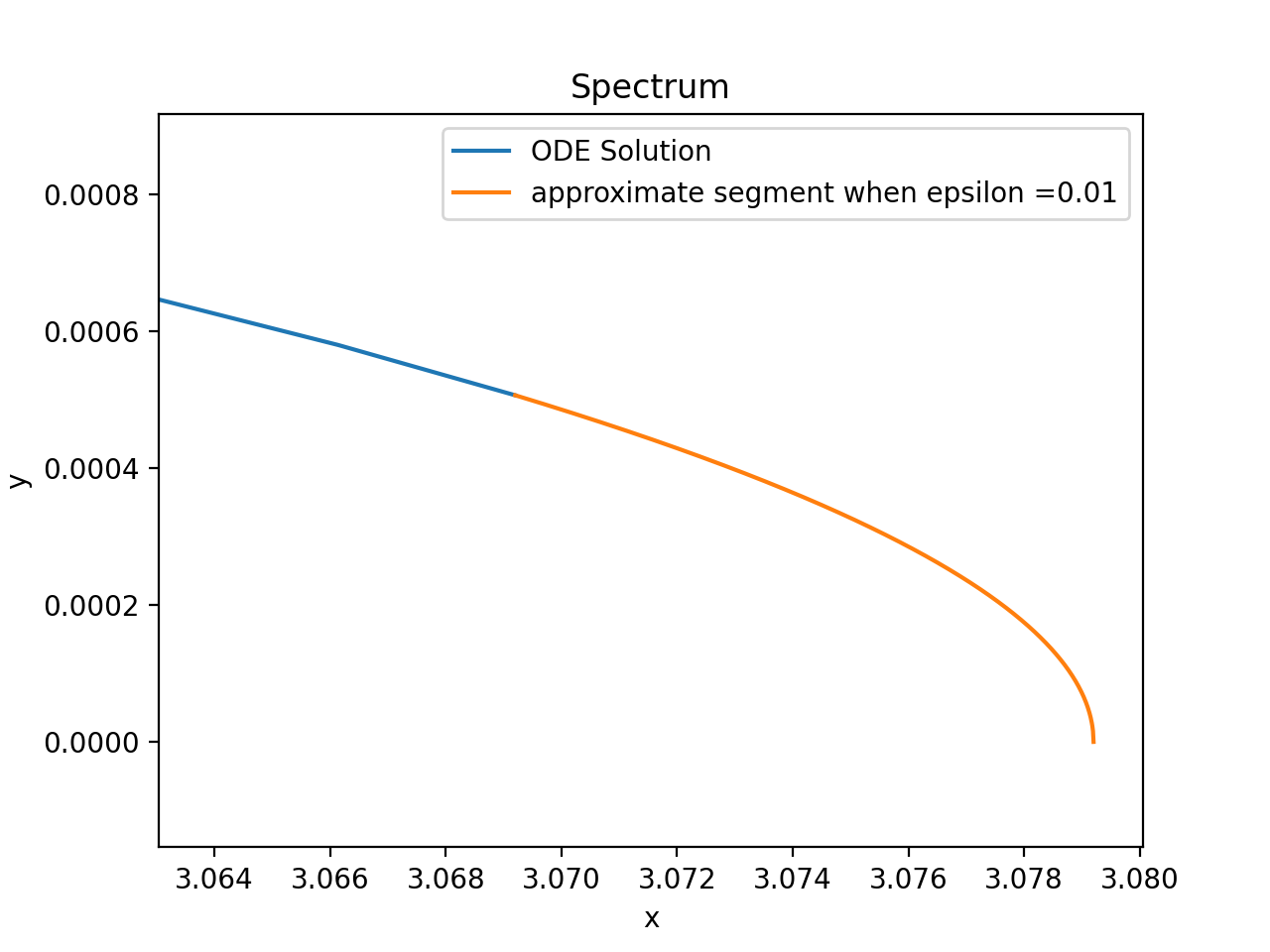}
      \caption{Zoom in at the tail segment.}
      \label{fig:eps=0.01 zoom}
    \end{subfigure}
    \caption{The plot of the spectrum where $x>0$.}
    \label{fig:eps-0.01}
\end{figure}

\end{description}

\begin{figure}[hbt!]
    \centering
    \begin{subfigure}[t]{.45\textwidth}
      \includegraphics[width=1\linewidth]{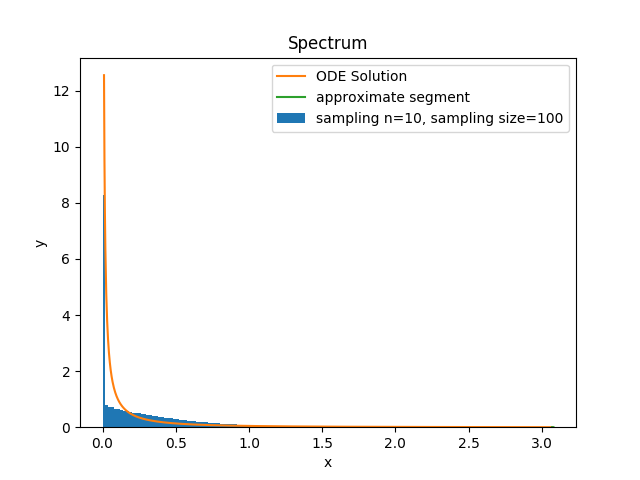}
      \caption{Sampling of the singular values of $M_n$ for $n=10$.}
      \label{fig:sampling zz n=10}
    \end{subfigure}
    \begin{subfigure}[t]{.45\textwidth}
      \includegraphics[width=1\linewidth]{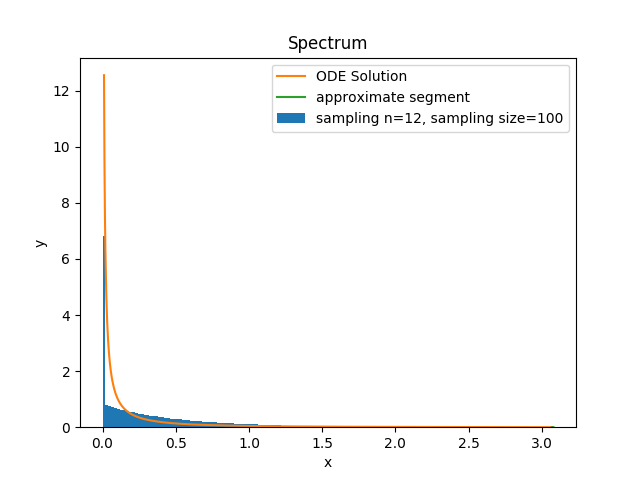}
      \caption{Sampling of the singular values of $M_n$ for $n=12$.}
      \label{fig:sampling zz n=12}
    \end{subfigure}
    \caption{The approximated spectrum of the singular values of the 3-layer Z-shaped graph matrix and some samplings of the singular values of $M_n$, for $n=10$ and $n=12$.}
    \label{fig:sampling zz}
\end{figure}

To test this solution experimentally, we can sample from the distribution of singular values of $M_n$ by sampling a random graph $G$, computing the resulting matrix $M_n(G)$, and computing its singular values. See \Cref{fig:sampling zz} for a plot of the approximated spectrum together with the empirical distribution of the singular values of $M_n$ with $n=10$ and $n=12$, respectively (where we sampled $100$ random graphs $G$).

\bibliography{s0_main}

\begin{appendix}
\section{Dominant Constraint Graphs are Well-Behaved}\label{append:well-behaved}

\setlength{\parskip}{1.5mm}
\setlength{\baselineskip}{1.3em}

In this section, we prove that dominant constraint graphs on $H(\a_{Z_m},2k)$ and $H(\a_{Z_m},2k) \cup H(\a_{Z_m},2k)$ are well-behaved. More precisely, we show the following theorems.
\begin{thm}\label{thm:dominantwellbehaved}
    If $C \in \mathcal{C}_{H(\a_{Z_m},2k)}$ and $\val(C) \neq 0$ then $C$ has at least $m(k - 1)$ edges. Moreover, if $\val(C) \neq 0$ and $C$ is not well-behaved then $C$ has at least $m(k-1) + 1$ edges.
\end{thm}
\begin{thm}\label{thm:variancewellbehaved}
If $C \in \mathcal{C}_{H(\a_{Z_m},2k) \cup H(\a_{Z_m},2k)}$ and $\val(C) \neq 0$ then $C$ has at least $2m(k - 1)$ constraint edges. Moreover, if $\val(C) \neq 0$ but the restriction of $C$ to either copy of $H(\a_{Z_m},2k)$ has value $0$ (i.e., $\val(C) \neq 0$ because of constraint edges between the two copies of $H(\a_{Z_m},2k)$) then $C$ has at least $2m(k - 1) + 2$ edges.
\end{thm}
\subsection{The set of graphs \texorpdfstring{$R(H(\alpha_{Z_m},2k))$}{R(H(Z(m),2k))}}
The high level idea for proving Theorems \ref{thm:dominantwellbehaved} and \ref{thm:variancewellbehaved} is as follows. We show that for any constraint graph which has nonzero value and does not have too many constraint edges, we can find a vertex $v \in H(\alpha_{Z_m},2q)$ which is not incident to any constraint edges or spokes of odd multiplicity. This implies that there must be a constraint edge between the two neighbors of $v$ in its wheel. We can then merge the two neighbors of $v$ together, delete $v$ and all edges incident to $v$, and repeat this process. When we are done, we will be left with a cycle of length $2$ for each wheel. For $H(\alpha_{Z_m},2q)$, this implies that the constraint graph is either well-behaved or has an extra edge. For $H(\alpha_{Z_m},2q) \cup H(\alpha_{Z_m},2q)$, this implies that the restriction of the constraint graph to each copy of $H(\alpha_{Z_m},2q)$ has nonzero value.

In order to implement this strategy, we need to analyze the more general class of graphs $H$ which can be obtained by starting with $H(\alpha_{Z_m},2q)$, iteratively taking vertices which are not incident to any spokes with odd multiplicity, merging their neighbors together, and deleting the vertex and all edges incident to it.

\begin{defn}\label{defn:R-H}
\begin{figure}[hbt!]
    \centering
    \begin{subfigure}[c]{.63\textwidth}
        \includegraphics[width=1\linewidth]{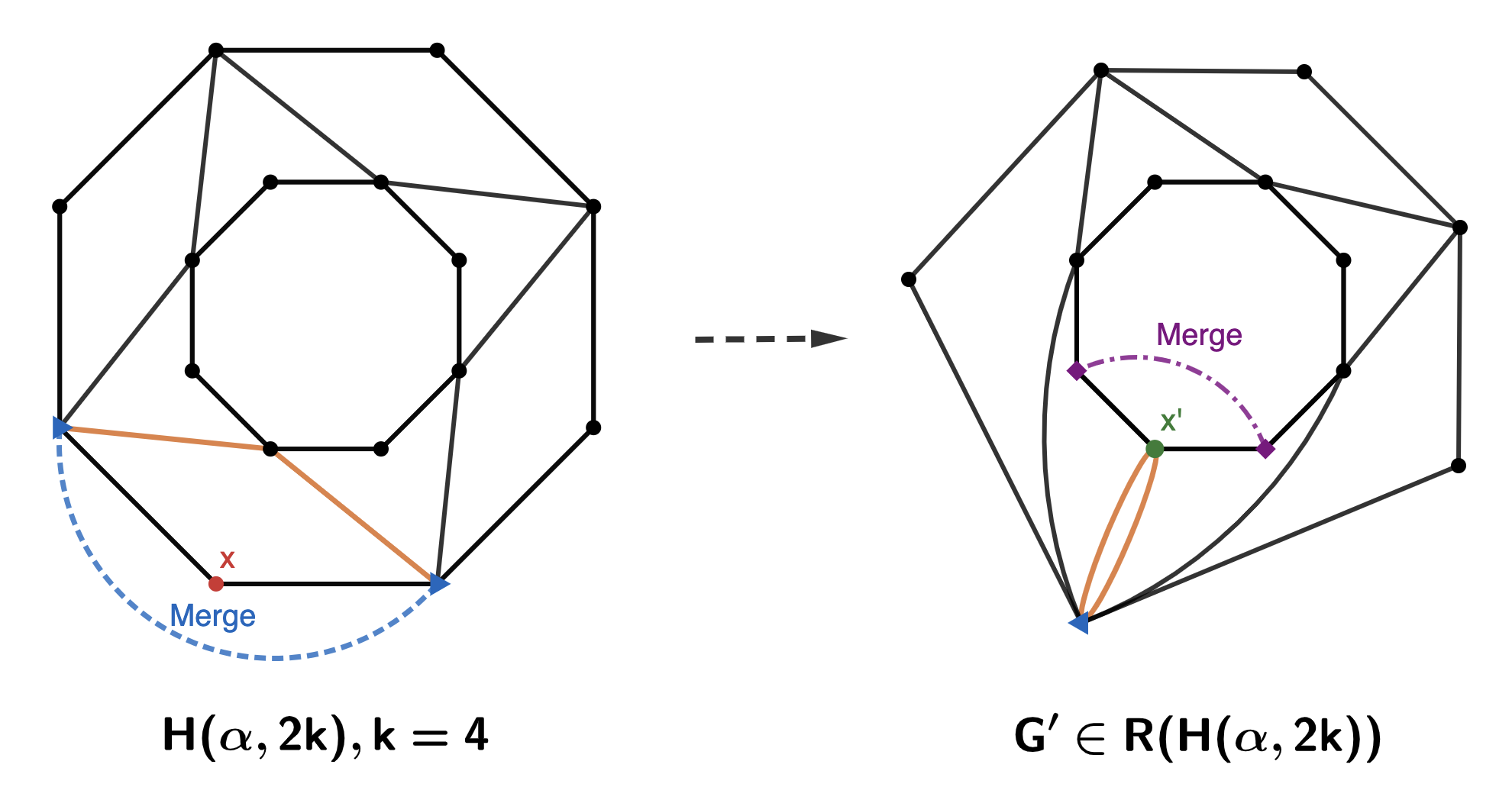}
    \end{subfigure}
    \begin{subfigure}[c]{.34\textwidth}
        \includegraphics[width=1\linewidth]{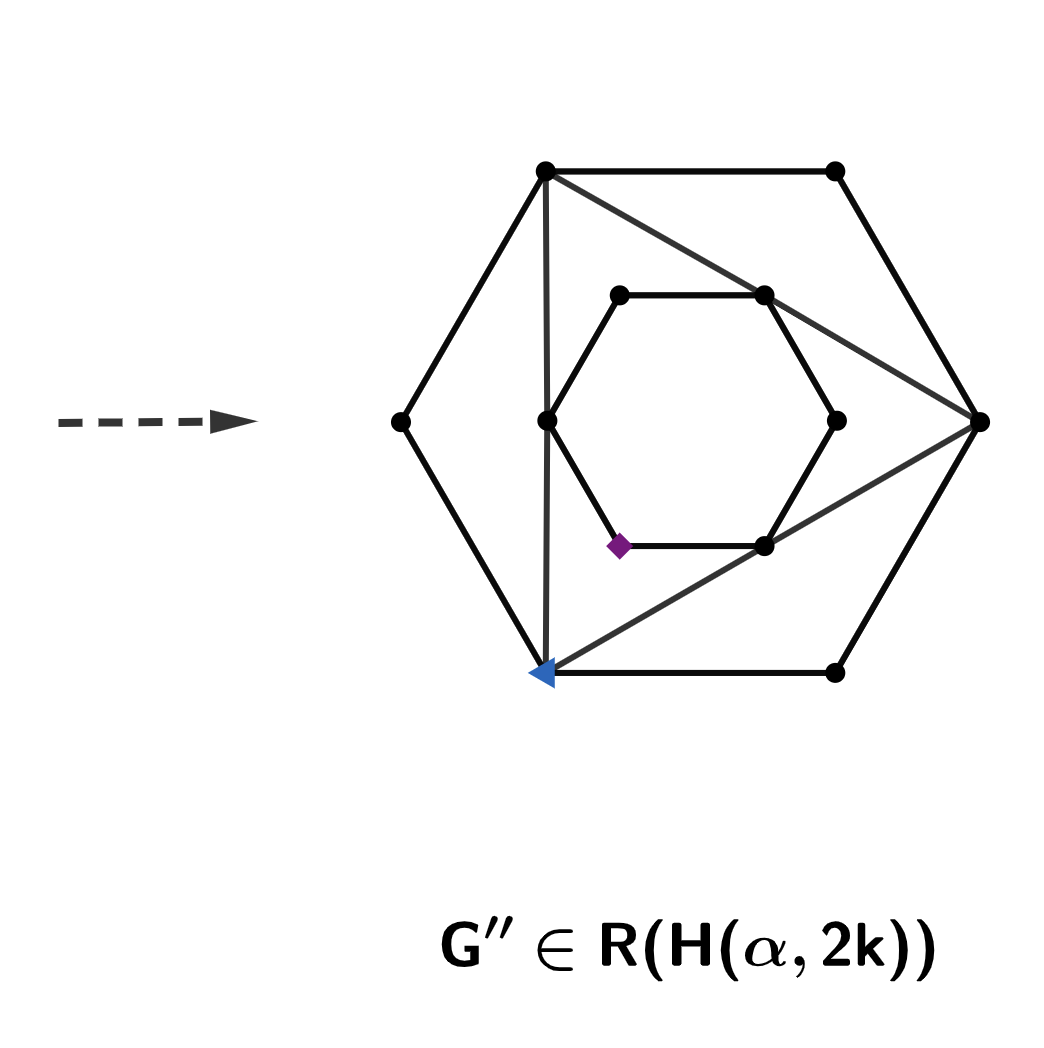}
    \end{subfigure}
    \caption{Illustration of \Cref{defn:R-H}.}
    \label{fig:R-H}
\end{figure}

Define $R(H(\alpha_{Z_m},2q))$ to be the set of graphs $H$ which can be obtained by starting from $H(\alpha_{Z_m},2q)$ and repeatedly applying the following contraction operation:
\begin{enumerate}
\item Choose a vertex $x \in V(H)$  such that $x$ is in a wheel with at least $4$ vertices and all spokes incident to $x$ have even multiplicity. Merge the two neighbors of $x$ in its wheel and delete all edges incident to $x$.
\end{enumerate}

See \Cref{fig:R-H} for an illustration.
\end{defn}
To analyze what happens when we apply contraction operations to $H(\alpha_{Z_m},2k)$, it is useful to consider the $k$ copies of $\a_{Z_m}$ and the $k$ copies of $\a_{Z_m}^T$ which $H(\alpha_{Z_m},2k)$ is built from.

\begin{defn}\label{defn:H-alpha-i}
\begin{figure}[hbt!]
    \centering
    \includegraphics[scale = 0.35]{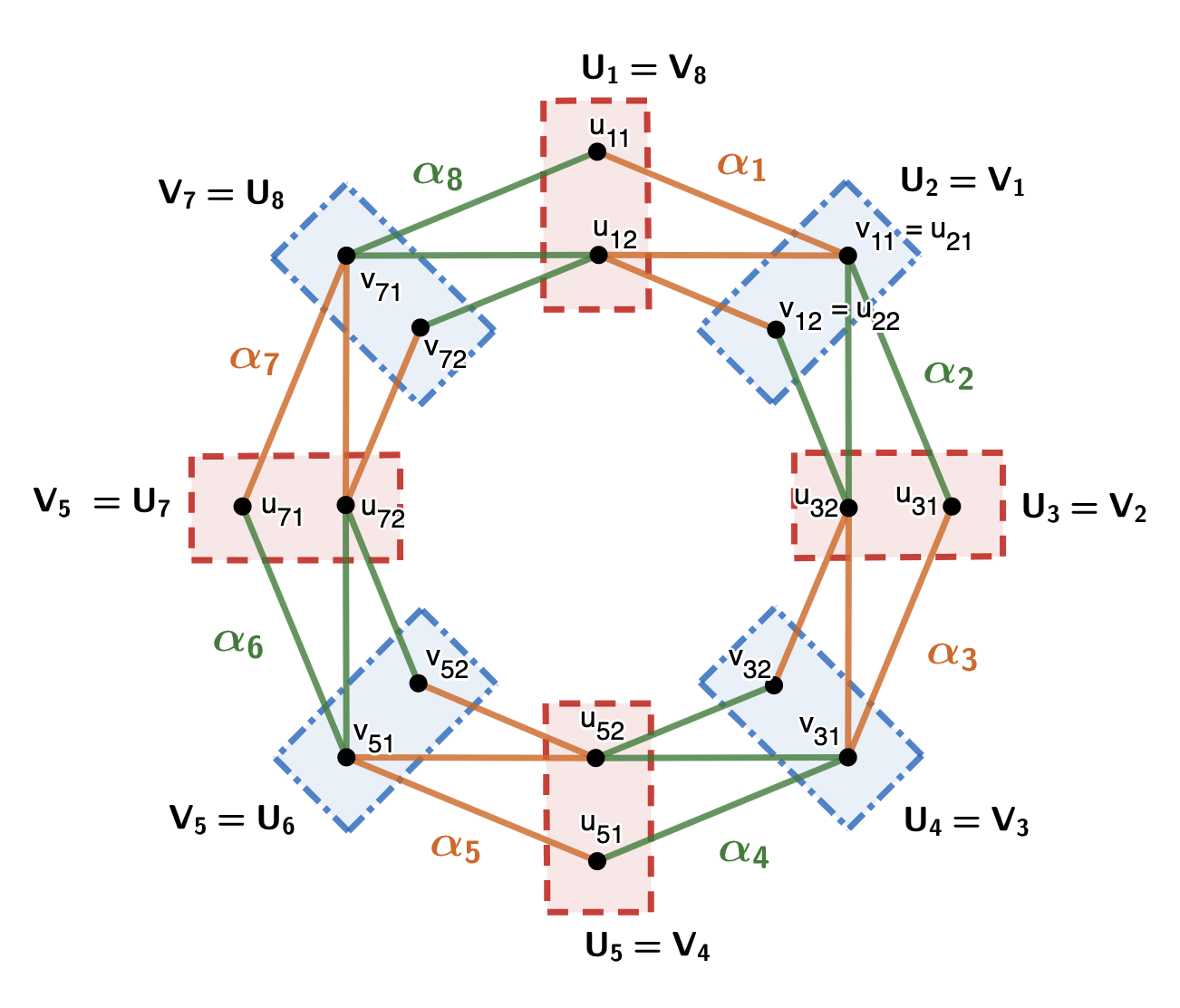}
    \caption{Illustration of \Cref{defn:H-alpha-i}: $v_{ij} = u_{(i+1)j}$ and $V_i = U_{j+1}$ for $j\in[2q]$.}
    \label{fig:H-alpha-i}
\end{figure}

For each $i \in [2q]$, we take $\a_i$ to be the following shape:
\begin{enumerate}
\item If $i$ is odd then we take $U_{\a_i} = \{u_{ij}: j \in [m]\}$, $V_{\a_i} = \{v_{ij}: j \in [m]\}$, and $E(\a_i) = \{\{u_{i(j+1)},v_{ij}\}: j \in [m-1]\} \cup \{\{u_{ij}, v_{ij}\}: j\in [m]\}$. In other words, we take $\a_i$ to be a copy of $\a_{Z_m}$.
\item If $i$ is even then we take $U_{\a_i} = \{u_{ij}: j \in [m]\}$, $V_{\a_i} = \{v_{ij}: j \in [m]\}$, and $E(\a_i) = \{\{u_{ij},v_{i(j+1)}\}: j \in [m-1]\} \cup \{\{u_{ij}, v_{ij}\}: j\in [m]\}$. In other words, we take $\a_i$ to be a copy of $\a_{Z_m}^T$.
\end{enumerate}
See \Cref{fig:H-alpha-i} for an illustration.
\end{defn}

As noted in \Cref{def:copies}, $H(\a_{Z_m},2q)$ can be constructed by concatenating $\a_1,\ldots,\a_{2q}$ and then merging $V_{\a_{2q}}$ and $U_{\a_1}$. More precisely, for each $i \in [2q-1]$ and $j \in [m]$, we set $v_{ij} = u_{(i+1)j}$. We then set $v_{(2q)j} = u_{1j}$ for all $j \in [m]$.

We now consider how $\a_1,\ldots,\a_{2q}$ are affected when we perform the contraction operations.
\begin{prop}\label{prop:R-H-alpha}
For all $H \in R(H(\alpha_{Z_m},2k))$, we can construct $H$ by concatenating $\a'_1,\ldots,\a'_{2k}$ and then merging $V_{\a'_{2k}}$ and $U_{\a'_1}$ where each $\a'_i$ is obtained from $\a_i$ by taking the following steps:
\begin{enumerate}
\item For each $j \in [m]$ such that $u_{ij}$ was deleted but $v_{ij}$ was not deleted, we delete the edges incident to $u_{ij}$ and replace $u_{ij}$ with $v_{ij}$ in $U_{\a_i}$. Note that after this replacement, $v_{ij} \in U_{\a_i} \cap V_{\a_i}$.
\item For each $j \in [m]$ such that $v_{ij}$ was deleted but $u_{ij}$ was not deleted, we delete the edges incident to $v_{ij}$ and replace $v_{ij}$ with $u_{ij}$ in $V_{\a_i}$. Note that after this replacement, $u_{ij} \in U_{\a_i} \cap V_{\a_i}$.
\item For each $j \in [m]$ such that $u_{ij}$ and $v_{ij}$ were both deleted, we delete all edges incident to $u_{ij}$ and $v_{ij}$ and replace $u_{ij}$ and $v_{ij}$ with $u_{i'j}$ in $U_{\a_i}$ and $V_{\a_i}$ where $i' \in [2k]$ is the largest index such that $i' < i$ and $u_{i'j}$ was not deleted. If this $i'$ does not exist because $u_{1j},\ldots,u_{(i-1)j}$ were all deleted then we instead take $i' \in [2k]$ to be the largest index such that $u_{i'j}$ was not deleted.

Note that after this replacement, $u_{i'j} \in U_{\a_i} \cap V_{\a_i}$. Also note that $u_{i'j} = v_{i''j}$ where $i'' \in [2k]$ is the smallest index such that $i'' > i$ and $v_{i''j}$ was not deleted. If this $i''$ does not exist because $v_{(i+1)j},\ldots,u_{(2k)j}$ were all deleted then we instead take $i'' \in [2k]$ to be the smallest index such that $v_{i'j}$ was not deleted.

\begin{figure}[hbt!]
    \centering
    \begin{subfigure}[c]{.67\textwidth}
        \includegraphics[width=1\linewidth]{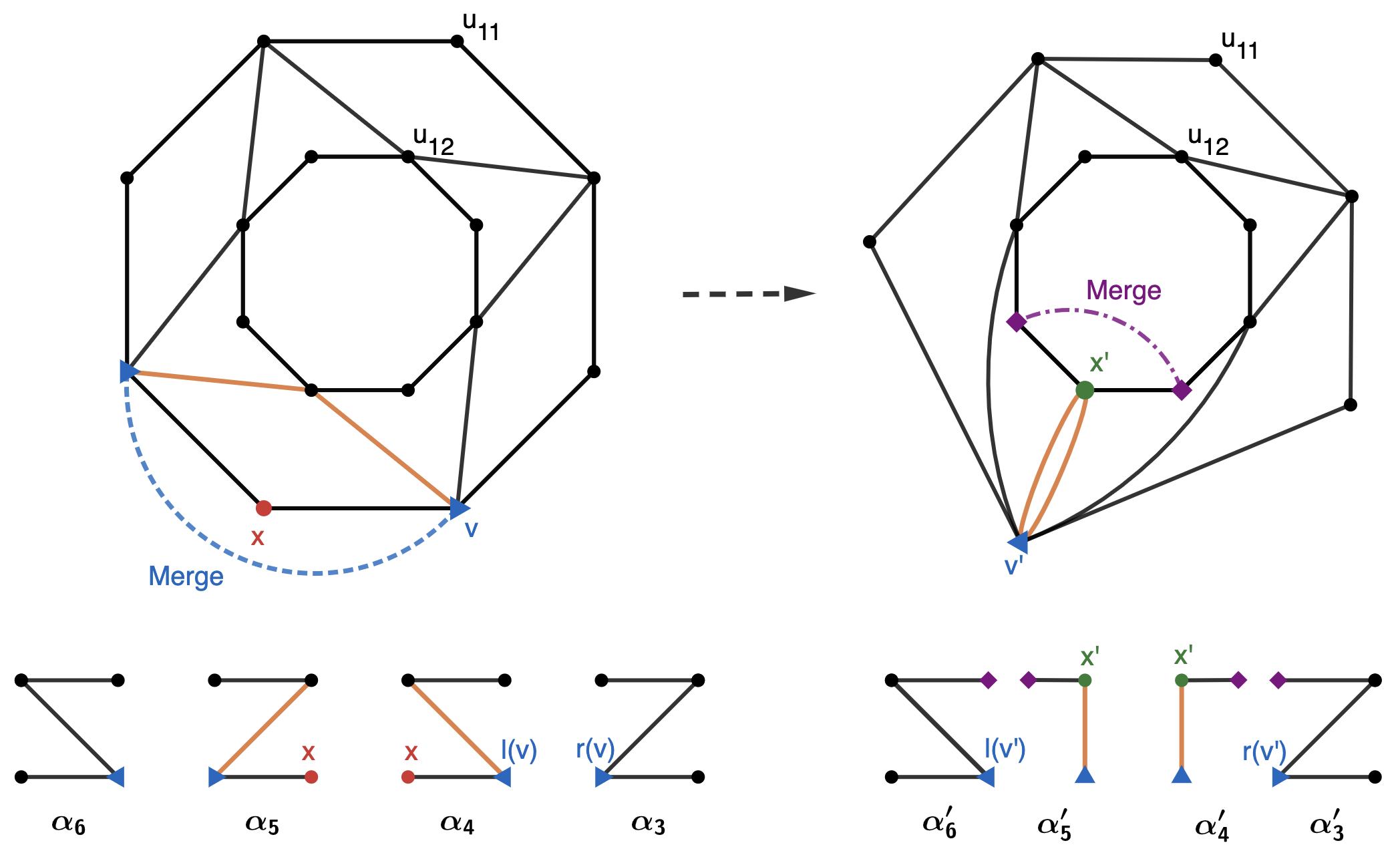}
    \end{subfigure}
    \begin{subfigure}[c]{.32\textwidth}
        \includegraphics[width=1\linewidth]{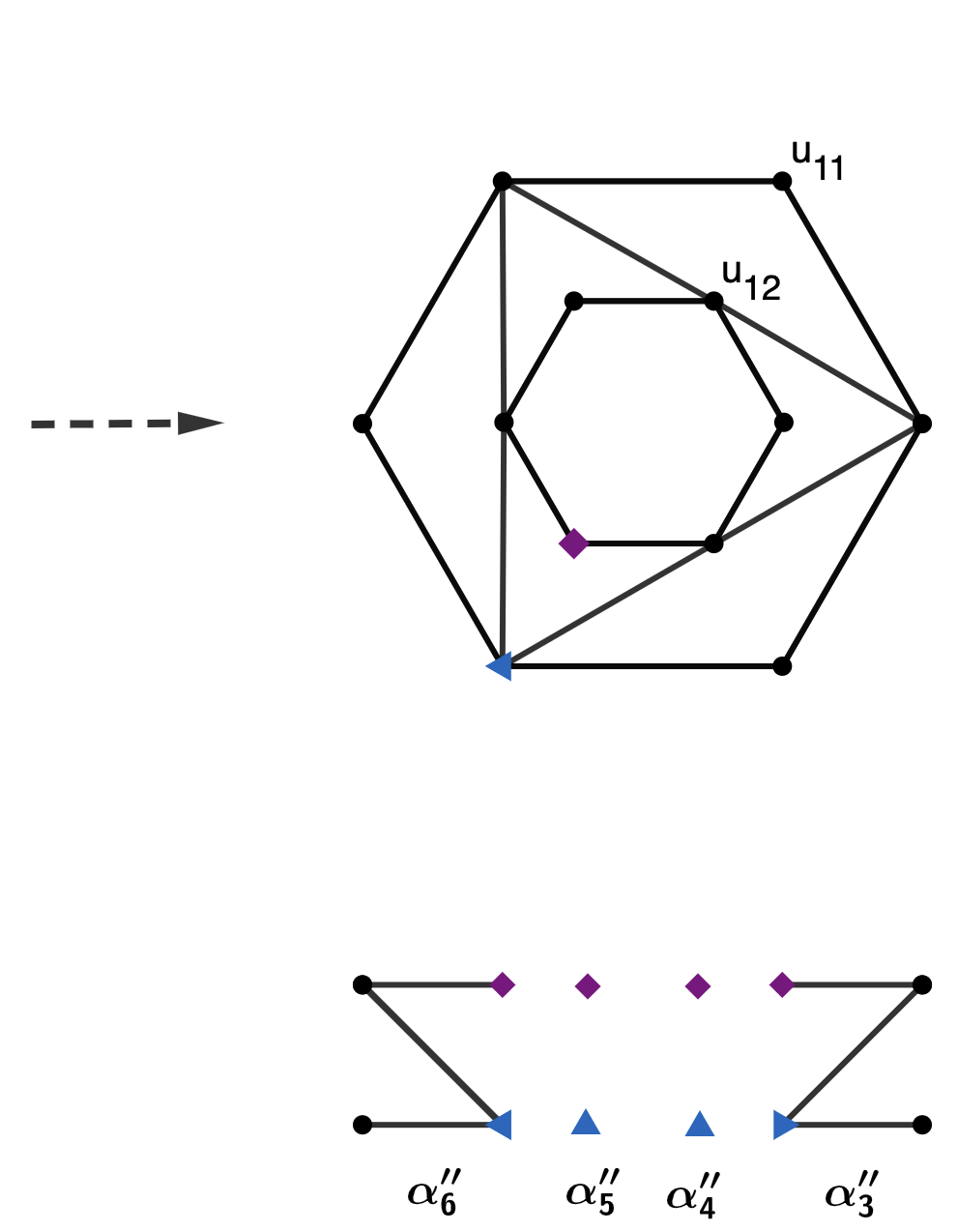}
    \end{subfigure}
    \caption{Illustration of \Cref{prop:R-H-alpha} and \Cref{defn:range-v}. For $v$ in the first figure, $l(v) = 4$ and $r(v) = 3$; for $v'$ in the second figure, $l(v') = 6$ and $r(v') = 3$. }
    \label{fig:R-H-alpha}
\end{figure}
\end{enumerate}
\end{prop}

\begin{defn}\label{defn:range-v}
Given a vertex $v \in V(H)$, we define $l(v) \in [2k]$ to be the index such that $v \in V_{\a'_{l(v)}} \setminus U_{\a'_{l(v)}}$. Similarly, we define $r(v) \in [2k]$ to be the index such that $v \in U_{\a'_{r(v)}} \setminus V_{\a'_{r(v)}}$.

If $l(v) < r(v)$ then we take the range of $v$ to be $range(v) = [l(v),r(v)]$. If $l(v) > r(v)$ then we take the range of $v$ to be $range(v) = [l(v),2k] \cup [1,r(v)]$.
\end{defn}

\begin{prop}
For all $H \in R(H(\alpha_{Z_m},2k))$, for all $j \in [m]$ and all vertices $v$ in the $j$th wheel of $H$,
\begin{enumerate}
\item $range(v)$ is well-defined. For all $i \notin range(v)$, $v \notin V(\a'_i)$.
\item $v = v_{l(v)j} = u_{r(v)j}$ and $r(v) - l(v) \equiv 1 \mod 2$.
\item $u_{l(v)j} \neq v$ and there is an edge between $u_{l(v)j}$ and $v$ in $\a'_{l(v)}$. Similarly, $v \neq v_{r(v)j}$ and there is an edge between $v$ and $v_{r(v)j}$ in $\a'_{r(v)}$.
\item For all $i \in range(v) \setminus \{l(v),r(v)\}$, $v \in U_{\a'_i} \cap V_{a'_i}$.
\end{enumerate}
\end{prop}
\begin{proof}
    We prove this proposition by induction on the number of contraction operations. For the base case, observe that for all $i \in [2k-1]$ and all $j \in [m]$, for each vertex $v = v_{ij} = u_{(i+1)j}$, $l(v) = i$, $r(v) = i+1$, $range(v) = [i,i+1]$, and the third and fourth statements hold. Similarly, for all $j \in [m]$, for each vertex $v = v_{2kj} = u_{1j}$, $l(v) = 2k$, $r(v) = 1$, $range(v) = \{2k\} \cup \{1\}$, and the third and fourth statements hold.
    
    For the inductive step, assume that the proposition is true for all $H \in R(H(\alpha_{Z_m},2q))$ such that $H$ was obtained from $H(\alpha_{Z_m},2q)$ using at most $k$ contractions. Now consider an $H' \in R(H(\alpha_{Z_m},2q))$ which is obtained from $H(\alpha_{Z_m},2q)$ using $k+1$ contractions, a $j \in [m]$, and a vertex $v$ in the $j$th wheel of $H'$. Let $H \in R(H(\alpha_{Z_m},2q))$ be the graph we obtain after the first $k$ contractions. By the inductive hypothesis, $range(v)$ is well-defined for $H$ and satisfies all of the given properties. Let $i_1 = l(v)$ and $i_2 = r(v)$ (with respect to $H$). Now consider what happens after the final contraction. There are $3$ cases to consider.
    \begin{enumerate}
    \item Neither $u_{{i_1}j}$ nor $v_{{i_2}j}$ is deleted in the final contraction. In this case, $range(v)$ is unchanged and still satisfies all of the given properties.
    \item $v_{{i_2}j}$ is deleted in the final contraction. In this case, let $i_3 = r(v_{{i_2}j})$ (with respect to $H$), let $w = v_{{i_3}j}$, and let $i_4 = r(w)$ (with respect to $H$).

    By the inductive hypothesis, we have that in $H$, 
    \begin{enumerate}
    \item For all $i \in range(v) \setminus \{i_1, i_2\}$, $v \in U_{\a'_i} \cap V_{\a'_i}$.
    \item For all $i \in range(v_{{i_2}j}) \setminus \{i_2, i_3\}$, $v_{{i_2}j} \in U_{\a'_i} \cap V_{\a'_i}$.
    \item For all $i \in range(w) \setminus \{i_3, i_4\}$, $w \in U_{\a'_i} \cap V_{\a'_i}$.
    \item $u_{{i_1}j} \neq v$ and $w \neq v_{{i_4}j}$.
    \end{enumerate}
    The contraction merges $v$ and $w$ and deletes $v_{{i_2}j}$. This replaces $v_{{i_2}j}$ with $v = w$ in $V_{\a'_{i_2}}$, $U_{\a'_{i_3}}$, and $U_{\a'_{i}} \cap V_{\a'_{i}}$ for all $i \in range(v_{{i_2}j}) \setminus \{i_2, i_3\}$. Thus, after the contraction, $range(v) = [i_1,i_4]$ (or $[i_4,2q] \cup [1,i_1]$ if $i_4 < i_1$) and $v \in U_{\a'_{i}} \cap V_{\a'_{i}}$ for all $i \in range(v) \setminus \{i_1, i_4\}$.

    It is not hard to check the remaining properties. Since $i_2 - i_1 \equiv 1 \mod 2$, $i_3 - i_2 \equiv 1 \mod 2$, and $i_4 - i_3 \equiv 1 \mod 2$, $i_4 - i_1 \equiv 1 \mod 2$. Since $i_1 \neq i_4$, neither $u_{{i_1}j}$ nor $v_{{i_4}j}$ is equal to $v$.
    \item $u_{{i_1}j}$ is deleted in the final contraction. This case can be handled in a similar way as the previous case.
    \end{enumerate}
\end{proof}

We now analyze the spokes of $H$. 
\begin{lemma}\label{lem:spokeproperties}
For all $H \in R(H(\alpha_{Z_m},2k))$, any $j \in [m]$, and any vertex $v$ in the $j$th wheel of $H$,
\begin{enumerate}
\item If $j = m$ and $l(v)$ is odd or $j = 1$ and $l(v)$ is even then $v$ is not incident to any spokes.
\item If $j > 1$ and $l(v)$ is even then all spokes incident to $v$ are between the $(j-1)$th and $j$th wheels of $H$. Moreover, 
\begin{enumerate}
\item[1.] If neither $u_{l(v)(j-1)}$ nor $v_{r(v)(j-1)}$ was deleted and $u_{l(v)(j-1)} \neq v_{r(v)(j-1)}$ then the spokes $\{u_{l(v)(j-1)},v\}$ and $\{v,v_{r(v)(j-1)}\}$ have odd multiplicity and all other spokes incident to $v$ have even multiplicity.
\item[2.] If neither $u_{l(v)(j-1)}$ nor $v_{r(v)(j-1)}$ was deleted and $u_{l(v)(j-1)} = v_{r(v)(j-1)}$ then all spokes incident to $v$ have even multiplicity.
\item[3.] If at least one of $u_{l(v)(j-1)}$ and $v_{r(v)(j-1)}$ was deleted then both $u_{l(v)(j-1)}$ and $v_{r(v)(j-1)}$ were deleted (in fact, $u_{l(v)(j-1)}$ and $v_{r(v)(j-1)}$ were merged before being deleted) and all spokes incident to $v$ have even multiplicity.
\end{enumerate}
\item If $j < m$ and $l(v)$ is odd then all spokes incident to $v$ are between the $j$th and $(j+1)$th wheels of $H$. Moreover, 
\begin{enumerate}
\item[1.] If neither $u_{l(v)(j+1)}$ nor $v_{r(v)(j+1)}$ was deleted and $u_{l(v)(j+1)} \neq v_{r(v)(j+1)}$ then the spokes $\{u_{l(v)(j+1)},v\}$ and $\{v,v_{r(v)(j+1)}\}$ have odd multiplicity and all other spokes incident to $v$ have even multiplicity.
\item[2.] If neither $u_{l(v)(j+1)}$ nor $v_{r(v)(j+1)}$ was deleted and $u_{l(v)(j+1)} = v_{r(v)(j+1)}$ then all spokes incident to $v$ have even multiplicity.
\item[3.] If at least one of $u_{l(v)(j+1)}$ and $v_{r(v)(j+1)}$ was deleted then both $u_{l(v)(j+1)}$ and $v_{r(v)(j+1)}$ were deleted (in fact, $u_{l(v)(j+1)}$ and $v_{r(v)(j+1)}$ were merged before being deleted) and all spokes incident to $v$ have even multiplicity.
\end{enumerate}
\end{enumerate}
\end{lemma}
\begin{lemma}\label{lem:equalityimplications}
For all $H \in R(H(\alpha_{Z_m},2k))$, any $j \in [m-1]$, and any $i,i' \in [2k]$ such that $i$ is even and $i'$ is odd, if $v_{ij} = u_{i'j}$ then either $v_{i(j+1)} = u_{i'(j+1)}$ or $v_{i(j+1)}$ and $u_{i'(j+1)}$ were merged and then deleted.

Similarly, for all $H \in R(H(\alpha_{Z_m},2k))$, any $j \in [2,m]$, and any $i,i' \in [2k]$ such that $i$ is odd and $i'$ is even, if $v_{ij} = u_{i'j}$ then either $v_{i(j-1)} = u_{i'(j-1)}$ or $v_{i(j-1)}$ and $u_{i'(j-1)}$ were merged and then deleted.
\end{lemma}
\begin{proof}[Proof of Lemmas \ref{lem:spokeproperties} and \ref{lem:equalityimplications}] 
    We prove these two lemmas by induction on the number of contraction operations. The base case $H(\alpha_{Z_m},2k)$ is trivial. For the inductive step, assume that Lemmas \ref{lem:spokeproperties} and \ref{lem:equalityimplications} are true for all $H \in R(H(\alpha_{Z_m},2k))$ such that $H$ was obtained from $H(\alpha_{Z_m},2k)$ using at most $x$ contractions. Now consider an $H' \in R(H(\alpha_{Z_m},2k))$ which is obtained from $H(\alpha_{Z_m},2k)$ using $x+1$ contractions.
    
    There are three ways that Lemmas \ref{lem:spokeproperties} and \ref{lem:equalityimplications} could potentially fail:
    \begin{enumerate}
    \item For some vertex $v \in V(H)$, the other endpoint of one of the spokes incident to $v$ is merged with another vertex (which may be the other endpoint of a different spoke incident to $v$) and the spokes incident to $v$ no longer satisfy the statements in Lemma \ref{lem:spokeproperties}.
    \item Some vertex $v \in V(H)$ is merged with another vertex $w$ and the spokes incident to $v = w$ no longer satisfy the statements in Lemma \ref{lem:spokeproperties}.
    \item Some vertex $v$ in the $j$th wheel of $H$ is merged with another vertex $w$ but the corresponding vertices in the $(j-1)$th or $(j+1)$th wheel of $H$ are not equal/deleted.
    \end{enumerate}
    Let $H \in R(H(\alpha_{Z_m},2q))$ be the graph we obtain after the first $x$ contractions. We show that these three possibilities do not happen. Given $j \in [m]$ and a vertex $v$ in the $j$th wheel of $H$, let $i_1 = l(v)$, and let $i_2 = r(v)$. There are several cases to consider.
    \begin{enumerate}
        \item $i_1 = l(v)$ is odd, $j < m$, and the final contraction deletes a vertex $v'$ in the $(j+1)$th wheel of $H$ and merges its two neighbors $u'$ and $w'$. If $v' = u_{i'(j+1)}$ where $i'$ is odd then $u'$ and $w'$ cannot be endpoints of spokes incident to $v$ but $v'$ may be an endpoint of a spoke incident to $v$. We have the following cases:
        \begin{enumerate}
        \item If $v'$ is equal to $u_{i_1(j+1)}$ or $v_{i_2(j+1)}$ then since $v'$ was not incident to any spokes of odd multiplicity, we must have had that $u_{i_1(j+1)} = v_{i_2(j+1)}$ and this spoke had even multiplicity. Thus, after the contraction, both $u_{i_1(j+1)}$ and $v_{i_2(j+1)}$ are deleted (and they were merged before being deleted).
        \item If $v'$ is not equal to $u_{i_1(j+1)}$ or $v_{i_2(j+1)}$ but is an endpoint of a spoke incident to $v$ then deleting $v'$ deletes this spoke. Since this spoke had even multiplicity, the statements in Lemma \ref{lem:spokeproperties} still hold.
        \item If $v'$ is not an endpoint of a spoke incident to $v$ then the multiplicities of the spokes incident to $v$ are unaffected.
        \end{enumerate}
        If $v' = u_{i'(j+1)}$ where $i'$ is odd then $v'$ cannot be an endpoint of a spoke incident to $v$ but $u'$ and/or $w'$ may be an endpoint of a spoke incident to $v$. We have the following cases:
        \begin{enumerate}
            \item If $u' = u_{i_1(j+1)}$ and $w' = v_{i_2(j+1)}$ or then after the contraction, the spoke $u_{i_1(j+1)} = v_{i_2(j+1)}$ and the spoke $\{u_{i_1(j+1)},v\}$ has even multiplicity.
            \item If $u' \notin \{u_{i_1(j+1)},v_{i_2(j+1)}\}$ or $w' \notin \{u_{i_1(j+1)},v_{i_2(j+1)}\}$ then whether or not $u_{i_1(j+1)} = v_{i_2(j+1)}$ and the parities of the multiplcities of the spokes incident to $v$  are not affected by the contraction.
        \end{enumerate}
        \item $i_1 = l(v)$ is even, $j > 1$, and the final contraction deletes a vertex $v'$ in the $(j-1)$th wheel of $H$ and merges its two neighbors $u'$ and $w'$. This case can be handled in a similar way as the previous case.
        \item $i_1 = l(v)$ is odd and the final contraction deletes the vertex $v_{{i_2}j}$. In this case, let $i_3 = r(v_{{i_2}j})$, let $w = v_{{i_3}j}$, and let $i_4 = r(w)$.
        
        If $j = m$ then all of the vertices which were merged to form $v$ and all of the vertices which were merged to form $w$ were of the form $v_{im}$ where $i$ is odd. Since none of these vertices were incident to any spokes, after the contraction, $v = w$ is not incident to any spokes. 
    
        If $j < m$ then by the inductive hypothesis, all spokes incident to $v$ and $w$ are between the $j$th and $(j+1)$th wheels of $H$. Moreover, before the contraction,
        \begin{enumerate}
            \item[1.] All spokes incident to $v$ except for $\{u_{{i_1}(j+1)},v\}$ and $\{v,v_{{i_2}(j+1)}\}$ have even multiplicity. Similarly, all spokes incident to $w$ except for $\{u_{{i_3}(j+1)},w\}$ and $\{w,v_{{i_4}(j+1)}\}$ have even multiplicity.
            \item[2.] Since $v_{{i_2}j} = u_{{i_3}j}$ and $i_2$ is even, either $v_{{i_2}(j+1)} = u_{{i_3}(j+1)}$ or $v_{{i_2}(j+1)}$ and $u_{{i_3}(j+1)}$ were merged and then deleted.
        \end{enumerate}
        There are several possible cases.
        \begin{enumerate}
            \item If $u_{{i_1}(j+1)}$, $v_{{i_2}(j+1)}$, $u_{{i_3}(j+1)}$, or $v_{{i_4}(j+1)}$ was deleted then $u_{{i_1}(j+1)}$, $v_{{i_2}(j+1)}$, $u_{{i_3}(j+1)}$, and $v_{{i_4}(j+1)}$ must have been all merged together and then deleted.
            \item If $u_{{i_1}(j+1)} = v_{{i_2}(j+1)}$ and $u_{{i_3}(j+1)} = v_{{i_4}(j+1)}$ then $u_{{i_1}(j+1)} = v_{{i_2}(j+1)} = u_{{i_3}(j+1)} = v_{{i_4}(j+1)}$ and all spokes incident to $v$ have even multiplicity after the contraction.
            \item If $u_{{i_1}(j+1)} = v_{{i_2}(j+1)}$ and $u_{{i_3}(j+1)} = v_{{i_4}(j+1)}$ then $u_{{i_1}(j+1)} = v_{{i_2}(j+1)} = u_{{i_3}(j+1)} \neq v_{{i_4}(j+1)}$. In this case, the spokes $\{u_{{i_1}(j+1)},v\}$ and $\{w,v_{{i_4}(j+1)}\}$ are distinct and have odd multiplicity.
            \item If $u_{{i_1}(j+1)} \neq v_{{i_2}(j+1)}$ and $u_{{i_3}(j+1)} = v_{{i_4}(j+1)}$ then $u_{{i_1}(j+1)} \neq v_{{i_2}(j+1)} = u_{{i_3}(j+1)} = v_{{i_4}(j+1)}$. In this case, the spokes $\{u_{{i_1}(j+1)},v\}$ and $\{w,v_{{i_4}(j+1)}\}$ are distinct and have odd multiplicity.
            \item If $u_{{i_1}(j+1)} \neq v_{{i_2}(j+1)}$ and $u_{{i_3}(j+1)} \neq v_{{i_4}(j+1)}$ then after the contraction, the spoke $\{v,v_{{i_2}(j+1)}\} = \{u_{{i_2}(j+1)},w\}$ has even multiplcity and is distinct from the spokes $\{u_{{i_1}(j+1)},v\}$ and $\{w,v_{{i_4}(j+1)}\}$. In this case, if $u_{{i_1}(j+1)} \neq v_{{i_4}(j+1)}$ then after the contraction, the spokes $\{u_{{i_1}(j+1)},v\}$ and $\{w,v_{{i_4}(j+1)}\}$ are distinct and have odd multiplicity. If $u_{{i_1}(j+1)} \neq v_{{i_4}(j+1)}$ then after the contraction, the spoke $\{u_{{i_1}(j+1)},v\} = \{w,v_{{i_4}(j+1)}\}$ has even multiplicity.

            Finally, we observe that for all $j > 1$ (including $j = m$), either $v_{i_1(j-1)} = u_{{i_4}(j-1)}$ or $v_{i_1(j-1)}$ and $u_{{i_4}(j-1)}$ were merged and then deleted. To see this, note that since $i_1$ is odd and $v_{{i_1}j} = u_{{i_2}j}$, by the inductive hypothesis, either $v_{{i_1}(j-1)} = u_{{i_2}(j-1)}$ or $v_{{i_1}(j-1)}$ and $u_{{i_2}(j-1)}$ were merged and then deleted. Similarly, since $i_3$ is odd and $v_{{i_3}j} = u_{{i_4}j}$, either $v_{{i_3}(j-1)} = u_{{i_4}(j-1)}$ or $v_{{i_3}(j-1)}$ and $u_{{i_4}(j-1)}$ were merged and then deleted. Finally, since $i_2$ is even, and $v_{{i_2}j} = u_{{i_3}j}$ was deleted, we must have that either $u_{{i_2}(j-1)} = v_{{i_3}(j-1)}$ or $u_{{i_2}(j-1)}$ and $v_{{i_3}(j-1)}$ were merged and then deleted.

            Putting these observations together, either $v_{i_1(j-1)} = u_{{i_4}(j-1)}$ or $v_{i_1(j-1)}$ and $u_{{i_4}(j-1)}$ were merged and then deleted, as needed.
        \end{enumerate}
        \item The remaining cases are as follows:
        \begin{enumerate}
            \item $i_1 = l(v)$ is odd and the final contraction deletes the vertex $u_{{i_1}j}$.
            \item $i_2 = l(v)$ is even and the final contraction deletes the vertex $v_{{i_2}j}$.
            \item $i_1 = l(v)$ is even and the final contraction deletes the vertex $u_{{i_1}j}$.
        \end{enumerate} 
        These cases can be handled in a similar way as the previous case.
    \end{enumerate}
\end{proof}

\subsection{Proof of Theorems \ref{thm:dominantwellbehaved} and \ref{thm:variancewellbehaved}}
We now prove Theorems \ref{thm:dominantwellbehaved} and \ref{thm:variancewellbehaved}. To do this, we use the decomposition of $H  \in R(H(\a_{Z_m},2k))$ as the composition of $\a'_1,\ldots,\a'_{2k}$ and we use ideas from Appendix B of \cite{AMP20}.

We first define the set of allowable constraint graphs on a multi-graph $H \in R(H(\a_{Z_m},2k))$. These constraint graphs are the possible constraint graphs which could result from starting with a constraint graph $C \in \mathcal{C}_{\(\a_Z,2k\)}$ and applying the contraction operations described in Definition \ref{defn:R-H}.
\begin{defn}
Given an $H \in R(H(\a_{Z_m},2k))$, we define $\mathcal{C}_H$ to be the set of constraint graphs on $H$ which keep all of the vertices in each $\a'_i$ distinct. More precisely, we require that for all $i \in [2k]$ and all $j,j' \in [m]$,
\begin{enumerate}
\item If $j' \neq j$ and $u_{ij},u_{ij'}$ were not deleted (being merged with other vertices is okay) then there is no path of constraint edges between $u_{ij}$ and $u_{ij'}$.
\item If $j' \neq j$ and $v_{ij},v_{ij'}$ were not deleted (being merged with other vertices is okay) then there is no path of constraint edges between $v_{ij}$ and $v_{ij'}$.
\item If $u_{ij},v_{ij'}$ were not deleted (being merged with other vertices is okay) then there is no path of constraint edges between $u_{ij}$ and $v_{ij'}$.
\end{enumerate}
\end{defn}
\begin{prop}\label{lem:edgedistinctness}
For all $H \in R(H(\a_{Z_m},2k))$ and $i \in [2k]$, $E(\a'_i) \subseteq E(\a_i)$. Moreover, for all $C \in \mathcal{C}_H$, even after all of constraint edges are contracted, all of the edges in $E(\a'_i)$ are distinct.
\end{prop}
\begin{proof}
Observe that the only way a new vertex $u_{i'j}$ can be introduced to $V(\a'_i)$ is if $u_{ij}$ and $v_{ij}$ were both deleted first. If so, all edges incident to $u_{ij}$ and $v_{ij}$ were deleted and the new vertex $u_{i'j}$ is not incident to any edges in $E(\a'_i)$. The moreover statement follows from the definition of $\mathcal{C}_H$.
\end{proof}
We now prove Theorem \ref{thm:dominantwellbehaved} by proving the following more general theorem.
\begin{thm}\label{thm:generaldominantwellbehaved}
Given $H \in R(H(\a_{Z_m},2k))$, if $C \in \mathcal{C}_H$ and $\val(C) \neq 0$ then $|E(C)| \geq \frac{|V(H)|}{2} - m$. Moreoever, if $C$ is not well-behaved then $|E(C)| \geq \frac{|V(H)|}{2} - m + 1$
\end{thm}
\begin{proof}
For each $i \in [2k]$, we define the following vertex separator $S_i(C)$ of $\a'_i$. These separators are similar to the  separators used in Appendix B of \cite{AMP20}.
\begin{defn}
Given an $H \in R(H(\a_{Z_m},2k))$, we define $B = \{v \in V(H): l(v) > r(v)\}$ to be the set of vertices in $H$ which appear in $V_{\a'_{2q}} = U_{\a'_1}$.
\end{defn}
\begin{defn}
Given an $H \in R(H(\alpha_{Z_m},2k))$ and a $C \in \mathcal{C}_H$, for each $i \in [2k]$, we define the set $S_i(C)$ to be the set of vertices in $V(\a'_i)$ which are in $B$ or appear both earlier and later in $H$. More precisely, for each $v \in V(\a'_i)$, $v \in S_i(C)$ if and only if at least one of the following holds:
\begin{enumerate}
\item $v \in U_{\a'_i} \cap V_{\a'_i}$
\item $v \in B$ or there is a path of constraint edges between $v$ and a vertex $w \in B$.
\item $l(v) = i$ and there is a path of constraint edges between $v$ and a vertex $u$ such that $l(u) < i$.
\item $r(v) = i$ and there is a path of constraint edges between $v$ and a vertex $w$ such that $r(w) > i$.
\end{enumerate}
\end{defn}
\begin{prop}
If $\val(C) \neq 0$ then for all $i \in [2q]$, $S_i(C)$ is a vertex separator of $\alpha'_{i}$.
\end{prop}
\begin{proof}
Assume that there is an edge $e = \{u,v\} \in E(\a'_i)$ such that $u,v \notin S_i(C)$. We must have that $u = u_{ij}$ and $v = v_{ij'}$ for some $i \in [2q]$ and $j,j' \in [m]$ such that $|j' - j| \leq 1$. Since $u = u_{ij} \notin S_i(C)$, $r(u) = i$, $l(u) < r(u)$, and there is no path of constraint edges between $u$ and any vertex $u'$ such that $r(u') > r(u)$ or $l(u') > r(u')$. Thus, for all $i' > i$, $u \notin V(\a'_{i'})$. Similarly, since $v = v_{ij'} \notin S_i(C)$, $l(v) = i$, $l(v) < r(v)$, and there is no path of constraint edges between $v_{ij'}$ and any vertex $u'$ such that $l(u') < l(v)$ or $l(u') > r(u')$. Thus, for all $i' < i$, $v \notin V(\a'_{i'})$. This implies that $e \notin E(\a'_{i'})$ for any $i' \neq i$ so $e$ only appears once in $H/C$ and thus $\val(C) = 0$, which is a contradiction. 
\end{proof}
These vertex separators are useful because they allow us to upper bound $|V(H/C)|$ which in turn allows us to lower bound $|E(C)|$.
\begin{prop}
$|V(H/C)| = |V(H)| - |E(C)| = \frac{\sum_{i = 1}^{2k}{|V(\a'_{i}) \setminus S_i(C)|}}{2} + |B/C|$ where $|B/C|$ is the number of vertices in $B$ after we contract the constrint edges in $C$.
\end{prop}
\begin{proof}
Observe that in the sum $\sum_{i = 1}^{2k}{|V(\a'_{i}) \setminus S_i(C)|}$, each vertex $v \in V(H/C) \setminus B$ is counted twice, once for the first $i$ such that $v \in V(\a'_i)$ and once for the last $i$ such that $v \in V(\a'_i)$.
\end{proof}
\begin{defn}
Given an $H \in R(H(\alpha_{Z_m},2k))$ and a constraint graph $C$ on $H$, we say that an edge $\{u_{ij},v_{ij}\} \in E(\a'_i)$ is a critical edge if $u_{ij},v_{ij} \in S_i(C)$.
\end{defn}
\begin{cor}\label{cor:constraintedgeslowerbound}
If $H \in R(H(\alpha_{Z_m},2k))$ and $C \in \mathcal{C}_H$ is a constraint graph such that $\val(C) \neq 0$ then $|E(C)| \geq \frac{|V(H)|}{2} - |B/C| + \frac{\# \text{ of critical edges}}{2}$.
\end{cor}
\begin{proof}
Observe that 
\begin{equation*}
|E(C)| = |V(H)| - \frac{\sum_{i = 1}^{2q}{|V(\a'_{i}) \setminus S_i(C)|}}{2} - |B/C|
\end{equation*}
and 
\begin{equation*}
\sum_{i = 1}^{2q}{|S_i(C)|} = 2qm + \# \text{ of critical edges}
\end{equation*}
\end{proof}
To prove Theorem \ref{thm:generaldominantwellbehaved}, it is sufficient to show that if $H$ does not contain a vertex which is not incident to any constraint edges or spokes of odd multiplicity then either $|B/C| < m$ or there are at least two critical edges. 
\begin{defn}
Given $H \in R(H(\alpha_{Z_m},2k))$ and $C \in \mathcal{C}_H$, we say that a vertex $v \in V(H)$ is removable if it is contained in a wheel of $H$ which has at least $4$ vertices and it is not incident to any constraint edges in $C$ or spokes of odd multiplicity in $H$.
\end{defn}
\begin{lemma}\label{lem:findingkey}
For all $H \in R(H(\alpha_{Z_m},2k))$ and $C \in \mathcal{C}_H$, for all $i \in [2k]$ and $j \in [m]$, if $l(v_{ij}) = i$, $v_{ij} \in S_i(C)$, and the $j$th wheel of $H$ has at least $4$ vertices then at least one of the following is true:
\begin{enumerate}
\item There exists an $i' \in [2,i]$ and a $j' \in [m]$ such that $r(u_{i'j'}) = i'$, $l(u_{i'j'}) < i'$, and $u_{i'j'}$ is a removable vertex. Note that $u_{i'j'} \notin B$.
\item There exists an $i' \in [i]$ and a $j' \in [m]$ such that $\{u_{i'j'},v_{i'j'}\}$ is a critical edge.
\end{enumerate}
Similarly, for all $H \in R(H(\alpha_{Z_m},2k))$ and $C \in \mathcal{C}_H$, for all $i \in [2k]$ and $j \in [m]$, if $r(u_{ij}) = i$, $u_{ij} \in S_i(C)$, and the $j$th wheel of $H$ has at least $4$ vertices then at least one of the following is true:
\begin{enumerate}
\item There exists an $i' \in [i,2k-1]$ and a $j' \in [m]$ such that $l(v_{i'j'}) = i'$, $r(v_{i'j'}) > i'$, and $v_{i'j'}$ is a removable vertex. Note that $v_{i'j'} \notin B$.
\item There exists an $i' \in [i,2k]$ and a $j' \in [m]$ such that $\{u_{i'j'},v_{i'j'}\}$ is a critical edge.
\end{enumerate}
\end{lemma}
\begin{proof}
We only prove the first statement as the proof of the second statement is similar.

We prove this statement by induction. For the base case $i = 1$, observe that we automatically have that $u_{1j} \in S_1(C)$. Since $l(v_{1j}) = 1$ and $v_{1j} \in S_1(C)$, $\{u_{1j},v_{1j}\}$ is a critical edge. For the inductive step, assume that the statement is true for $i \leq x$ and assume that $i = x+1$. Let $u = u_{ij}$ and let $i' = l(u)$. There are a few cases to consider:
\begin{enumerate}
\item If $u = u_{ij} \in S_i(C)$ then $\{u_{ij},v_{ij}\}$ is a critical edge.
\item If $u = v_{i'j} \in S_{i'}(C)$ then the result follows from the inductive hypothesis (we know that $i' < i$ as otherwise we would have that $u \in S_i(C)$).
\item If $u$ is not incident to any constraint edges and is not incident to any spokes with odd multiplicity then we are done as $u$ is the vertex we are looking for.
\item If $u$ is not incident to any constraint edges, $u$ is incident to spokes with odd multiplicity, and $i'$ is odd (which happens if and only if $i$ is even) then by Lemma \ref{lem:spokeproperties}, the only spokes incident to $u$ with odd multiplicity are $\{u_{i'(j+1)},v_{i'j}\}$ and $\{u_{ij},v_{i(j+1)}\}$. Since $u = v_{i'j} = u_{ij}$ is not incident to any constraint edges, this means that the vertices $u = u_{i'j}$, $u_{i'(j+1)}$, $v = v_{ij}$, and $v_{i(j+1)}$ are paired up by paths of constraint edges.

Let $w = v_{i(j+1)}$ and let $i'' = l(w)$. Note that there cannot be a path of constraint edges between $v = v_{ij}$ and $w$ as $v$ and $w$ were both part of $\a_i$ so they must be distinct. We must have that $i' \leq l(w) \leq i$ as the spokes $\{u_{i'(j+1)},u\}$ and $\{u,w\}$ are distinct. Thus, regardless of whether there is a path of constraint edges between $w$ and $u_{i'j}$ or a path of constraint edges between $w$ and $u_{i'(j+1)}$, $w \in S_{i''}(C)$.

If $i'' < i$ then the result follows from the inductive hypothesis (note that wheel $j+1$ must have at least $4$ vertices as otherwise all spokes incident to wheel $j+1$ would have even multiplicity). If $i'' = i$ then we can repeat the entire argument using $w$ as the starting vertex (this must eventually terminate as $j$ increases each time the argument is repeated).
\item If $u$ is not incident to any constraint edges, $u$ is incident to spokes with odd multiplcity, and $i'$ is even then we can use a similar argument as the previous case.
\item If $u$ is incident to a constraint edge but $u = v_{i'j} \notin S_{i'}(C)$ and $u = u_{ij} \notin S_i(C)$ then the other endpoint of the constraint edge incident to $u$ must be a vertex $w$ such that $l(w) > l(u)$ and $r(w) < r(u) = i$. Moreover, $w \in S_{l(w)}(C)$ and $w \in S_{r(w)}(C)$. Letting $i'' = l(w)$ and $j'$ be the wheel containing $w$, the result follows by applying the inductive hypothesis on $w$ (unless wheel $j'$ only has two vertices in which case $\{u_{i''j'},v_{i''j'}\}$ and $\{u_{r(w)j'},v_{r(w)j'}\}$ are both critical edges).
\end{enumerate}
\end{proof}
\begin{cor}\label{cor:findingkeycorollary}
For all $H \in R(H(\alpha_{Z_m},2q))$  such that at least one wheel of $H$ has at least $4$ vertices and all $C \in \mathcal{C}_{H}$, there are a total of at least two of the following:
\begin{enumerate}
\item Removable vertices.
\item Critical edges.
\item Paths of constraint edges between vertices of $B$.
\end{enumerate}
\end{cor}
\begin{proof}
We first observe that there is always a critical edge or a removable vertex which is not in $B$. To see this, choose a $j$ such that $W_j$ has at least $4$ vertices and let $u$ be the vertex in $W_j \cap B$. Letting $i = l(u)$, observe that $u = v_{ij} \in S_i(C)$. By Lemma \ref{lem:findingkey}, we can find a critical edge or a remoavable vertex which is not in $B$, as needed.

This implies that if there is a path of constraint edges between two vertices of $B$ then the result is true. Thus, we can assume that $C$ does not have a path of constraint edges between any two vertices of $B$.

We now observe that if $C$ contains a path of constraint edges between a vertex $u \in B$ and a vertex $v \notin B$ then there must be at a total of at least two critical edges and/or removable vertices. To see this, let $W_j$ be the wheel containing $v$, let $i_1 = l(v)$ and let $i_2 = r(v)$. Observe that $i_1 < i_2$. If $W_j$ only contains two vertices then both $\{u_{{i_1}j},v_{{i_1}j}\}$ and $\{u_{{i_2}j},v_{{i_2}j}\}$ are critical edges. If $W_j$ contains at least 4 vertices then by Lemma \ref{lem:findingkey}, either there exists an $i' \in [2,i_1]$ and a $j' \in [m]$ such that $u_{i'j'}$ is a removable vertex or there exists an $i' \in [i_1]$ and a $j' \in [m]$ such that $\{u_{i'j'},v_{i'j'}\}$ is a critical edge (or both). Similarly, by Lemma \ref{lem:findingkey}, either there exists an $i'' \in [i_2,2k-1]$ and a $j'' \in [m]$ such that $v_{i''j''}$ is a removable vertex or there exists an $i'' \in [i_2,2k]$ and a $j'' \in [m]$ such that $\{u_{i''j''},v_{i''j''}\}$ is a critical edge (or both).

The remaining case is when none of the vertices in $B$ are incident to a constraint edge. In this case, let $u \in B$ be a vertex such that $u$ is contained in a wheel with at least $4$ vertices and $r(u) + 2k - l(u)$ is minimized. If $u$ is not incident to any spokes with odd multiplicity then $u$ is a removable vertex. If $u$ is incident to spokes with odd multiplicity, let $W_j$ be the wheel containing $u$. We have the following cases:
\begin{enumerate}
\item If $l(u)$ is even and $r(u)$ is odd then by Lemma \ref{lem:spokeproperties}, the spokes of odd multiplicity which $u$ is incident to are $\{u_{l(u)(j-1)},u\}$ and $\{u,v_{r(u)(j-1)}\}$. Since $u$ is not incident to any constraint edges, $u_{l(u)(j-1)}$ and $v_{r(u)(j-1)}$ must both be incident to constraint edges. Letting $u'$ be the vertex in $W_
{j-1}$ which is in $B$, we must have that $l(u) \leq l(u')$ and $r(u') \leq r(u)$. Since we chose $u$ to minimize $r(u) + 2k - l(u)$, we must have that $l(u') = l(u)$ and $r(u') = r(u)$. Thus, we can repeat this argument with $u'$. Since $j$ decreases every time we repeat this argument, we cannot repeat this argument indefinitely so we must eventually reach a removable vertex in $B$. Since we always have either a critical edge or a removable vertex which is not in $B$, the result follows.
\item If $l(u)$ is odd and $r(u)$ is even then we can use a similar argument. By Lemma \ref{lem:spokeproperties}, the spokes of odd multiplicity which $u$ is incident to are $\{u_{l(u)(j+1)},u\}$ and $\{u,v_{r(u)(j+1)}\}$. Since $u$ is not incident to any constraint edges, $u_{l(u)(j+1)}$ and $v_{r(u)(j+1)}$ must both be incident to constraint edges. Letting $u'$ be the vertex in $W_
{j+1}$ which is in $B$, we must have that $l(u) \leq l(u')$ and $r(u') \leq r(u)$. Since we chose $u$ to minimize $r(u) + 2k - l(u)$, we must have that $l(u') = l(u)$ and $r(u') = r(u)$. Thus, we can repeat this argument with $u'$. Since $j$ increases every time we repeat this argument, we cannot repeat this argument indefinitely so we must eventually reach a removable vertex in $B$. Since we always have either a critical edge or a removable vertex which is not in $B$, the result follows.
\end{enumerate}
\end{proof}
we can now prove Theorem \ref{thm:generaldominantwellbehaved} by induction on $\frac{|V(H)|}{2} - m$. The base case $\frac{|V(H)|}{2} - m = 0$ is trivial as $|E(C)| \geq 0$ and if $C$ is not well-behaved then $|E(C)| \geq 1$. For the inductive step, assume that Theorem \ref{thm:generaldominantwellbehaved} is true for all $H \in R(H(\a_{Z_m},2k))$ such that $\frac{|V(H)|}{2} - m = x$. Consider an $H \in R(H(\a_{Z_m},2k))$ such that $\frac{|V(H)|}{2} - m = x+1$. 

If $H$ contains a removable vertex $v$ then there must be a path of constraint edges between the two neighbors of $v$ so we can merge these two vertices together and delete $v$. This gives us a constraint graph $H'$ and a constraint graph $C' \in \mathcal{C}_{H}$ with $|V(H')| = |V(H)| - 2$ and $|E(C')| = |E(C) - 1$. The result now follows from the inductive hypothesis.

If $H$ does not contain a removable vertex $v$ then by Corollary \ref{cor:findingkeycorollary} and Corollary \ref{cor:constraintedgeslowerbound}, $|E(C)| \geq \frac{|V(H)|}{2} - m + 1$, as needed. 
\end{proof}
We now prove Theorem \ref{thm:variancewellbehaved} by proving the following more general theorem.
\begin{thm}\label{thm:generalvariancewellbehaved}
Given $H_1,H_2 \in R(H(\a_{Z_m},2k))$, if $C$ is a constraint graph on $H_1 \cup H_2$ such that $\val(C) \neq 0$ and the restrictions $C_1,C_2$ of $C$ to $H_1$ and $H_2$ are in $\mathcal{C}_{H_1}$ and $\mathcal{C}_{H_2}$ respectively then $|E(C)| \geq \frac{|V(H_1)|}{2} + \frac{|V(H_2)|}{2} - 2m$. Moreover, if $\val(C_1) = 0$ or $\val(C_2) = 0$ then $|E(C)| \geq \frac{|V(H_1)|}{2} + \frac{|V(H_2)|}{2} - 2m + 2$.
\end{thm}
\begin{proof}[Proof sketch]
We prove this theorem by using the same strategy we used to prove \ref{thm:generaldominantwellbehaved}. We adapt this strategy to analyze $H_1 \cup H_2$ as follows
\begin{enumerate}
\item We shift the indices of $H_2$ to be in $[2k+1,4k]$ rather than $[1,2k]$. This places all of $H_2$ to the right of $H_1$.
\item We take $B = \{v \in V(H_1) \cup V(H_2): l(v) > r(v)\}$. Note that we can decompose $B$ as $B = B_1 \cup B_2$ where $B_1 = \{v \in V(H_1): l(v) > r(v)\}$ and $B_2 = \{v \in V(H_2): l(v) > r(v)\}$.
\item We define the vertex separator $S_i(C)$ for $\a'_i$ in the same way as before.
\item The analogue of Corollary {cor:constraintedgeslowerbound} is that if $H_1,H_2 \in R(H(\alpha_{Z_m},2k))$ and $C \in \mathcal{C}_H$ is a constraint graph such that $\val(C) \neq 0$ and the restrictions $C_1,C_2$ of $C$ to $H_1$ and $H_2$ are in $\mathcal{C}_{H_1}$ and $\mathcal{C}_{H_2}$ respectively then 
\begin{equation*}
|E(C)| \geq \frac{|V(H_1)| + |V(H_2)|}{2} - |B/C| + \frac{\# \text{ of critical edges}}{2}.    
\end{equation*}
\item Lemma \ref{lem:findingkey} holds for both $H_1$ and $H_2$ (where the $i$ indices are increased by $2k$ for $H_2$).
\end{enumerate}
We now show the following analogue of Corollary \ref{cor:findingkeycorollary}.
\begin{cor}\label{cor:findingkeycorollaryanalogue}
Given $H_1,H_2 \in R(H(\a_{Z_m},2k))$ such that both $H_1$ and $H_2$ have at least one wheel which has at least $4$ vertices, if $C$ is a constraint graph on $H_1 \cup H_2$ such that $\val(C) \neq 0$ and the restrictions $C_1,C_2$ of $C$ to $H_1$ and $H_2$ are in $\mathcal{C}_{H_1}$ and $\mathcal{C}_{H_2}$ respectively then at least one of the following must hold:
\begin{enumerate}
\item $H_1$ or $H_2$ contains a removable vertex.
\item $2m - |B/C| + \frac{\# \text{ of critical edges}}{2} \geq 2$.
\end{enumerate}
\end{cor}
\begin{proof}
We use similar reasoning as we used to prove Corollary \ref{cor:findingkeycorollary}. In particular, we make the following observations. If neither $H_1$ nor $H_2$ contains a removable vertex then
\begin{enumerate}
\item By the same reasoning as the first paragraph of the proof of Corollary \ref{cor:findingkeycorollary}, both $H_1$ and $H_2$ contain at least one critical edge.
\item If there is a path of constraint edges from a vertex $u \in B$ to a vertex $v \in V(H_1) \setminus B$ then by the same reasoning as the third paragraph of the proof of Corollary \ref{cor:findingkeycorollary}, $H_1$ must contain at least two critical edges. Since $H_2$ must contain at least one critical edge and the number of critical edges must be even, there must be at least $4$ critical edges.
\item Similarly, if there is a path of constraint edges from a vertex $u \in B$ to a vertex $v \in V(H_2) \setminus B$ then $H_1$ must contain at least two critical edges. Since $H_1$ must contain at least one critical edge and the number of critical edges must be even, there must be at least $4$ critical edges.
\item If there is a path of constraint edges between two vertices $u$ and $v$ in $B$ then since both $H_1$ and $H_2$ must contain at least one critical edge, we have that $2m - |B/C| + \frac{\# \text{ of critical edges}}{2} \geq 2$.
\item If none of the vertices in $B$ are incident to a constraint edge then following the same reasoning we used in the proof of Corollary \ref{cor:findingkeycorollary}, we can find a removable vertex in both $H_1$ and $H_2$. This contradicts the assumption that neither $H_1$ nor $H_2$ has a removable vertex.
\end{enumerate}
\end{proof}
we can now prove Theorem \ref{thm:generaldominantwellbehaved} by induction on $\frac{|V(H_1)|}{2} + \frac{|V(H_2)|}{2} - 2m$. For the base case, observe that if $\frac{|V(H_1)|}{2} + \frac{|V(H_2)|}{2} - 2m = 0$ then the result holds as $|E(C)| \geq \frac{|V(H_1)|}{2} + \frac{|V(H_2)|}{2} - 2m = 0$ and $\val(C_1) = \val(C_2) = 1$ as all edges in $E(H_1)$ and $E(H_2)$ have even multiplicity. For the inductive step, assume the result is true for $H_1,H_2$ such that $\frac{|V(H_1)|}{2} + \frac{|V(H_2)|}{2} - 2m = x$ and consider an $H_1,H_2$ such that $\frac{|V(H_1)|}{2} + \frac{|V(H_2)|}{2} - 2m = x+1$. 

If all wheels of $H_1$ only have two vertices then $\frac{|V(H_1)|}{2} = m$ and all edges of $H_1$ appear with even multiplicity. In this case, $\val(C_1) = \val(C_2) = 1$ as the edges in $H_1$ are already paired up with each other so $C_2$ must pair up the edges of $H_2$. By Theorem \ref{thm:generaldominantwellbehaved}, $|E(C)| \geq E(C_2) \geq \frac{|V(H_2)|}{2} - m = \frac{|V(H_1)|}{2} + \frac{|V(H_2)|}{2} - 2m$. Similar logic applies if all wheels of $H_2$ only have two vertices.

If both $H_1$ and $H_2$ have a wheel which has at least $4$ vertices and either $H_1$ or $H_2$ has a removable vertex, we can delete this vertex, merge the two neighbors of this vertex, and use the inductive hypothesis. If neither $H_1$ nor $H_2$ has a removable vertex then by Corollary \ref{cor:findingkeycorollaryanalogue}, $|E(C)| \geq \frac{|V(H_1)| + |V(H_2)|}{2} - 2m$ and if $\val(C_1) \neq 0$ or $\val(C_2) \neq 0$, $|E(C)| \geq \frac{|V(H_1)| + |V(H_2)|}{2} - 2m + 2$, as needed.
\end{proof}

\end{appendix}

\end{document}